\date{December 16, 2009} 

\documentclass[10pt]{article}

\usepackage{amssymb}
\usepackage{amsmath}
\input{liemacs.sty} 
\usepackage{glossary} 
\makeglossary      

\newcommand{\rO}{\mathop{\rm O{}}\nolimits}
\renewcommand{\mlabel}{\label} 

\begin{document} 



\title{Semibounded Representations\\ and Invariant 
Cones in \break Infinite Dimensional Lie Algebras} 
\author{Karl-Hermann Neeb
\begin{footnote}{
Fachbereich Mathematik, TU Darmstadt, Schlossgartenstrasse 7, 
64289-Darmstadt, Germany; neeb@mathematik.tu-darmstadt.de}
\end{footnote}
\begin{footnote}{Supported by DFG-grant NE 413/7-1, Schwerpunktprogramm 
``Darstellungstheorie''.} 
\end{footnote}
}

\maketitle
\date 
\begin{abstract} A unitary representation 
of a, possibly infinite dimensional, Lie group 
$G$ is called semi-bounded if the corresponding 
operators $i\dd\pi(x)$ from the derived representations 
are uniformly bounded from above on some non-empty open subset 
of the Lie algebra $\g$. In the first part of the present paper  
we explain how this concept leads to a fruitful 
interaction between the areas of infinite dimensional 
convexity, Lie theory, symplectic geometry (momentum maps) 
and complex analysis. Here 
open invariant cones in Lie algebras play a central role and 
semibounded representations have interesting connections to 
$C^*$-algebras which are quite different from the classical 
use of the group $C^*$-algebra of a finite dimensional 
Lie group. The second half is devoted
to a detailed discussion of semibounded representations 
of the diffeomorphism group of the circle, the 
Virasoro group, the metaplectic representation on the 
bosonic Fock space and the spin representation on fermionic 
Fock space. 
\end{abstract}

\section{Introduction} \mlabel{sec:1}

In the unitary representation theory of a finite dimensional 
Lie group $G$ a central tool is the convolution 
algebra $L^1(G)$, resp., its enveloping $C^*$-algebra $C^*(G)$, 
whose construction is based on the Haar measure, 
whose existence follows from the 
local compactness of $G$. Since the non-degenerate representations 
of $C^*(G)$ are in one-to-one correspondence to continuous unitary 
representations of $G$, the full power of the rich theory of $C^*$-algebras 
can be used to study unitary representations of $G$. 
To understand and classify irreducible unitary representations, 
the crucial methods are usually based on the fine structure theory of 
finite dimensional Lie groups, such as Levi and Iwasawa decompositions. 
Both methods are no longer available for infinite dimensional 
Lie groups because they are not locally compact and there is no general 
structure theory available. 

However, there are many interesting classes of infinite dimensional 
Lie groups which possess a rich unitary representation theory. 
Many of these representation show up naturally in various contexts 
of mathematical physics (\cite{Ca83}, 
\cite{Mick87}, \cite{Mick89}, \cite{PS86}, 
\cite{SeG81}, \cite{CR87}, \cite{Se58}, \cite{Se78}, 
\cite{Bak07}). 
The representations arising in mathe\-matical physics, resp., 
Quantum Field Theory are often characterized by the 
requirement that the Lie algebra $\g = \L(G)$ 
of $G$ contains an element $h$, corresponding 
to the Hamiltonian of the underlying physical system, for which the 
spectrum of the operator $i\cdot \dd \pi(h)$ 
in the ``physically relevant'' representations $(\pi, \cH)$ 
is non-negative. These representations are called {\it positive energy 
representations} (cf.\ \cite{Se67}, \cite{Bo96}, \cite{SeG81}, \cite{FH05}). 

To develop a reasonably general powerful theory of unitary representations 
of infinite dimensional Lie groups, new approaches have to be developed 
which do neither rest on a fine structure theory nor on the existence of 
invariant measures. 
In this note we describe a systematic approach which is 
very much inspired by the concepts and requirements of 
mathematical physics and which provides a unifying framework 
for a substantial class of representations and several interesting phenomena. 
Due to the lack of a general structure theory, 
one has to study specific classes 
of representations. Here we focus on 
{\it semibounded representations}. Semiboundedness is a stable version of the 
positive energy condition. It 
means that the selfadjoint operators $i\dd\pi(x)$ from the derived 
representation are uniformly bounded below for all $x$ in some 
non-empty open subset of $\g$. 
Our long term goal is to understand the decomposition 
theory and the irreducible semibounded representations by 
their geometric realizations.  

The theory of semibounded unitary representations 
combines results, concepts and methods from several branches of 
mathematics: the theory of convex sets and functions in locally  
convex spaces, infinite dimensional Lie theory, 
symplectic geometry (momentum maps, coadjoint orbits) 
and complex geometry (infinite dimensional K\"ahler manifolds 
and complex semigroup actions). 
In Sections 2-5 below, we describe the relevant aspects of these four 
areas and recall some basic results from \cite{Ne08, Ne09a}. 
A crucial new point is that our approach 
provides a common functional analytic environment 
for various important classes of unitary representations  
of infinite dimensional Lie groups. 
That it is now possible to study semibounded representations 
in this generality 
is  due to the recent progress in infinite dimensional 
Lie theory with fundamental achievements in the past decade. 
For a detailed survey we refer to \cite{Ne06}. A comprehensive 
exposition of the theory will soon be available in our 
monograph with H.~Gl\"ockner \cite{GN09}. 

For finite dimensional Lie groups semibounded unitary representations 
are well understood. In \cite{Ne00} they are called 
``generalized highest weight representations'' because 
the irreducible ones permit a classification in terms of their highest weight 
with respect to a root decomposition of a suitable quotient 
algebra (see Remark~\ref{rem:findim} for more details on this case). 
This simple picture does not carry over to infinite 
dimensional groups.

We now describe our setting in some more detail. 
Based on the notion of a smooth map between open subsets of a 
locally convex space one obtains the concept of a locally 
convex manifold and hence of a locally convex Lie group 
(cf.\ \cite{Ne06}, \cite{Mil84}, \cite{GN09}). 
In Section~\ref{sec:3} we discuss some key examples. 
Let $G$ be a (locally convex) Lie group and $\g$ be its 
Lie algebra. 
For a unitary representation $(\pi,\cH)$ of $G$ we write 
$\pi^v(g) := \pi(g)v$ for its orbit maps  
and call the representation $(\pi, \cH)$ of $G$ {\it smooth} 
if the space 
$$\cH^\infty := \{ v \in \cH \: \pi^v \in C^\infty(G,\cH)\}$$ 
of smooth vectors is dense in $\cH$. 
Then all operators $i \dd \pi(x)$, $x \in \g$, are 
essentially selfadjoint and crucial  
information on their spectrum is contained in the momentum set $I_\pi$ 
of the representation, which is a weak-$*$-closed convex subset 
of the topological dual $\g'$. 
It is defined as the weak-$*$-closed convex hull of the image of the 
momentum map on the projective space of $\cH^\infty$: 
$$ \Phi_\pi \: \bP({\cal H}^\infty)\to \g' \quad \hbox{ with } \quad 
\Phi_\pi([v])(x) 
= \frac{1}{i}  \frac{\la  \dd\pi(x).v, v \ra}{\la v, v \ra}\quad \mbox{ for } 
[v] = \C v. $$
As a weak-$*$-closed convex subset, $I_\pi$ is completely determined 
by its support functional 
$$ s_\pi \: \g \to \R \cup \{\infty\}, \quad s_\pi(x) 
= - \inf \la I_\pi,x \ra = \sup(\Spec(i\dd \pi(x))). $$
It is now natural to study those representations for which 
$s_\pi$, resp., the set $I_\pi$, contains the most significant 
information, and these are precisely the semibounded ones. 
As we shall see in Remark~\ref{rem:4.7}, the geometry of the sets 
$I_\pi$ is closely connected to invariant cones, so that 
we have to take a closer 
look at infinite dimensional Lie algebras containing 
open invariant convex cones $W$ which are pointed in the sense that they 
do not contain any affine line. 

For finite dimensional Lie algebras, there is a well developed 
structure theory of invariant convex cones 
(\cite{HHL89}) and even a characterization 
of those finite dimensional Lie algebras 
containing pointed invariant cones (\cite{Ne94}, 
\cite{Neu99}; see also \cite{Ne00} for a self-contained 
exposition). As the examples described in Section~\ref{sec:6} show, 
many key features of the finite dimensional theory survive, but 
a systematic theory of open invariant cones remains to be developed. 
A central point of the present note is to exploit properties of 
open invariant cones for the theory of semibounded 
representations, in particular 
to verify that certain unitary representations are 
semibounded. In 
Section~\ref{sec:7} we discuss two aspects of semiboundedness 
in the representation theory of $C^*$-algebras, namely 
the restrictions of algebra representations to the unitary group 
$\U(\cA)$ and covariant representations with respect to a Banach--Lie group 
acting by automorphisms on~$\cA$. 

Section~\ref{sec:8} is devoted to a detailed analysis 
of invariant convex cones in the Lie algebra 
$\cV(\bS^1)$ of smooth vector fields on the circle and its 
central extension, the Virasoro algebra $\vir$. 
In particular we show that, up to sign, there are only two open 
invariant cones in $\cV(\bS^1)$. From this insight we derive that 
the group $\Diff(\bS^1)_+$ of orientation preserving diffeomorphisms 
of $\bS^1$ has no non-trivial semibounded unitary representations, 
which is derived from the triviality of all unitary highest weight 
modules (\cite{GO86}). As one may expect from its importance 
in mathematical physics, the situation is different 
for the Virasoro group $\Vir$. 
For $\Vir$ we prove a convexity theorem for 
adjoint and coadjoint orbits which provides complete 
information on invariant cones and permits us to 
determine the momentum sets of the unitary highest weight representations. 
In particular, we show that these, together with their duals, are precisely 
the irreducible semibounded unitary representations. 
Our determination of the momentum sets uses the complex analytic tools 
from Section~\ref{sec:5}, which lead to a realization 
in spaces of holomorphic 
sections on the complex manifold $\Diff(\bS^1)_+/\bS^1$. 
This manifold has many interesting realizations. 
In string theory it occurs as as a space 
of complex structures  on the based loop space 
$C^\infty_*(\bS^1,\R)$ (\cite{BR87}), 
and Kirillov and Yuriev 
realized it as a space of univalent holomorphic functions 
on the open complex unit disc (\cite{Ki87}, \cite{KY87}). 
Its close relative $\Diff(\bS^1)/\PSL_2(\R)$ can be identified 
with the space of Lorentzian metrics on the one-sheeted hyperboloid 
(cf.\ \cite{KS88}), which leads in particular to an interpretation 
of the Schwarzian derivative in terms of a conformal factor. 

In Sections~\ref{sec:9} and \ref{sec:10} we continue our 
discussion of important examples with the automorphism groups 
$\Sp(\cH)$ of the canonical commutation relations (CCR) and 
$\rO(\cH)$ of the canonical anticommutation relations (CAR). 
Geometrically, $\Sp(\cH)$ is the group of real linear automorphisms 
of a complex Hilbert space $\cH$ preserving the imaginary part 
$\omega(x,y) := \Im\la x,y\ra$ of the scalar product and 
$\rO(\cH)$ is the group of real linear automorphisms preserving its real 
part $\beta(x,y) := \Re \la x,y\ra$ (cf.\ \cite{BR97}). 
Section~\ref{sec:9} is dedicated to the Fock representation of the 
(CCR). Here we start with the unitary representation $(W, S(\cH))$ of the 
Heisenberg group $\Heis(\cH)$ on the symmetric/bosonic Fock space 
$S(\cH)$. The group $\Sp(\cH)$ acts naturally by automorphisms 
$\alpha_g$ on $\Heis(\cH)$ and $W \circ \alpha_g$ is equivalent to 
$W$ if and only if $g$ belongs to the {\it restricted symplectic 
group} $\Sp_{\rm res}(\cH)$, i.e., its antilinear part 
$g_2$ is Hilbert--Schmidt. 
Since  
the Fock representation of $\Heis(\cH)$ is irreducible, this leads to a 
projective representation of $\Sp_{\rm res}(\cH)$ on $S(\cH)$. 
We write $\hat\Sp_{\rm res}(\cH)$ for the corresponding central 
$\T$-extension; the {\it metaplectic group} 
(\cite{Se59}, \cite{Sh62}, \cite[Sect.~5]{SeG81}). 
Using a quite general smoothness 
criterion (Theorem~\ref{thm:unirep-cenext}), we show that this group 
carries a natural Banach--Lie group structure and that its canonical 
unitary representation on $S(\cH)$ (the metaplectic representation) 
is smooth. From an explicit formula for the corresponding 
Lie algebra cocycle, we derive a natural 
presentation of this group as a quotient of a semidirect 
product. We further show that the metaplectic 
representation is semibounded and determine the cone of semibounded 
elements. The combined representation of the semidirect product 
$\Heis(\cH) \rtimes \hat\Sp_{\rm res}(\cH)$ is also semibounded 
and irreducible, and, using the tools from Section~\ref{sec:5}, we show 
that its momentum set is the closed convex hull of a single coadjoint 
orbit. 

In Section~\ref{sec:10} we then turn 
to the (CAR), for which we consider the representation on 
the fermionic Fock space $\Lambda(\cH)$. 
Here we likewise obtain a projective representation of the 
restricted orthogonal group 
$\OO_{\rm res}(\cH)$ (\cite{ShSt65}). With completely analogous 
arguments we then show that the corresponding central extension 
$\hat\OO_{\rm res}(\cH)$ is a Lie group, the {\it metagonal group}, 
whose representation on $\Lambda(\cH)$ (the spin representation) 
is smooth and semibounded. We also prove that the 
momentum set of the even spin representation 
is the weak-$*$-closed convex hull of a single coadjoint orbit. 
Our discussion of the bosonic and fermionic Fock representations 
are very much inspired by the construction of the metaplectic 
representation in \cite{Ve77} and the presentation in 
\cite{Ot95}. For finer results on cocycles and connections to physics 
we refer to \cite{Lm94}, \cite{Ru77} and \cite{SeG81}. 
For a more detailed discussion 
of the metaplectic and the metagonal group, see \cite{Ve90} 
and \cite[Sect.~IV.2]{Ne02a}. 

The example of the orthogonal group $\OO_{\rm res}(\cH)$ 
is particularly instructive because it shows very naturally 
how semibounded representations enter the scene for 
infinite dimensional analogs, such as $\rO_{\rm res}(\cH)$,  
of compact groups, for which 
one rather expects to see bounded representations such as 
the spin representation of $\rO_1(\cH)$ (cf.\ \cite{Ne98}). 
In many cases, such as for $\rO_1(\cH)$, 
the class of groups with bounded representations is too restrictive 
to do justice to the underlying geometry. What makes the larger 
groups more interesting is the rich supply of exterior automorphisms 
and the existence of non-trivial central extensions encoding 
relevant geometric information. 

As our examples show, 
the process of second quantization, i.e., passing from a 
one-particle Hilbert space to a many particle space, 
destroys norm continuity for the representation of the 
automorphism groups. What survives is semiboundedness for 
the centrally extended groups acting on the many particle spaces. 
This is closely related to the fact that 
Lie groups of non-unitary maps on a Hilbert space 
$\cH$, such as $\Sp_{\rm res}(\cH)$ have unitary representations on the 
corresponding many particle spaces. In physics language this means that 
``unphysical symmetries'' of the one particle Hilbert space 
may lead to symmetries of the many particle 
space (cf. \cite{Ru77}). In the context of finite dimensional 
Lie groups the analogous phenomenon is that 
non-compact matrix groups have non-trivial infinite dimensional 
unitary representations. 

We hope that the detailed discussion of three major classes of representations 
in Sections~\ref{sec:8}-\ref{sec:10} demonstrate the close interactions 
between convex geometry, complex analysis and Lie theory in the context 
of semibounded unitary representations. Presently, this theory is still 
in its infancy, but a general picture appears to evolve. 
One major point is that understanding semibounded unitary representations 
requires a good deal of knowledge on open invariant cones in the 
corresponding Lie algebra $\g$ and, what is closely related, 
information on the set $\g'_{\rm seq}$ 
of {\it semi-equicontinuous} coadjoint orbits $\cO_\lambda$, i.e., 
orbits for which the function $x \mapsto \sup \cO_\lambda(x)$ is bounded 
on some non-empty open subset of $\g$. Note that for any semibounded 
representation the momentum set $I_\pi$ is contained in 
$\g'_{\rm seq}$. In many important cases the methods developed in 
Section~\ref{sec:5} provide a complete description of $I_\pi$ in terms 
of generating coadjoint orbits. 

The main new results and aspects presented in this paper are: 
\begin{itemize}
\item Section~\ref{sec:5} provides tools to verify that unitary 
representations can be realized in spaces of holomorphic 
sections of line bundles and to calculate their momentum sets. 
Here Theorem~\ref{thm:5.7} is the key tool. 
\item The smoothness of the action of a 
Banach--Lie group on the space of smooth vectors proved in \cite{Ne10} 
is applied in our context in two essential ways. 
In Section~\ref{sec:4} it provides a Hamiltonian action on 
$\bP(\cH^\infty)$ for a unitary representation of a Banach--Lie group, 
and in Section~\ref{sec:7} it relates a $C^*$-dynamical system 
$(\cA, H, \alpha)$ and corresponding covariant representations 
satisfying a spectral condition to semibounded representations 
of  semidirect product Lie groups of the form $\U(\cA^\infty) \rtimes H$.  
This leads in particular to the remarkable conclusions in 
Theorem~\ref{thm:7.4}, which are based on elementary properties of 
invariant cones. 
\item The analysis of the convexity properties of adjoint and coadjoint 
orbits of $\cV(\bS^1)$ and $\vir$ is new. It leads in particular to an 
identification of the class of irreducible semibounded representations 
with the unitary highest weight representations and their duals 
(Theorem~\ref{thm:8.22}). 
\item The insights that the Fock representations of
 $\Sp_{\rm res}(\cH)$ and $\rO_{\rm res}(\cH)$ are semibounded 
seems to be new, and so are the results on adjoint and coadjoint orbits 
of the corresponding Banach--Lie groups. 
\end{itemize}

For finite dimensional groups, the first systematic 
investigation of unitary representations $(\pi, \cH)$ 
for which the cone $\{ x \in \g \: s_\pi(x) \leq 0\}$ 
is non-trivial for non-compact simple Lie groups (which are necessarily 
hermitian by \cite{Vin80}) has been undertaken in the pioneering work of 
G.~Olshanski (\cite{Ols82}). Based on the powerful structure theory for 
invariant cones developed in \cite{HHL89} by Hofmann, Hilgert and 
Lawson, we were eventually able to develop a general theory 
for semibounded representations of finite dimensional 
Lie groups, including a classification and a disintegration 
theory (\cite{Ne00}). We hope that one can also develop a similarly 
rich theory of complex semigroups and holomorphic extensions, 
 so that $C^*$-algebraic tools 
become available to deal with direct integrals of semibounded 
representations. In \cite{Ne08} we undertook some first steps in 
this direction, including a complete theory for the abelian 
case (cf.~Theorem~\ref{thm:ab}). 
What is needed here is a good theory of analytic vectors, 
which becomes a tricky issue for infinite dimensional Lie groups. 
Up to now, existence of analytic vectors is only known for very 
special classes of groups such as certain direct limits (\cite{Sa91}) 
and the canonical commutator relations in 
Quantum Field Theory (\cite{Re69}, \cite{He71}). 

If $K$ is a compact simple Lie group and $\cL(K) := C^\infty(\bS^1,K)$ 
the corresponding loop group, then the group 
$T_r := \T = \R/\Z$ acts smoothly by rotations on $\cL(K)$ and also on 
its canonical central extension $\tilde\cL(K)$ by $\T$, which 
leads to the ``smooth version'' $\hat\cL(K) := \tilde\cL(K)\rtimes\bT$ 
of affine Kac--Moody groups. For these groups positive energy 
representations are defined by requiring the 
spectrum of the generator of $T_r$ to be bounded below in 
the representation. 
Similarly, one defines positive energy representations 
of the group $\Diff(\bS^1)_+$ of orientation preserving diffeomorphisms 
of the circle. 
Various aspects of the theory of irreducible positive energy representations 
were developed in \cite{SeG81}, \cite{GW84}, \cite{PS86} and 
\cite{TL99a, TL99b}, but only in \cite[Sects.~9.3,11.4; Prop.~11.2.5]{PS86} 
one finds some attempts towards a decomposition theory. 

Positivity conditions for spectra 
also play a key role in I.~E.~Segal's concept 
of physical representations of the full unitary group $\U(\cH)$ 
of a Hilbert space~$\cH$, endowed with the strong operator topology 
(\cite{Se57}). 
Here the positivity requirements 
even imply boundedness of the representation, discrete decomposability 
and even a classification of the irreducible representations. 

\tableofcontents
\printglossary 

\section{Semi-equicontinuous convex sets} \mlabel{sec:2}

Let $E$ be a real locally convex space and $E'$ be its topological dual, i.e., 
the space of continuous linear functionals on~$E$. 
We write $\la \alpha, v \ra = \alpha(v)$ for the natural pairing 
$E' \times E \to \R$ and endow $E'$ with the weak-$*$-topology, i.e., 
the coarsest topology for which all linear maps 
$$ \eta_v \: E' \to \R, \quad \eta_v(\alpha) := \alpha(v) $$
are continuous. For a subset $X \subeq E'$, the set 
$$ B(X) := \{ v \in E \: \inf\la X, v \ra > - \infty\} $$ 
\glossary{name=$B(X)$,description={$=\{ v \in E \: \inf\la X, v \ra > - \infty\}$ for $X \subeq E'$}} is a convex cone which coincides with the domain 
of the {\it support function} 
$$ s_X \: E \to \R \cup \{ \infty\}, \quad s_X(v) 
:= - \inf \la X, v \ra = \sup \la X, - v\ra $$
of $X$ in the sense that $B(X) = s_X^{-1}(\R)$. As a sup 
of a family of continuous linear functionals, the function 
$s_X$ is convex, lower semicontinuous and positively homogeneous. 

\begin{defn} We call a subset $X \subeq E'$ {\it semi-equicontinuous} if 
$s_X$ is bounded on some non-empty open subset of $E$ (cf.\ \cite{Ne09a}). 
This implies in particular that the cone $B(X)$ has interior points 
and even that $s_X$ is continuous on $B(X)^0$ 
(\cite[Prop.~6.8]{Ne08}). 
\end{defn} 

If the space $E$ is barrelled, which includes in particular 
Banach and Fr\'echet spaces, we have the following handy 
criterion for semi-equicontinuity. We only have to apply 
\cite[Thm.~6.10]{Ne08} to the lower semicontinuous function 
$s_X$. 

\begin{prop} \mlabel{prop:seqcrit} 
If $E$ is barrelled, then $X \subeq E'$ is semi-equicontinuous if and 
only if $B(X)$ has interior points.
\end{prop}

\begin{rem}\mlabel{rem:1.1} (a) The notion of semi-equicontinuity 
generalizes the notion of equicontinuity, which is equivalent to 
$s_X$ being bounded on some $0$-neighborhood of $E$. 
In fact, the boundedness of $s_X$ on some symmetric 
$0$-neighborhood $U = - U$ means that there exists a $C > 0$ with 
$s_X(\pm v) \leq C$ for $v \in U.$
This is equivalent to 
$|\alpha(v)| \leq C$ for $v\in U, \alpha \in X,$
which means that $X$ is equicontinuous. 

(b) If $Y := \oline{\conv}(X)$ denotes the weak-$*$-closed convex 
hull of $X$, then $s_X = s_Y$, 
and, using the Hahn--Banach Separation Theorem, 
$Y$ can be reconstructed from $s_Y$ 
by 
$$Y = \{ \alpha \in E' \: (\forall v \in B(Y))\, 
\alpha(v) \geq \inf \la Y,v \ra = - s_Y(v)\}.$$
If, in addition, the interior $B(Y)^0$ is non-empty, then we even have 
$$Y = \{ \alpha \in E' \: (\forall v \in B(Y)^0)\, 
\alpha(v) \geq \inf \la Y,v \ra = - s_Y(v)\}$$
(\cite[Prop.~6.4]{Ne08}). 
\end{rem}

\begin{defn} (a) 
For a convex subset $C \subeq E$ we put 
$$ \lim(C) := \{ x \in E \: C + x \subeq C \} $$ 
\glossary{name={$\lim(C)$},description={recession cone of convex set $C$}}
\glossary{name={$H(C)$},description={translations preserving convex set $C$}}
and 
$$  
H(C) := \lim(C) \cap -\lim(C) := \{ x \in E \: C + x = C \}. $$
Then $\lim(C)$ is a convex cone and $H(C)$ a linear subspace of $E$. 

(b) A convex cone $W \subeq E$
is called {\it pointed} if $H(W) = \{0\}$. 

(c) For a subset $C \subeq E$, 
\glossary{name={$C^\star$},description={dual cone of $C$}}
$$  C^\star := \{ \alpha \in E' \: \alpha(C) \subeq \R_+\} $$
is called the {\it dual cone} and for a subset $X \subeq E'$, we 
define the {\it dual cone} by 
$$  X^\star := \{ v \in E \: \la X, v \ra \subeq \R_+\}. $$
\end{defn}

\begin{exs} \mlabel{ex:2.3} (a) If $E$ is a Banach space, then the unit ball 
$$X := \{ \alpha \in E' \: \|\alpha\| \leq 1 \} $$ 
in $E'$ is equicontinuous because the Hahn--Banach Theorems imply that 
$s_X(v) = \|v\|$ for $v\in E$. 

(b) If $\eset\not=\Omega \subeq E$ is an open convex cone, then 
its dual cone $\Omega^\star$ 
is semi-equicontinuous because we have 
$s_{\Omega^\star} = 0$ on $\oline\Omega = (\Omega^\star)^\star$ and 
$s_{\Omega^\star} = \infty$ on the complement of this closed cone. 
\end{exs}

We have just seen that open convex cones lead to semi-equicontinuous sets.
There is also a partial converse: 

\begin{rem} \mlabel{rem:conered} 
Let $X \subeq E'$ be a semi-equicontinuous set 
and $\tilde E := E \oplus \R$. For the set 
$$\tilde X := X \times \{1\} \subeq \tilde E' 
\quad \mbox{ we then have } \quad  s_{\tilde X}(v,t) := t + s_X(v), $$
so that the boundedness of $s_X$ on some non-empty open subset of $E$ implies 
that the interior of the dual cone $\tilde X^\star \subeq \tilde E$ 
is non-empty. In view of Example~\ref{ex:2.3}(b), this 
means that $X$ is semi-equicontinuous if and only if it can be embedded 
into the dual of some open convex cone in $\tilde E$. 
\end{rem}

The following observation shows 
that semi-equicontinuous convex sets share many  
important properties with compact convex sets 
(cf.\ \cite[Prop.~6.13]{Ne08}): 

\begin{proposition} \mlabel{prop:locomp} 
Let $X \subeq E'$ be a 
non-empty weak-$*$-closed convex subset and $v \in E$ such that 
the support function $s_X$ is bounded above on some neighborhood of $v$. 
Then $X$ is weak-$*$-locally compact, 
the function 
$$\eta_v \: X \to \R, \quad \eta_v(\alpha) := \alpha(v)$$ 
is proper, and there exists an extreme point $\alpha \in X$ 
with $\alpha(v) = \min \la X, v \ra$. 
\end{proposition}

We conclude this section with some elementary 
properties of convex subsets of locally convex spaces. 

The following lemma (\cite[Cor.~II.2.6.1]{Bou07}) is often useful: 
\begin{lem}
  \mlabel{lem:bou} 
For a convex subset $C$ of a locally convex space $E$ the following 
assertions hold: 
\begin{description}
\item[\rm(i)] $C^0$ and $\oline C$ are convex. 
\item[\rm(ii)] $\oline{C}^0 = C^0$ and if $C^0 \not=\eset$, then 
$\oline{C^0} = \oline C$. 
\end{description}
\end{lem}

\begin{lem} \mlabel{lem:limcone} 
If $\eset\not=C \subeq E$ 
is an open or closed convex subset, then the following assertions hold: 
\begin{description}
  \item[\rm(i)] $\lim(C) = \lim(\oline C)$ is a closed convex cone. 
  \item[\rm(ii)] $\lim(C) = \{ v \in E \: v = \lim_{n \to \infty} t_n c_n,
c_n \in C, t_n \to 0, t_n \geq 0 \}$. 
If $t_j c_j \to v$ holds for a net with $t_j \geq 0$ and 
$t_j \to 0$ and $c_j \in C$, then also $v \in \lim(C)$.  
  \item[\rm(iii)] If $c \in C$ and $d \in E$ satisfy 
$c + \N d \subeq C$, then $d \in \lim(C)$. 
\item[\rm(iv)] $c + \R d \subeq C$ implies $d \in H(C)$. In particular, 
$H(C) = \{0\}$ if and only if $C$ contains no affine lines. 
\item[\rm(v)] $H(C)$ is closed and the subset $C/H(C) \subeq E/H(C)$ 
contains no affine lines. 
\item[\rm(vi)] $B(C)^\star = \lim(C)$ and $B(C)^\bot = H(C)$. 
\end{description}
\end{lem}

\begin{prf} (cf.\ \cite[Prop.~6.1]{Ne08}) (i) If $C$ is open, 
then $C = (\oline C)^0$ by Lemma~\ref{lem:bou}, and 
thus 
$x + C \subeq C$ is equivalent to $x + \oline C \subeq \oline C$. 
In particular, $\lim(C)$ is closed. 

 (ii) If $c \in C$ and $x \in \lim(C)$, then $c + n x \in C$
for $n \in \N$ and  $ \frac{1}{n}(c + nx)\to~x. $ 
If, conversely, $x = \lim t_j c_j$ with $t_j \to 0$,
$t_j \geq 0, c_j \in C$, and $c \in \oline C$, then 
$(1 - t_j) c + t_j c_j \to c + x \in \oline C$
implies that $\oline C + x \subeq \oline C$, i.e.\ $x \in 
\lim(\oline C) = \lim(C)$. 

 (iii) In view of $\frac{1}{n}(c + n d) \to d$, (ii) implies 
$d \in \lim(C)$. 

(iv) immediately follows from (iii). 

(v) The closedness of $H(C) = \lim(C) \cap - \lim(C)$ follows from (i). 
Therefore the quotient topology on $E/H(C)$ is Hausdorff, so that 
the quotient topology turns $E/H(C)$ into a locally 
convex space. Let $q \: E \to E/H(C)$ denote the quotient map. 
If $y + \R d \subeq q(C) = C/H(C)$ is an affine line and 
$y = q(x), d = q(c)$, then $x + \R c \subeq q^{-1}(C/H(C)) = C$ 
implies $c \in H(C)$ by (iv), which leads to $d = q(c) = 0$. 
Hence $C/H(C)$ contains no affine lines. 

(vi) From 
$$\oline C = \{ v \in E \: (\forall \alpha \in B(C))\ 
\alpha(v) \geq \inf \alpha(C)\} $$
(a consequence of the Hahn--Banach Separation Theorem), 
we derive that \break 
$B(C)^\star \subeq \lim(\oline C) =\lim(C)$. 
Conversely, $\lim(C) \subeq B(C)^\star$ follows immediately 
from $c + \lim(C) \subeq C$ for each $c \in C$. 
\end{prf} 

\begin{rem} If $X \subeq E'$ is semi-equicontinuous, 
then $B(X)$ has interior points, so that 
$H(X) \subeq B(X)^\bot = \{0\}$ 
(Lemma~\ref{lem:limcone}(vi)) implies $H(X) = \{0\}$, i.e., 
$X$ contains no affine lines. 
If, conversely, $H(X) = \{0\}$ and 
$\dim V < \infty$, then it follows from 
\cite[Prop.~V.1.15]{Ne00} that $X$ is semi-equicontinuous. 
Therefore closed convex subsets of finite dimensional vector 
spaces are semi-equicontinuous if and only if they contain no affine 
lines. 
\end{rem}

For later applications, we record the following fact on fixed point 
projections for actions of compact groups. 

\begin{prop}\mlabel{prop:project}
Let $K$ be a compact group acting continuously on the 
complete locally convex space $E$ by the representation 
$\pi \: K \to \GL(E)$. 

{\rm(a)} If $\Omega \subeq E$ is an open or closed  
$K$-invariant convex subset, then $\Omega$ is invariant under the 
fixed point projection 
$$ p(v) := \int_K \pi(k)v\, d\mu_K(k), $$
where $\mu_K$ is a normalized Haar measure on $K$. 

{\rm(b)} If $C \subeq E'$ is a weak-$*$-closed convex 
$K$-invariant subset, then 
$C$ is invariant under the adjoint 
$$ p'(\lambda)v := \lambda(p(v)) = \int_K \lambda(\pi(k)v)\, d\mu_K(k) $$
of~$p$.
\end{prop} 

\begin{prf} (a) The existence of the integrals defining the 
projection $p$ follows from the 
completeness of $E$ (cf.\ \cite[Prop.~3.30]{HoMo98}). 
Let $\lambda \in B(\Omega) \subeq E'$ 
be a continuous linear functional bounded below on 
$\Omega$. 

If $\Omega$ is open, we then have 
$\lambda(\pi(k)v) > \inf \lambda(\Omega)$ for every $k \in K$, 
so that 
$$ \lambda(p(v)) = \int_K \lambda(\pi(k)v)\, d\mu_K(k) 
> \inf \lambda(\Omega). $$
In view of the Hahn--Banach Separation Theorem, this implies that 
$p(v) \in \oline\Omega$ (Remark~\ref{rem:1.1}(b)). If 
$p(v) \in \partial \Omega$, then 
\cite[Prop.~II.5.2.3]{Bou07} implies the existence of 
$\lambda \in B(\Omega) = B(\oline\Omega)$ 
with $\lambda(p(v)) = \min \lambda(\oline\Omega) = 
\inf \lambda(\Omega),$ a contradiction. 
Therefore  $p(v) \in \oline\Omega^0 = \Omega$ 
(cf.\ Lemma~\ref{lem:bou}). 

If $\Omega$ is closed, then the preceding argument 
implies $\lambda(p(v)) \geq \inf\lambda(\Omega)$, and hence that 
$p(v) \in \Omega$ by the Separation Theorem. 

(b) Now let $C \subeq E'$ be weak-$*$-closed and $G$-invariant. 
For each $v \in E$ and $\lambda \in C$, we then have 
$$p'(\lambda)(v) =  \lambda(p(v)) = \int_K \lambda(\pi(k)v)\, d\mu_K(k) 
\geq \inf \la C, v \ra, $$
so that the Hahn--Banach Separation Theorem shows that $p'(\lambda) \in C$.  
\end{prf}

\section{Infinite dimensional Lie groups} \mlabel{sec:3}

In this section we provide the definition of a locally convex Lie group 
and present several key examples that will show up later in our discussion 
of semibounded representations. 

\begin{defn} 
(a) Let $E$ and $F$ be locally convex spaces, $U
\subeq E$ open and $f \: U \to F$ a map. Then the {\it derivative
  of $f$ at $x$ in the direction $h$} is defined as 
$$ \dd f(x)(h) := (\partial_h f)(x) := \derat0 f(x + t h) 
= \lim_{t \to 0} \frac{1}{t}(f(x+th) -f(x)) $$
whenever it exists. The function $f$ is called {\it differentiable at
  $x$} if $\dd f(x)(h)$ exists for all $h \in E$. It is called {\it
  continuously differentiable}, if it is differentiable at all
points of $U$ and 
$$ \dd f \: U \times E \to F, \quad (x,h) \mapsto \dd f(x)(h) $$
is a continuous map. Note that this implies that the maps 
$\dd f(x)$ are linear (cf.\ \cite[Lemma~2.2.14]{GN09}). 
The map $f$ is called a {\it $C^k$-map}, $k \in \N \cup \{\infty\}$, 
if it is continuous, the iterated directional derivatives 
$$ \dd^{j}f(x)(h_1,\ldots, h_j)
:= (\partial_{h_j} \cdots \partial_{h_1}f)(x) $$
exist for all integers $j \leq k$, $x \in U$ and $h_1,\ldots, h_j \in E$, 
and all maps \break $\dd^j f \: U \times E^j \to F$ are continuous. 
As usual, $C^\infty$-maps are called {\it smooth}. 

  (b) If $E$ and $F$ are complex locally convex spaces, then a map $f$ is 
called {\it complex analytic} if it is continuous and for each 
$x \in U$ there exists a $0$-neighborhood $V$ with $x + V \subeq U$ and 
continuous homogeneous polynomials $\beta_k \: E \to F$ of degree $k$ 
such that for each $h \in V$ we have 
$$ f(x+h) = \sum_{k = 0}^\infty \beta_k(h), $$
as a pointwise limit (\cite{BoSi71}). 
The map $f$ is called {\it holomorphic} if it is $C^1$ 
and for each $x \in U$ the 
map $\dd f(x) \: E \to F$ is complex linear (cf.\ \cite[p.~1027]{Mil84}). 
If $F$ is sequentially complete, then $f$ is holomorphic if and only if 
it is complex analytic (cf.\ \cite{Gl02}, \cite[Ths.~3.1, 6.4]{BoSi71}). 

(c) If $E$ and $F$ are real locally convex spaces, 
then we call $f$ {\it real analytic}, resp., $C^\omega$, 
if for each point $x \in U$ there exists an open neighborhood 
$V \subeq E_\C$ and a holomorphic map $f_\C \: V \to F_\C$ with 
$f_\C\res_{U \cap V} = f\res_{U \cap V}$  (cf.\ \cite{Mil84}). 
The advantage of this definition, which differs from the one in 
\cite{BoSi71}, is that it works nicely for non-complete spaces, 
any analytic map is smooth, 
and the corresponding chain rule holds without any condition 
on the underlying spaces, which is the key to the definition of 
analytic manifolds (see \cite{Gl02} for details).
\end{defn}

Once one has introduced the concept of a smooth function 
between open subsets of locally convex spaces, it is clear how to define 
a locally convex smooth manifold. 
A {\it (locally convex) Lie group} $G$ is a group equipped with a 
smooth manifold structure modeled on a locally convex space 
for which the group multiplication and the 
inversion are smooth maps. We write $\1 \in G$ for the identity element and 
$\lambda_g(x) = gx$, resp., $\rho_g(x) = xg$ for the left, resp.,
right multiplication on $G$. Then each $x \in T_\1(G)$ corresponds to
a unique left invariant vector field $x_l$ with 
$x_l(g) := T_\1(\lambda_g)x, g \in G.$
The space of left invariant vector fields is closed under the Lie
bracket of vector fields, hence inherits a Lie algebra structure. 
In this sense we obtain on $\g := T_\1(G)$ a continuous Lie bracket which
is uniquely determined by $[x,y]_l = [x_l, y_l]$ for $x,y \in \g$. 
We shall also use the functorial 
notation $\L(G) := (\g,[\cdot,\cdot])$ 
for the Lie algebra of $G$ and, accordingly, 
$\L(\phi) = T_\1(\phi)\: \L(G_1) \to \L(G_2)$ 
for the Lie algebra morphism associated to 
a morphism $\phi \: G_1 \to G_2$ of Lie groups. 
Then $\L$ defines a functor from the category 
of locally convex Lie groups to the category of locally convex topological 
Lie algebras. The adjoint action of $G$ on $\L(G)$ is defined by 
$\Ad(g) := \L(c_g)$, where $c_g(x) = gxg^{-1}$. This action is smooth and 
each $\Ad(g)$ is a topological isomorphism of $\L(G)$. The coadjoint action 
on the topological dual space $\L(G)'$ is defined by 
$$\Ad^*(g)\alpha := \alpha \circ \Ad(g)^{-1}$$ 
and the maps $\Ad^*(g)$ are continuous with respect to the 
weak-$*$-topology on $\L(G)'$, but in general the coadjoint 
action of $G$ is {\sl not} continuous with respect to this topology. 
If $\g$ is a Fr\'echet, resp., a Banach space, then 
$G$ is called a {\it Fr\'echet-}, resp., a 
{\it Banach--Lie group}. 

A smooth map $\exp_G \: \L(G) \to G$  is called an {\it exponential function} 
if each curve $\gamma_x(t) := \exp_G(tx)$ is a one-parameter group 
with $\gamma_x'(0)= x$. The Lie group $G$ is said to be 
{\it locally exponential} 
if it has an exponential function for which there is an open $0$-neighborhood 
$U$ in $\L(G)$ mapped diffeomorphically by $\exp_G$ onto an 
open subset of $G$. All Banach--Lie groups 
are locally exponential (\cite[Prop.~IV.1.2]{Ne06}). 
Not every infinite dimensional Lie group has an exponential 
function (\cite[Ex.~II.5.5]{Ne06}), but exponential functions 
are unique whenever they exist. 

In the context of unitary representation theory, the exponential function 
permits us to associate to each element $x$ of the Lie algebra a unitary 
one-parameter group $\pi_x(t) := \pi(\exp_G tx)$. 
We therefore {\bf assume} in the following that 
$G$ has an exponential function. 


\begin{exs} \mlabel{ex:liegrp} 
Here are some important examples of infinite dimensional 
Lie groups that we shall encounter below. 

(a) (Unitary groups) If $\cA$ is a unital $C^*$-algebra, then 
its unitary group 
$$ \U(\cA) := \{ g \in \cA \: g^* g = gg^* = \1\} $$
is a Banach--Lie group with Lie algebra 
$$ \fu(\cA) = \{ x \in \cA \: x^* = - x\}. $$
In particular, the unitary group 
$\U(\cH) = \U(B(\cH))$ of a complex Hilbert space $\cH$ is of this form, 
and we write $\fu(\cH)$ for its Lie algebra.  

As we shall see below (Definition~\ref{def:7.3}), 
in some situations one is forced to consider more 
general classes of algebras: A locally convex topological 
unital algebra $\cA$ is called a {\it continuous inverse algebra} 
if its unit group $\cA^\times$ is open and the inversion map 
$a \mapsto a^{-1}$ is continuous. This condition already implies 
that $\cA^\times$, endowed with the canonical manifold 
structure as an open subset, is a Lie group 
(cf.\ \cite{Gl02}). If, in addition, $\cA$ is a complex algebra and 
$*$ a continuous algebra involution, then 
$$ \U(\cA) := \{ g \in \cA \: g^* g = gg^* = \1\} $$
is a closed subgroup, and even a submanifold, as is easily seen with the 
Cayley transform $c(x) = (\1 -x)(\1 + x)^{-1}$. It defines an involutive  
diffeomorphism of some open neighborhood of $\1$ onto some 
open neighborhood of $0$ satisfying $c(x)^* = c(x^*)$ and $c(x^{-1}) = -c(x)$. 
In particular $c(g)$ is skew-hermitian if and only if $g$ is unitary, 
and since $\U(\cA)$ is a subgroup of $\cA^\times$, this argument 
shows that it actually is a Lie subgroup.

(b) (Schatten class groups) If $\cH$ is a real or complex Hilbert space and 
$B_p(\cH)$ 
\glossary{name={$B_p(\cH)$},description={Schatten ideal of order $p$}}
denotes the $p$-Schatten ideal $(p \geq 1)$ with 
the norm $\|A\|_p := \tr((A^*A)^{p/2})^{1/p}$, then 
$$ \U_p(\cH) := \U(\cH) \cap (\1 + B_p(\cH)) $$
is a Banach--Lie group with Lie algebra  
\glossary{name={$\fu_p(\cH)$},description={skew-hermitian elements 
of $B_p(\cH)$; Lie algebra of $\U_p(\cH)$}}
$$ \fu_p(\cH) := \fu(\cH) \cap B_p(\cH) $$
(cf.\ \cite{Mick89}, \cite{Ne04}). 

(c) (Restricted groups) If $P$ is an orthogonal projection 
on $\cH$ and $G \subeq \GL(\cH)$ a subgroup, then we call 
$$ G_{\rm res} := \{ g \in G \: [g,P] \in B_2(\cH)\} $$
the corresponding restricted group. Using the fact that $B_2(\cH)$ 
is invariant under left and right multiplication 
with elements of $G$, it is easy to 
see that this is indeed a subgroup, and in many cases it carries a natural 
Banach--Lie group structure. Writing $\cH = \im(P) \oplus \ker(P)$ and, 
accordingly, operators on $\cH$ as $(2 \times 2)$-matrices, 
then \
$$ g= \pmat{a & b \\ c &d} \in G_{\rm res} 
\quad \Longleftrightarrow \quad b,c \in B_2(\cH). $$
In particular, we have the {\it restricted unitary group} 
$$ \U_{\rm res}(\cH,P) := \{ g \in \U(\cH) \: [g,P] \in B_2(\cH)\}. $$
\glossary{name={$\U_{\rm res}(\cH,P)$},description={restricted 
unitary group defined by $P$}}

(d) If $\cH$ is a complex Hilbert space, then the scalar product 
$\la \cdot, \cdot \ra$ (always assumed to be linear in the first 
component), defines two real bilinear forms 
$$ \beta(x,y) := \Re \la x, y \ra 
\quad \mbox{ and } \quad 
\omega(x,y) := \Im \la x,y \ra = \Re \la x,Iy\ra,$$
where $\beta$ is symmetric and 
$\omega$ is skew-symmetric. 
Writing $\cH_\R$ for the underlying real Banach space, we thus 
obtain the {\it symplectic group}  \glossary{name={$\cH_\R$},description={real Hilbert space underlying the complex Hilbert space $\cH$}}
\glossary{name={$\Sp(\cH)$},description={symplectic 
group of real Hilbert space $\cH_\R$}} 
\begin{align*}
\Sp(\cH) 
&:= \Sp(\cH_\R,\omega) := \{ g \in \GL(\cH_\R) \: 
(\forall v, w \in \cH_\R)\, \omega(gv,gw) = \omega(v,w) \} \\
&=    \{ g \in \GL(\cH_\R) \: g^\top I g = I\} 
\end{align*} 
\glossary{name={$\rO(\cH)$},description={orthogonal 
group of real Hilbert space $\cH_\R$}} 
and the {\it orthogonal group} 
\begin{align*} 
\OO(\cH) 
&:= \OO(\cH_\R,\beta) := \{ g \in \GL(\cH_\R) \: 
(\forall v, w \in \cH_\R)\, \beta(gv,gw) = \beta(v,w) \} \\
&=    \{ g \in \GL(\cH_\R) \: g^\top g = \1\}.  
\end{align*}
There are two important variants of these groups, namely the 
Hilbert--Lie groups 
$$ \Sp_2(\cH) := \Sp(\cH) \cap (\1 + B_2(\cH_\R)), \quad 
\OO_2(\cH) := \OO(\cH) \cap (\1 + B_2(\cH_\R)) $$
and the restricted groups 
$$ \Sp_{\rm res}(\cH) := 
\{ g \in \Sp(\cH) \: [g,I] \in B_2(\cH_\R) \} $$
and $$ \OO_{\rm res}(\cH) := 
\{ g \in \OO(\cH) \: [g,I] \in B_2(\cH_\R) \}, $$
where $Iv := iv$ denotes the complex structure on $\cH_\R$ defining the 
complex Hilbert space~$\cH$. 
\glossary{name={$\Sp_{\rm res}(\cH)$},description={restricted 
symplectic group of $\cH_\R$}} 
\glossary{name={$\OO_{\rm res}(\cH)$},description={restricted 
orthogonal group of $\cH_\R$}} 

The groups 
$\Sp(\cH_\R, \omega)$ and $\OO(\cH_\R, \beta)$ play a key role 
in Quantum Field Theory as the automorphism groups of the 
canonical commutation relations (CCR) and the 
canonical anticommutation relations (CAR) (cf.~\cite{BR97}). However, 
only the corresponding restricted groups, resp., their 
 central extensions, have corresponding 
unitary representations (cf.\ Sections~\ref{sec:9} and \ref{sec:10}).

(e) The group $\Diff(M)^{\rm op}$ of diffeomorphisms of a compact manifold 
$M$ is a Lie group with respect to the 
group structure defined by $\phi \cdot \psi := \psi \circ \phi$. 
Its Lie algebra is the space $\cV(M)$ of smooth vector fields on $M$ 
with respect to the natural Lie bracket. 
The use of the opposite group simplifies many formulas 
and minimizes the number of negative signs. 
In particular, it implies that the exponential function is given by 
the time $1$-flow and not by its inverse. 
As we shall see below, this convention also 
leads to simpler formulas because the action of $\Diff(M)^{\rm op}$ 
by pullbacks is a left action. 

(f) If $K$ is a Lie group and $M$ is a compact manifold, then 
the space $C^\infty(M,K)$ of smooth maps is a Lie group 
with Lie algebra $C^\infty(M,\fk)$, where $\fk$ is the 
Lie algebra of $K$. 

(g) A domain $\cD$ in the complex Banach space $V$ is said to be 
{\it symmetric} if there exists for each point $x \in \cD$ 
a biholomorphic involution $s_x \in \Aut(\cD)$ for which 
$x$ is an isolated fixed point, or, equivalently 
$T_x(s_x) = - \id$. It is called a {\it symmetric Hilbert domain} 
if, in addition, $V$ is a complex Hilbert space. 
Then the group $\Aut(\cD)$ of biholomorphic automorphisms 
of $\cD$ carries a natural Banach--Lie group structure 
(cf.\ \cite[Sect.~V]{Ka97}, \cite{Up85}, \cite[Thm.~V.11]{Ne01a}).

Here is a typical example. 
If $\cH_\pm$ are two complex Hilbert spaces 
and $B_2(\cH_+, \cH_-)$ is the Hilbert space of Hilbert--Schmidt operators 
from $\cH_+$ to $\cH_-$, then 
$$ \cD := \{ z \in B_2(\cH_+, \cH_-) \: \|z\| < 1 \} $$
is a symmetric Hilbert domain, where $\|\cdot\|$ denotes the operator norm, 
which is smaller than the Hilbert--Schmidt norm $\|\cdot\|_2$. 
In particular, $\cD$ is unbounded if both spaces are infinite dimensional. 
On the Hilbert space $\cK := \cH_- \oplus \cH_+$ we define a 
hermitian form by 
$\gamma((v,w), (v',w')) := \la v,v' \ra - \la w,w'\ra$ and 
write 
$$ \U(\cH_+, \cH_-) \:= \{ g \in \GL(\cK) \: (\forall x,y\in \cK)\,  
\gamma(gx,gy) = \gamma(x,y)\} $$
for the corresponding {\it pseudo-unitary group}. 
From the projection 
$P(v_-,v_+) := (v_-,0)$ we now obtain a 
restricted pseudo-unitary group 
\glossary{name={$\U_{\rm res}(\cH_+, \cH_-)$},description={
restricted pseudo-unitary group}} 
$$\U_{\rm res}(\cH_+, \cH_-) 
= \{ g \in \U(\cH_+, \cH_-) \: \|[g,P]\|_2 < \infty\} $$ 
(cf.\ (c) above), and this group acts on $\cD$ by 
$$ g.z = \pmat{a & b \\ c &d}.z := (az + b)(cz+d)^{-1}. $$
The subgroup $\T \1$ of $\U_{\rm res}(\cH_+, \cH_-)$ acts trivially, 
and we thus obtain an isomorphism 
$$ \U_{\rm res}(\cH_+, \cH_-)/\T \1 \cong \Aut(\cD)_0$$  
(combine \cite[Sect.~5]{Ka97} on the automorphisms of the completion 
with respect to the operator norm with the extension result in 
\cite[Thm.~V.11]{Ne01a}). 

(h) In (g), the case $\cH=\cH_+ = \cH_-$ is of particular interest 
because it leads to two other series of irreducible 
symmetric Hilbert domains generalizing the operator balls 
in symmetric, resp., skew-symmetric matrices. 
Let $\sigma$ be an antilinear isometric involution on $\cH$ 
and define $x^\top := \sigma x^* \sigma$ for $x \in B(\cH)$. 
Then the map $\iota \: \cH \to \cK, v \mapsto \frac{1}{\sqrt 2}(v, \sigma v)$ 
is an antilinear isometric embedding whose image is an orthogonal 
direct sum $\cK = \iota(\cH) \oplus i \iota(\cH)$, so that 
$\cK \cong (\cH_\R)_\C$ as Hilbert spaces, where the complex conjugation 
on $\cK$ is given by $\tau(v,w) = (\sigma w, \sigma v)$. 
By complex linear extension, 
we thus obtain an embedding $\gamma \: \Sp(\cH)\into \GL(\cK)$, whose image 
preserves the real subspace $\iota(\cH)$ and 
the complex bilinear skew-symmetric form $\omega_\C$ obtained by 
complex bilinear extension of $\omega(x,y) = \Im \la x,y\ra$. 
Since the complex bilinear form 
$$ \omega_\C((x,y), (x',y')) := 
\la x, \sigma y' \ra - \la x', \sigma y\ra $$ 
on $\cK$ satisfies 
$$ \omega_\C((x,\sigma x), (y, \sigma y)) 
= \la x, y\ra - \la y,x \ra 
= 2i \Im \la x, y \ra = 2i \omega(x,y), $$
the subgroup $\gamma(\Sp(\cH))$ is contained in $\Sp(\cK,\omega_\C)$. 
From the fact that it also commutes with $\tau$ 
one easily derives that it also preserves the canonical 
hermitian form 
$\gamma((x,y), (x',y')) = \la x,x'\ra - \la y,y'\ra$, and this leads 
to an isomorphism 
\begin{equation}
  \label{eq:spcomplex}
\Sp(\cH) \cong \gamma(\Sp(\cH)) = \U(\cH, \cH) \cap \Sp(\cK, \omega_\C) 
\end{equation}
(cf.\ \cite[Rem.~I.2]{Ne02a}, \cite[Sect.~IV]{NO98}). 
This in turn leads to isomorphisms 
$$ \Sp_{\rm res}(\cH) \cong \U_{\rm res}(\cH, \cH) \cap \Sp(\cK, \omega_\C)$$
and 
\begin{equation}
  \label{eq:spresliealg}
\fsp_{\rm res}(\cH) \cong 
\Big\{ \pmat{ a & b \\ b^* & -a^\top} \: a \in \fu(\cH), b = b^\top 
 \in B_2(\cH)\Big\}. 
\end{equation}
From that one further derives that $\Sp_{\rm res}(\cH)$ 
acts by holomorphic automorphisms on the symmetric Hilbert domain 
$$ \cD_s := \{ z \in B_2(\cH) \: z^\top = z, \|z\| < 1  \} $$
by fractional linear transformations, 
and we thus obtain isomorphisms 
$$\Sp_{\rm res}(\cH)/\{\pm \1\} \cong \Aut(\cD)_0 
\quad \mbox{ and } \quad 
\Sp_{\rm res}(\cH)/\U(\cH) \cong \cD_s. $$

(i) For the complex symmetric bilinear form 
$$\beta_\C((x,y), (x',y')) := 
\la x, \sigma y' \ra + \la x', \sigma y\ra $$ 
we similarly obtain 
\begin{equation}
  \label{eq:ocomplex}
\rO(\cH) \cong \U(\cH, \cH) \cap 
\rO(\cK, \beta_\C) 
\end{equation}
which in turn leads to 
$$ \rO_{\rm res}(\cH) \cong \U_{\rm res}(\cH, \cH) \cap \rO(\cK, \beta_\C)$$
and 
\begin{equation}
  \label{eq:oresliealg}
\fo_{\rm res}(\cH) \cong 
\Big\{ \pmat{ a & b \\ -b^* & -a^\top} \: a \in \fu(\cH), b^\top = - b
 \in B_2(\cH)\Big\}. 
\end{equation}
\end{exs}

Since we shall need it several times in the following, we recall some 
basic facts on the adjoint and the coadjoint action on a centrally extended 
Lie algebra. 

\begin{rem} \mlabel{rem:centext} 
Let $\hat\g = \R \oplus_\omega \g$ be  a central 
extension of the Lie algebra $\g$ defined by the 
$2$-cocycle $\omega$, i.e.,  
$$ [(z,x), (z',x')] = (\omega(x,x'), [x,x']) 
\quad \mbox{ for } \quad z,z' \in \R, x,x' \in \g. $$
Then the adjoint action of $\hat\g$ factors through a representation 
$\ad_{\hat\g} \: \g \to\der(\hat\g)$, given by 
$\ad_{\hat\g}(x)(z,y) = (\omega(x,y), [x,y])$, which implies that 
$\g \to \g', x \mapsto i_x\omega = \omega(x,\cdot)$ is a $1$-cocycle. 

If $G$ is a corresponding Lie group to which the action of 
$\g$ on $\hat\g$ integrates as a smooth linear action, 
then it is of the form 
\begin{equation}
  \label{eq:centad}
\Ad_{\hat\g}(g)(z,y) = (z + \Theta_\omega(g)(\Ad(g)y), \Ad(g)y) 
= (z - \Theta_\omega(g^{-1})(y), \Ad(g)y), 
\end{equation}
where $\Theta_\omega \: G \to \g'$ is a $1$-cocycle with 
$T_\1(\Theta_\omega)x = i_x\omega$. 
The uniqueness of the representation of  $G$ on $\hat\g$ follows from 
the general fact that, for a connected Lie group, smooth representations 
are uniquely determined by their derived representations 
(\cite[Rem.~II.3.7]{Ne06}). The existence for simply 
connected $G$ follows from \cite[Prop.~VII.6]{Ne02b}. 
The corresponding dual representation of $G$ on 
$\hat\g' \cong \R \times \g'$ is then given by 
\begin{equation}
  \label{eq:centcoad}
\Ad^*_{\hat\g}(g)(z,\alpha) = (z, \Ad^*(g)\alpha -  z\Theta_\omega(g)).
\end{equation}

As these formulas show, the passage from the adjoint and coadjoint 
action of $\g$ to the $G$-action on $\hat\g$ is completely encoded in 
the cocycle $\Theta_\omega \: G \to \g'$. 
\end{rem}

\section{Momentum sets of smooth unitary representations} \mlabel{sec:4}

In this section we introduce the concept of a semibounded unitary 
representation of a Lie group~$G$. A key tool to study these representations 
is the momentum map $\Phi \: \bP(\cH^\infty) \to \g'$. According to 
Theorem~\ref{thm:4.4}, this map is also a momentum map in the 
classical sense of differential geometry, provided $G$ is a Banach--Lie 
group. 

\begin{defn} A {\it unitary representation of $G$} is a pair 
$(\pi, {\cal H})$ of a complex Hilbert space $\cH$ and a group homomorphism 
$\pi \: G \to \U(\cH)$. It is said to be {\it continuous} 
if the action map 
$G \times \cH \to \cH, (g,v) \mapsto \pi(g)v$ is continuous, which, 
since $G$ acts by isometries, is equivalent to the continuity of  
the orbit maps $\pi^v \: G \to \cH, g \mapsto \pi(g)v.$ 
We write 
$${\cal H}^\infty := \{ v \in \cH \: \pi^v \in C^\infty(G,\cH) \} $$ 
for the subspace of {\it smooth vectors}. 
The representation 
$(\pi, {\cal H})$ is said to be 
{\it smooth} if ${\cal H}^\infty$ is dense in ${\cal H}$.
On  $\cH^\infty$ the 
derived representation $\dd\pi$ of the Lie algebra $\g = \L(G)$ 
is defined by 
\glossary{name={$\dd\pi(x)$},description={derived representation}}
$$ \dd\pi(x)v := \derat0 \pi(\exp_G tx)v $$ 
(cf.\ \cite[Rem.~IV.2]{Ne01b}). 
If $(\pi, \cH)$ is smooth, then the 
invariance of $\cH^\infty$ under $\pi(G)$ implies that 
the operators $i \dd\pi(x)$, $x \in \g$, on this space 
are essentially selfadjoint (cf.\ \cite[Lemma 5.6]{Ne08}, 
\cite[Thm.~VIII.10]{RS75}). 
\end{defn} 

\begin{rem} (a) If $(\pi, \cH)$ is a smooth unitary representation, 
then the space $\cH^c := \{ v \in \cH \: 
\pi^v \in C(G,\cH) \} \supeq \cH^\infty$ of continuous vectors 
is dense, and since this 
space is closed (by the uniform boundedness of $\pi(G)$), it 
follows that $\cH = \cH^c$. This in turn implies that the 
$G$-action on $\cH$ is continuous. This means that smooth representations 
are in particular continuous. 

(b) If $G$ is finite dimensional, then 
G\aa{}rding's Theorem asserts that every continuous unitary representation 
of $G$ 
is smooth. However, this is false for infinite dimensional Lie 
groups. The representation of the additive 
Banach--Lie group $G := L^2([0,1],\R)$ on 
$\cH = L^2([0,1],\C)$ by $\pi(g)f := e^{ig}f$ is continuous 
with $\cH^\infty = \{0\}$ 
(\cite{BN08}). 
\end{rem}

\begin{defn} (a) Let $\bP({\cal H}^\infty) = \{ [v] := 
\C v \: 0 \not= v \in 
{\cal H}^\infty\}$ 
denote the projective space of the subspace ${\cal H}^\infty$ 
of smooth vectors. The map 
\glossary{name={$\Phi_\pi$},description={momentum map of $\pi$}}
$$ \Phi_\pi \: \bP({\cal H}^\infty)\to \g' \quad \hbox{ with } \quad 
\Phi_\pi([v])(x) 
= \frac{1}{i}  \frac{\la  \dd\pi(x).v, v \ra}{\la v, v \ra} $$
is called the {\it momentum map of the unitary representation $\pi$}.  
The right hand side is well defined because it
only depends on $[v] = \C v$. The operator $i\dd\pi(x)$ is symmetric so
that the right hand side is real, and since $v$ is a smooth vector, 
it defines a continuous linear functional on $\g$. 
We also observe that we have a natural action of $G$ on 
$\bP(\cH^\infty)$ by $g.[v] := [\pi(g)v]$, and the relation 
$$ \pi(g) \dd\pi(x) \pi(g)^{-1} = \dd\pi(\Ad(g)x) $$
immediately implies that $\Phi_\pi$ is equivariant with respect 
to the coadjoint action of $G$ on $\g'$. 

\glossary{name={$I_\pi$},description={$\oline{\conv}(\im(\Phi_\pi))$, momentum set of $\pi$}}
\par (b) The weak-$*$-closed convex hull 
$I_\pi \subeq \g'$ of the image of $\Phi_\pi$ is called the 
{\it (convex) momentum set of $\pi$}. In view of the equivariance 
of $\Phi_\pi$, it is an $\Ad^*(G)$-invariant subset of $\g'$. 
\end{defn}

For the following theorem, we recall the definition of a 
Hamiltonian action: 
\begin{defn} Let $\sigma \: G \times M \to M$ be a smooth 
action of the Lie group $G$ on the smooth manifold $M$ (neither is assumed 
to be finite dimensional) and suppose that $\omega$ is a $G$-invariant 
closed $2$-form on $M$. Then this action is said to be 
{\it Hamiltonian} if there exists a map 
$J \: \g \to C^\infty(M,\R)$ for which 
$\dd J_x = i_{\dot\sigma(x)}\omega$ holds for the derived 
Lie algebra homomorphism $\dot\sigma \: \g \to \cV(M)$, defined by 
$\dot\sigma(x)_p = - T_{(\1,p)}(\sigma)(x,0)$. 
The map $\Phi \: M \to \g', \Phi(m)(x) := J_x(m)$ 
is then called the corresponding {\it momentum map}. Of particular 
interest are momentum maps which are equivariant with respect to the 
coadjoint action. 
\end{defn}

\begin{thm}
  \mlabel{thm:4.4} 
Let $(\pi, \cH)$ be a unitary representation of a Banach--Lie group $G$. 
Then the seminorms 
$$ p_n(v) := \sup \{ \|\dd\pi(x_1) \cdots \dd\pi(x_n)v\| \: 
x_i \in \g, \|x_i\|\leq  1\} $$
define on $\cH^\infty$ the structure of a Fr\'echet space 
with respect to which the action of $G$ is smooth. 
Accordingly, the projective 
space $\bP(\cH^\infty)$ carries the structure of a complex 
Fr\'echet manifold on which $G$ acts smoothly by holomorphic maps. 
The Fubini--Study metric on $\bP(\cH)$ induces on 
$\bP(\cH^\infty)$ the structure of a weak K\"ahler manifold 
whose corresponding weak symplectic form $\Omega$ is given for any 
unit vector $v \in \cH^\infty$ by 
$$ \Omega_{[v]}(T_v(q)x, T_v(q)y) = -2\Im \la x, y\ra 
\quad \mbox{ for } \quad x,y \in v^\bot, $$ 
where $q \: \cH^\infty \setminus \{0\} \to \bP(\cH^\infty), v \mapsto [v],$ 
is the canonical projection. With respect to this symplectic 
form, the action of $G$ on $\bP(\cH^\infty)$ is Hamiltonian 
with momentum map $\Phi_\pi$. 
\end{thm}

\begin{prf} The smoothness of the action of $G$ on $\cH^\infty$ follows from 
Theorem~\ref{thm:4.1.3}. This implies that 
the natural charts defined by projections of affine 
hyperplanes define on $\bP(\cH^\infty)$ a complex manifold 
structure for which $G$ acts smoothly by holomorphic maps 
(see \cite[Prop.~V.2]{Ne01b} for details).  
Moreover, the Fubini--Study metric on $\bP(\cH)$ restricts to a 
(weak) K\"ahler metric on $\bP(\cH^\infty)$ which is invariant 
under the $G$-action. It is determined by the property that the 
projection map $q$ 
satisfies for any unit vector $v$ and $x,y \in v^\bot$ the relation 
$$ \la T_v(q)x, T_v(q)y\ra = \la x,y \ra. $$
In particular, we see that $\Omega$ defines a weak 
symplectic $2$-form on $\bP(\cH^\infty)$  
(cf.\ \cite[Sect.~5.3]{MR99}). For $x \in \g$, the smooth function 
$$ H_x([v]) := \frac{1}{i} \frac{\la \dd\pi(x)v, v \ra}{\la v,v\ra} $$
on $\bP(\cH^\infty)$ 
now satisfies for $y \in v^\bot$ and $\|v\|=1$: 
\begin{align*}
 \dd H_x([v])(T_v(q)y) 
&= 2\Re \la -i\dd\pi(x)v,y\ra 
= 2 \Im \la \dd\pi(x)v,y\ra \\
&= -\Omega_{[v]}(T_v(q)\dd\pi(x)v, T_v(q)y), 
\end{align*}
i.e., 
$\dd H_x = i_X \Omega$
for the smooth vector field on $\bP(\cH^\infty)$ 
defined by 
$X([v]) := -T_v(q)\dd\pi(x)v = \dot\sigma(x)([v])$. 
This means that the action of $G$ on $\bP(\cH^\infty)$ is Hamiltonian and 
$\Phi_\pi \: \bP(\cH^\infty) \to \g'$ is a corresponding momentum 
map. 
\end{prf}

The main new point of the preceding theorem is 
that it provides for any Banach--Lie group 
a natural analytic setup for which the momentum 
map $\Phi_\pi$ really is a momentum map for a smooth action on 
a weak symplectic manifold. The corresponding result 
for the action of the Banach--Lie group $\U(\cH)$ 
on $\bP(\cH)$ is well-known (cf.\ \cite[Ex.~11.4(h)]{MR99}). 

A key property of the momentum set $I_\pi$ 
is that it provides complete information on the extreme spectral 
values of the operators $i\dd\pi(x)$: 
\begin{equation}
  \label{eq:momspec}
\sup(\Spec(i\dd\pi(x))) = s_{I_\pi}(x)= - \inf \la I_\pi, x \ra 
 \quad \mbox{ for } \quad x \in \g 
\end{equation}
(cf.\ \cite[Lemma 5.6]{Ne08}). 
This relation immediately entails the equivalences 
\begin{equation}
  \label{eq:momspec2}
  x \in I_\pi^\star \quad \Longleftrightarrow \quad 
s_{I_\pi}(x) \leq 0 \quad \Longleftrightarrow \quad -i\dd\pi(x) \geq 0. 
\end{equation}
It is now natural to study those representations for which 
$s_\pi$, resp., the set $I_\pi$, contains the most significant 
information, which leads to the following concept: 

\begin{defn} A smooth unitary representation $(\pi, \cH)$ of a 
Lie group $G$ is said to be {\it semibounded} if its momentum 
set $I_\pi$ is semi-equicontinuous, i.e., if 
$$s_\pi(x) := s_{I_\pi}(x) = \sup(\Spec(i\dd\pi(x))) $$ 
is bounded from above on some non-empty open subset of 
$\g$. 
\glossary{name={$s_\pi(x)$},description={$= \sup(\Spec(i\dd \pi(x)))$, 
for smooth unit. rep. $\pi$}}
We call $(\pi, \cH)$ {\it bounded} if 
$I_\pi$ is equicontinuous. 

The representation $(\pi, \cH)$ is said to satisfy the 
{\it positive energy condition} with respect to some $d \in \g$ if 
$i\dd\pi(d) \geq 0$.  
It satisfies the {\it positive energy condition} with respect to some 
convex cone $C \subeq \g$, if 
$i\dd\pi(x) \geq 0$ holds for each $x \in C$ (cf.\ \cite{SJOSV78}). 
\end{defn}

For a semibounded representation $(\pi, \cH)$, 
the domain $s_\pi^{-1}(\R)$ of $s_\pi$ is a convex cone 
with non-empty interior and $s_\pi$ is continuous on this 
open cone (cf.\ \cite[Prop.~6.8]{Ne08}). 
Since the momentum set $I_\pi$ is invariant under the coadjoint action, 
the function $s_\pi$ and its domain are 
invariant under the adjoint action. This leads to two open 
invariant convex cones in $\g$: 
$$ W_\pi := \{ x \in \g \: s_\pi(x) < \infty\}^0 = B(I_\pi)^0 
\quad \mbox{ and } \quad 
C_\pi := \{ x \in \g \: s_\pi(x) \leq 0 \}^0, $$
where ${}^0$ denotes the interior of a 
set. 
\glossary{name={$W_\pi$},description={$=B(I_\pi)^0$, for unit. rep. $\pi$}}
\glossary{name={$C_\pi$},description={interior of $(I_\pi)^\star$, for unit. rep. $\pi$}}
Note that \eqref{eq:momspec2} implies that $C_\pi$ is the 
interior of the dual cone $I_\pi^\star$. 

We collect some basic properties of these two cones: 

\begin{prop} \mlabel{prop:repcone} 
For a smooth representation $(\pi, \cH)$ of $G$, 
the following assertions hold: 
\begin{description}
\item[\rm(i)]  If $C_\pi\not=\eset$, then $\pi$ is semibounded, and if 
$i\1 \in \dd\pi(\g)$, then the converse is also true. 
\item[\rm(ii)]  If $C_\pi\not=\eset$, then $H(C_\pi) = \ker \dd\pi$. 
In particular, $C_\pi$ is pointed if 
and only if the derived representation $\dd \pi$ is injective.  
\item[\rm(iii)] If $\fh \subeq \g$ is a subalgebra with the property 
that $\dd\pi\res_{\fh}$ is a bounded representation, then 
$\fh \subeq H(W_\pi)$. 
\item[\rm(iv)] Let $\eta \: H \to G$ be a morphism of Lie subgroups. 
If $\pi$ is semibounded, then 
$\pi_H := \pi \circ \eta$ is semibounded if $W_\pi \cap \L(\eta)\fh 
\not=\eset$, and then $I_{\pi_H} = \L(\eta)^*I_\pi$. 
In particular, if $H \subeq G$ is a Lie subgroup, 
then the restriction $\pi_H := \pi\res_H$ of $\pi$ 
to $H$ is semibounded 
if $W_\pi \cap \fh \not=\eset$, and then $I_{\pi_H} = I_\pi\res_{\fh}$. 
\end{description}
\end{prop}

\begin{prf} (i) If $C_\pi \not=\eset$, then the boundedness of 
$s_\pi$ on $C_\pi$ implies that $\pi$ is semibounded. 
Suppose, conversely, that $\pi$ is semibounded and 
$i\1 \in \dd\pi(\g)$. Then 
there exists an element $x \in \g$ for which 
$s_\pi$ is bounded from above on a neighborhood of $x$ by some $M \in \R$. 
Pick $z \in \g$ with $\dd\pi(z) = i \1$. 
Now $s_\pi(x + M z) \leq 0$ on this neighborhood, so that 
$C_\pi\not=\eset$. 

(ii) The relation $C_\pi = (I_\pi^\star)^0$ implies that 
\begin{equation}
  \label{eq:c-dual}
  I_\pi^\star = \oline{C_\pi}
\end{equation}
(Lemma~\ref{lem:bou}). If $C_\pi$ is non-empty, 
then Lemma~\ref{lem:limcone}(i) implies that 
$$H(C_\pi) = H(I_\pi^\star) = I_\pi^\star \cap - I_\pi^\star 
= I_\pi^\bot = \ker \dd\pi $$ 
because $\dd\pi(x) = 0$ is equivalent to 
$\Spec(i\dd\pi(x))) \subeq \{0\}$. 

(iii) That $\dd\pi\res_{\fh}$ is bounded implies that 
all operators $i\dd\pi(x)$, $x \in \fh$, are bounded, so that 
$\fh \subeq B(I_\pi) \subeq \oline{W_\pi} = \lim(W_\pi)$ 
(Lemma~\ref{lem:limcone}(i)), and since $\fh = - \fh$, 
the assertion follows. 

(iv) From \eqref{eq:momspec} 
it follows that $s_{\pi_H} = s_\pi \circ \L(\gamma)$, 
and if $W_\pi \cap \L(\gamma)\fh \not=\eset$, this function 
is bounded on some non-empty open subset of $\fh$. Therefore 
$\pi_H$ is semibounded. 

Pick $x \in \fh$ with $y := \L(\gamma)x \in W_\pi$. Then the evaluation map 
$\eta_y^\g \: I_{\pi} \to \R, \alpha \mapsto \alpha(x),$ 
is proper by Proposition~\ref{prop:locomp}, 
and this map factors through the adjoint map 
$\L(\eta)' \: \g' \to \fh', \alpha \mapsto\alpha \circ \L(\eta)$: 
$$ \eta_y^\g = \eta_x^\fh \circ \L(\eta)' \: I_\pi \to \R. $$
Since the weak-$*$-topology on $\fh'$ is Hausdorff, 
\cite[Ch.~I, \S 10, Prop.~1.5]{Bou89} implies that 
$\L(\eta)' \: I_{\pi} \to \fh'$ is a proper map. In particular, 
its image is closed, hence a weak-$*$-closed convex subset 
of $\fh'$. For each $x \in \fh$ we have 
$$ \sup \la I_{\pi_H}, x \ra 
= s_{\pi_H}(-x)= s_{\pi}(-\L(\eta)x) = \sup \la I_\pi, \L(\eta)x \ra 
= \sup \la \L(\eta)'I_{\pi}, x\ra, $$
so that the Hahn--Banach Separation Theorem implies that 
$I_{\pi_H}  = \L(\eta)'I_\pi$. 
\end{prf}

\begin{rem} \mlabel{rem:4.7} 
If $i\1 \not\in \dd\pi(\g)$, then we may consider the 
direct product Lie group $\hat G := G \times \T$, 
where we consider $\T$ as a subgroup of $\C^\times$, with the exponential 
function $\exp_\T(t) = e^{it}$ on $\L(\T) \cong \R$. 
Our representation $\pi$ now extends trivially to a smooth 
unitary representation of $\hat G$, defined by 
$$ \hat\pi(g,t) := t \pi(g). $$
Now $\hat\g \cong \g \times \R$, $\hat\g' \cong \g' \times\R$, 
and the momentum map of $\hat\pi$ is given by 
$$ \Phi_{\hat\pi}([v]) = (\Phi_\pi([v]), 1), $$ 
so that $I_{\hat\pi} = I_{\pi} \times \{1\}$. 
With Remark~\ref{rem:conered}, we now see that 
$\pi$ is semibounded if and only if $C_{\hat\pi} \not=\eset$. 
Note also that $\dd\hat\pi(0,1) = i \1$, so that 
Proposition~\ref{prop:repcone}(i) applies to this representation. 
\end{rem}

\begin{rem} (on bounded representations) (a) 
From \cite[Thm.~3.1]{Ne09a} we know that a smooth representation 
$(\pi, \cH)$ is bounded if and only if $\pi \: G \to \U(\cH)$ 
is a morphism of Lie groups, if $\U(\cH)$ is endowed with the 
Banach--Lie group structure defined by the operator norm. 
In view of \eqref{eq:momspec}, the boundedness of $(\pi, \cH)$ 
is equivalent to the boundedness of all operators 
$i\dd\pi(x)$ and the continuity of the derived representation 
$\dd\pi \: \g \to B(\cH)$ as a morphism of topological 
Lie algebras, where we identify the operator $\dd\pi(x)$ with 
its continuous extension to all of $\cH$. 

(b) In \cite[Prop.~3.5]{Ne09a} it is also shown that if 
$\pi$ is continuous with respect to the norm topology on 
$\U(\cH)$, then $\pi$ is automatically smooth if 
either $G$ is locally exponential or $\g$ is a barrelled space. 
In particular, for Banach--Lie groups $G$ the bounded representations 
are precisely the norm-continuous ones. The concept of a norm-continuous 
unitary representation also makes sense for arbitrary topological 
groups, but the refined concept of semiboundedness does not; it 
requires the Lie algebra for its definition. 

(c) If $(\pi, \cH)$ is a bounded unitary representation of a 
Lie group $G$, then the closure $\fh$ of $\dd\pi(\g) \subeq \fu(\cH)$ 
is a Banach--Lie algebra, and $H := \la \exp_{\UU(\cH)} \fh \ra \subeq 
\U(\cH)$ carries a natural Banach--Lie group structure with 
$\L(H) = \fh$ (\cite[Thm.~IV.4.9]{Ne06}). 
If $G$ is connected, then $\pi \: G \to \U(\cH)$ 
factors through the inclusion map $H \to \U(\cH)$ of this 
Banach--Lie group. In this sense bounded representations are 
a ``Banach phenomenon'' and all questions on this class of 
representations can be reduced to Banach--Lie groups. 

(d) If $(\pi, \cH)$ is a bounded representation for which 
$\ker(\dd\pi) = I_\pi^\bot = \{0\}$, then $\|\dd\pi(x)\|$ 
defines a $G$-invariant norm on $\g$. If $G$ is finite-dimensional, 
then the existence of an invariant norm implies that the 
Lie algebra $\g$ is compact. In particular, all its 
irreducible representations are finite dimensional. 
However, in the infinite dimensional context, there is a substantially 
richer supply of bounded unitary representations. For a detailed 
discussion of bounded unitary representations of the unitary 
groups $\U_p(\cH)$ (Examples~\ref{ex:liegrp}(b)), we refer to \cite{Ne98}. 
\end{rem}

\glossary{name={$A \< B$},description={$\inf\Spec(B-A)> 0$, for $A,B \in \Herm(\cH)$}}
\begin{ex} \mlabel{ex:4.9} 
Let $\cH$ be a complex Hilbert space and consider 
$\U(\cH)$ as a Banach--Lie group with Lie algebra $\fu(\cH)$. 
Writing $A \< B$ for $A, B \in B(\cH)$ if $B - A$ is a positive 
invertible operator, we see that for 
the identical representation $(\pi, \cH)$, the cone 
$$ C_\pi = \{ x \in \fu(\cH) \: i x \< 0\} $$
is non-empty. Since $\pi$ is bounded, we also have 
$W_\pi = \g$. The same holds for all representations 
$(\pi^{\otimes n}, \cH^{\otimes n})$ on the $n$-fold tensor 
power of $\cH$. 

The natural representations $\pi_s$ and $\pi_a$ of $\U(\cH)$ 
on the symmetric and antisymmetric 
Fock spaces  
\glossary{name={$S(\cH)$},description={bosonic/symmetric Fock space}} 
\glossary{name={$\Lambda(\cH)$},description={fermionic/alternating 
Fock space}} 
$$ S(\cH) := \hat\oplus_{n \in \N_0} S^n(\cH) \quad \mbox{ and } \quad 
\Lambda(\cH) := \hat\oplus_{n \in \N_0} \Lambda^n(\cH) $$ 
are direct sums of bounded representations, hence in particular smooth. 
They are not bounded, as the restriction to the subgroup 
$\T \1$ already shows. However, the relation 
$C_\pi \subeq C_{\pi_s}, C_{\pi_a}$ still implies that 
$\pi_s$ and $\pi_a$ are semibounded. We shall use this simple 
observation later in our proof that the metaplectic and the 
spin representation are semibounded (cf.\ Sections~\ref{sec:9}, \ref{sec:10}). 
\end{ex}

We conclude this section with a convenient tool to verify 
the existence of eigenvectors for semibounded 
unitary one-parameter groups. 

\begin{prop} \mlabel{prop:maxsmoothvect} 
Let $(\pi, \cH)$ be a
semibounded unitary representation of $G$ 
and $d \in W_\pi$ with $\exp_G(\lambda d) \in Z(G)$ for some 
$\lambda \in \R^\times$. 
Suppose that $\pi(\exp_G(\lambda d)) \in \T \1$, which is in 
particular the case if $\pi$ is irreducible. 
Then $i\dd\pi(d)$ is bounded from above, has discrete spectrum, 
and there exists a smooth vector $v \in\cH^\infty$ which is an 
eigenvector for the largest eigenvalue of $i\dd\pi(d)$. 
\end{prop}

\begin{prf} Replacing $d$ by $\lambda d$, we may assume that 
$\lambda = 1$. If $\pi$ is irreducible, Schur's Lemma 
implies that $\pi(Z(G)) \subeq \T \1$ and thus in particular 
$\pi(\exp_G d) \in \T\1$. Let 
$\beta \in \R$ with $\pi(\exp d) = e^{i\beta}\1$. 
Then $e^{\dd\pi(d)- i \beta\1} = \1$, 
so that \break $\Spec(i\dd\pi(d) + \beta \1) \subeq 2 \pi \Z$ 
is discrete. 

Let $P_n$ denote the orthogonal projection 
onto the eigenspace of $i\dd\pi(d)$ corresponding to the eigenvalue 
$2\pi n - \beta$. This projection is given by the following integral 
$$ P_n(v) 
= \int_0^1 e^{2\pi i nt} e^{-it \beta}\pi(\exp(td))v\, dt 
= \int_0^1 e^{it(2\pi n- \beta)}\pi(\exp(td))v\, dt. $$
If $v$ is a smooth vector for $G$, then the smoothness of the 
$\cH$-valued function 
$$ g\mapsto \pi(g) P_n(v) 
= \int_0^1 e^{it(2\pi n- \beta)}\pi\big(g\exp(td)\big)v\, dt $$
follows from the smoothness of the integrand as a function 
on $\R \times G$ by differentiation under the integral 
(cf.~\cite{GN09} for details). This implies that 
$P_n(\cH^\infty) \subeq \cH^\infty$. 

Therefore the density of $\cH^\infty$ in $\cH$ implies the 
density of $P_n(\cH^\infty)$ in $P_n(\cH)$. 
In particular, each non-zero eigenspace of $\dd\pi(d)$ contains 
smooth vectors. Finally, the boundedness of $\Spec(i\dd\pi(d))$ 
from above implies the existence of a maximal eigenvalue. 
\end{prf}

\section{Aspects of complex analysis} \mlabel{sec:5}

As we have already seen in \cite{Ne09a}, a closer analysis of 
semibounded representations requires a good understanding of the 
related complex geometric structures. Basically, the connection 
to complex analysis relies on the fact that if 
$A$ is a selfadjoint operator on a complex Hilbert space, 
then the associated unitary one-parameter group 
$\gamma_A(t) := e^{itA}$ extends to a strongly continuous homomorphism 
$\hat\gamma \: \C_+ := \R + i \R_+ \to B(\cH), z \mapsto e^{izA}$ 
\glossary{name={$\C_+$},description={closed upper half plane}}
(defined by a spectral integral) which 
is holomorphic on the open upper half plane if and only if 
$\Spec(A)$ is bounded from below. It is this key observation, 
applied in various situations, that leads to a variety of 
interesting tools and results, that we describe below. 

\begin{prob} If $(\pi, \cH)$ is a semibounded unitary representation, 
then we have a well-defined map 
$$ \hat\pi \: G \times W_\pi \to B(\cH), \quad 
(g,x) \mapsto g e^{i\dd\pi(x)}, $$
and if $G$ is finite dimensional and $\dd\pi$ injective, then 
we know from \cite{Ne00} that the product set 
$S_\pi := G \times W_\pi$ always carries a complex manifold structure, 
a holomorphic associative multiplication and an antiholomorphic involution 
$(g,x)^* := (g^{-1}, \Ad(g)x)$, turning $(S_\pi, *)$ into a complex 
involutive semigroup and $\hat\pi$ into a holomorphic representation. 
In particular, $\hat\pi(G \times W_\pi) = \pi(G) e^{i\dd\pi(W_\pi)}$ 
is an involutive subsemigroup of $B(\cH)$. 
This structure is extremely useful in the 
theory of semibounded representations and the related 
geometry. However, for infinite dimensional groups, our 
understanding of corresponding analogs is still quite rudimentary. 
\end{prob} 

As we shall see below, the preceding problem has a trivial solution 
for abelian groups.

\subsection{The abelian case} 

If $E$ is a locally convex space, then 
its additive group $G := (E,+)$ is a particularly 
simple locally convex Lie group. 
Then $\L(G) = E$ as an abelian Lie algebra, 
$\exp_G = \id_E$, and the coadjoint action is trivial. 

For such groups semibounded representations can be 
understood completely with classical methods, well known from 
the context of locally compact abelian groups. 
If $X \subeq E'$ is a weak-$*$-closed convex semi-equicontinuous 
subset, then Proposition~\ref{prop:locomp} implies that it is 
locally compact, so that the space 
$C_0(X)$ of continuous complex-valued functions on $X$ 
vanishing at infinity 
is a commutative $C^*$-algebra. 
Moreover, the interior $B(X)^0$ is a non-empty 
open convex cone in $X$, so that 
$S := E + i B(X)^0$ is a tube domain in the complexification 
$E_\C$. This open subset of $E_\C$ 
is a complex manifold, a semigroup with respect to 
addition, and $(x + iy)^* := -x + i y$ defines an antiholomorphic 
involution, turning it into a {\it complex involutive semigroup} 
(cf.~\cite{Ne08} for an extended discussion of this technique). 
This leads to a holomorphic homomorphism of semigroups 
$$\gamma \: S \to C_0(X), \quad 
\gamma(x + i y)(\alpha) := e^{i\alpha(x)-\alpha(y)} $$
(cf.\ Proposition~\ref{prop:locomp}), and one can show that 
every semibounded smooth representation $(\pi, \cH)$ of 
$G$ with $I_\pi \subeq X$ ``extends'' holomorphically 
to $S$ by \break 
$\hat\pi(x + iy) := \pi(x)e^{i\dd\pi(y)}$, which in 
turn extends to a representation $\rho$ of the $C^*$-algebra 
$C_0(X)$ on $\cH$. Using this correspondence and the Spectral Theorem 
for $C^*$-algebras, one arrives at the following 
spectral theorem (\cite[Thm.~4.1]{Ne09a}): 

\begin{thm} \mlabel{thm:ab} {\rm(Spectral Theorem 
for semibounded Representations)} 
For every regular Borel spectral measure 
$P$ on $X$, the prescription 
$\pi(v) := P(e^{i\eta_v})$
(where the right hand side denotes a spectral integral) 
defines a semibounded smooth representation of $G$ with 
$I_\pi \subeq X$. Conversely, any such semibounded representation has this 
form for a uniquely determined regular Borel spectral measure 
$P$ on $X$. 
\end{thm}

The preceding theorem provides a complete description of the 
semibounded representation theory of abelian Lie groups 
of the form $(E,+)$ and hence also of quotients thereof. 
In particular, every connected abelian Banach--Lie group is such a quotient. 
For recent results concerning extremely general spectral theorems 
for representations of commutative involutive algebras by unbounded 
operators  we refer to \cite{Th09}. 

\begin{rem} With Theorem~\ref{thm:ab} it is easy to see that 
the bounded unitary representations of $(E,+)$ are precisely those 
defined by spectral measures on (compact) equicontinuous 
subsets of the dual space. 

Similar characterizations are known for continuous isometric 
actions $\alpha \: G \to \Iso(E)$ 
of a locally compact abelian group $G$ on a Banach space $E$. 
One associates to such a representation its Arveson spectrum 
$\Sp(\alpha) \subeq \hat G$ and then the norm continuity 
of $\alpha$ is equivalent to the compactness of its 
spectrum (\cite[Thm.~8.1.12]{Ped89}).
\end{rem}

\begin{rem} The key point behind the Spectral Theorem~\ref{thm:ab} 
is that the semiboundedness provides a method to connect 
representations of $(E,+)$ to representations of a commutative 
$C^*$-algebra, so that the Spectral Theorem for commutative 
Banach-$*$-algebras provides the spectral measure on $X \subeq E'$. 

In general, continuous unitary representations of 
locally convex spaces can not be represented in terms of spectral 
measures on $E'$. This is closely related to the problem of writing 
the continuous positive definite functions 
$\pi^{v,v}(x) := \la \pi(x)v,v\ra$ as the Fourier transform  
$$ \hat\mu(x) = \int_{E'} e^{i\alpha(x)}\, d\mu(\alpha) $$
of some finite measure $\mu$ on $E'$. If $E$ is nuclear, 
then the Bochner--Minlos Theorem (\cite{GV64}) ensures the existence 
of such measures and hence also of spectral measures representing 
unitary representations. However, if $E$ is an infinite dimensional 
Banach space, then $E$ is not nuclear (\cite{GV64}), 
so that the Bochner--Minlos 
Theorem does not apply. Therefore it is quite remarkable that 
nuclearity assumptions are not needed to deal with semibounded 
representations. 
\end{rem}

\subsection{Hilbert spaces of holomorphic functions} 

In general, the momentum set of a semibounded representation is 
not easy to compute, but in many interesting 
situations it is the weak-$*$-closed convex hull 
of a single coadjoint orbit $\cO_\lambda \subeq \g'$. 
Here our intuition is guided by the finite dimensional case, 
which is by now well understood (cf.\ \cite{Ne00}).

\begin{rem}
  \mlabel{rem:findim} 
For a finite dimensional connected Lie group $G$,  
semibounded unitary representations are direct integrals 
of irreducible ones (\cite[Sect.~XI.6]{Ne00}) and the irreducible 
ones possess various kinds of nice structures. Here the key 
result is that for every irreducible semibounded 
representation $(\pi, \cH)$  and $x \in B(I_\pi)^0$, 
we can minimize the proper functional 
$\eta_x(\alpha)= \alpha(x)$ on the convex set $I_\pi$ 
(Proposition~\ref{prop:locomp}). One can even show that 
the minimal value is taken in $\Phi_\pi([v])$ for an analytic 
vector $v$, and from that one can derive that $v$ is an eigenvector 
for $i\dd\pi(x)$ (this is implicitly shown in 
\cite[Thm.~X.3.8]{Ne00}). 
Combing this with finite 
dimensional structure theory, based on the observation 
that the adjoint image of the centralizer $Z_G(x)$ has 
compact closure, one can show that if 
$\dd\pi$ is injective, then $\cH^\infty$ contains a dense 
highest weight module of the complexification 
$\g_\C$ (\cite[Thms.~X.3.9, XI.4.5]{Ne00}).   

If $[v_\lambda] \in \bP(\cH^\infty)$ is a highest weight vector 
for the highest weight representation 
$(\pi_\lambda, \cH_\lambda)$, then 
the corresponding $G$-orbit $G[v_\lambda]$ has the remarkable 
property that it is a complex homogeneous subspace of $\bP(\cH^\infty)$, 
and one can even show that it is the unique $G$-orbit in 
$\bP(\cH^\infty)$ with this property 
(\cite[Thm.~XV.2.11]{Ne00}). Its image 
$\cO_\lambda := \Phi_{\pi_\lambda}(G[v_\lambda]) \subeq \g^*$ is a coadjoint orbit 
satisfying \glossary{name={$\Ext(C)$},description={extreme points of convex set $C$}}
\begin{equation} 
  \label{eq:extorbit}
\cO_\lambda = \Ext(I_{\pi_\lambda}), \quad 
I_{\pi_\lambda} = \conv(\cO_\lambda) 
\quad \mbox{ and } \quad 
G[v_\lambda] = \Phi_{\pi_\lambda}^{-1}(\cO_\lambda) 
\end{equation}
(\cite[Thm.~X.4.1]{Ne00}), where we write 
$\Ext(C)$ for the set of extreme points of a convex set~$C$. 
Moreover, two irreducible semibounded representations are equivalent if 
and only if the corresponding momentum sets, resp., the coadjoint orbits 
$\cO_\lambda$ coincide (\cite[Thm.~X.4.2]{Ne00}).  
\end{rem}

In particular, a central feature of unitary highest weight representations 
is that the image of the highest weight orbit already determines the 
momentum set as its closed convex hull. 
It is therefore desirable to understand in which situations 
certain $G$-orbits in $\bP(\cH^\infty)$ 
already determine the momentum set as the closed convex hull 
of their image. 
As we shall see below, this situation frequently occurs if
$\cH$ consists of holomorphic functions on some complex manifold, 
resp., holomorphic sections of a line bundle. 

\begin{defn} Let $M$ be a complex manifold (modelled on a locally 
convex space) and $\cO(M)$ the space of holomorphic complex-valued 
functions on $M$. We write $\oline M$ for the conjugate complex manifold. 
A holomorphic function 
$$ K \: M \times \oline M \to \C $$
is said to be a {\it reproducing kernel} of a Hilbert subspace 
$\cH \subeq \cO(M)$ if for each $w\in M$ the function 
$K_w(z) := K(z,w)$ is contained in $\cH$ and satisfies 
$$ \la f, K_z \ra = f(z) \quad \mbox{ for } \quad z \in M, f \in \cH. $$
Then $\cH$ is called a {\it reproducing kernel Hilbert space} 
and since it is determined uniquely by the kernel $K$, it is 
denoted $\cH_K$ (cf.\ \cite[Sect.~I.1]{Ne00}). 
\glossary{name={$\cH_K$},description={Hilbert space with reproducing kernel $K$}}
\glossary{name={$\cO(M)$},description={holomorphic functions on $M$}}

Now let $G$ be a real Lie group and 
$\sigma \: G \times M  \to M, (g,m) \mapsto g.m$ 
be a smooth right action of $G$ on $M$ by holomorphic maps. Then 
 $(g.f)(m) := f(g^{-1}.m)$ defines a unitary representation 
of $G$ on a reproducing kernel Hilbert space $\cH_K \subeq \cO(M)$ 
if and only if the kernel $K$ is {\it invariant}: 
$$ K(g.z, g.w) = K(z,w) \quad \mbox{ for } \quad z,w \in M, g \in G. $$
In this case we call 
${\cal H}_K$ a $G$-invariant reproducing kernel Hilbert space 
and write $(\pi_K(g)f)(z) := f(g^{-1}.z)$ for the corresponding 
unitary representation of $G$ on~$\cH_K$. 
\end{defn} 

\begin{thm} {\rm(\cite[Thm.~2.7]{Ne09a})} \mlabel{thm:2.7} 
Let $G$ be a Fr\'echet--Lie group acting smooth\-ly 
by holomorphic maps on the complex manifold $M$ and 
$\cH_K \subeq \cO(M)$ be a $G$-invariant reproducing kernel 
Hilbert space on which the representation is semibounded, so that 
the cone $W_{\pi_K}$ is not empty.  
If, for each $x \in W_{\pi_K}$,  
the action $(m,t) \mapsto \exp_G(-tx).m$ 
of $\R$ on $M$ extends to a holomorphic action of the upper half plane 
$\C_+$ (continuous on $\C_+$ and holomorphic on its interior), then 
$$ I_{\pi_K} = \oline{\conv}(\Phi_{\pi_K}(\{ [K_m] \: K_m \not=0\})). $$
\end{thm}

\begin{ex} \mlabel{ex:5.7} Here is the prototypical example that 
shows why $x \in W_{\pi_K}$ is closely related to the existence of a 
holomorphic extension of the corresponding flow on $M$. 

Let $(\pi, \cH)$ be a continuous 
unitary representation of $G = \R$  and 
$\cH^\infty$ be the Fr\'echet space of smooth vectors on 
which $G$ acts smoothly (Theorem~\ref{thm:4.1.3}). 
Let $M := \oline{\cH^\infty}$ be the complex Fr\'echet manifold  
obtained by endowing $\cH^\infty$ with the opposite complex 
structure. Then $G$ acts on $M$ by holomorphic maps 
and the density of $\cH^\infty$ in $\cH$ yields an embedding 
$$\iota \: \cH \into \cO(M), \quad \iota(v)(m) = \la v,m\ra, $$
whose image is the reproducing kernel space $\cH_K$ with 
kernel $K(x,y) =\la y,x \ra.$ 
With $A := -i\dd\pi(1)$ we then have $\pi(t) = e^{itA}$, 
and if $\Spec(A)$ is bounded from below, then 
$\hat\pi \: \C_+ \to B(\cH), z \mapsto e^{zA}$ defines a holomorphic 
extension of the unitary representation to $\C_+$ 
(cf.\ \cite[Thm.~VI.5.3]{Ne00}). 
Since $\cH^\infty$ is invariant under every operator 
$\hat\pi(z)$, it is easy to see that 
$(z,m) \mapsto e^{-\oline z A}m$ defines a holomorphic 
action of $\C_+$ on $M$, extending the action of 
$\R$ given by $(t,m) \mapsto \pi(-t)m$.  
\end{ex}


\begin{rem} \mlabel{rem:5.8} Suppose that 
$q \: \bV \to M$ is a $G$-homogeneous holomorphic Hilbert bundle, 
i.e., $\bV$ carries a left action of $G$ by holomorphic 
bundle automorphisms and $q(g.z) = g.q(z)$ for $z \in \bV$, 
$g \in G$. We write $\Gamma(\bV)$ for 
the \glossary{name={$\Gamma(\bV)$},description={holomorphic sections of vector bundle $\bV$}}
space of holomorphic section of $\bV$. Then $G$ 
acts on $\Gamma(\bV)$ by $(g.s)(m) :=g.s(g^{-1}.m)$. 
We are interested in $G$-invariant Hilbert subspaces 
$\cH \subeq \Gamma(\bV)$ on which $G$ acts unitarily and for 
which the evaluation maps $\ev_m \: \cH \to \bV_m, s \mapsto s(m)$, 
$m \in M$, 
are continuous linear maps between Hilbert spaces. 

To see how Theorem~\ref{thm:2.7} applies in this situation, 
we realize $\Gamma(\bV)$ by holomorphic functions on the dual 
bundle $\bV^*$ whose fiber $(\bV^*)_m$ is the dual space  
of $\bV_m$. Each $s \in \Gamma(\bV)$ defines a holomorphic 
function $\hat s(\alpha_m) := \alpha_m(s(m))$ on $\bV^*$ which 
is fiberwise linear. We thus obtain an embedding 
$\Psi \: \Gamma(\bV) \to \cO(\bV^*)$ whose image consists of 
those holomorphic functions on $\bV^*$ which are fiberwise linear. 
Accordingly, $\Psi(\cH) \subeq\cO(\bV^*)$ is a reproducing kernel 
Hilbert space. The natural action of $G$ on $\bV^*$ is given by 
$(g.\alpha_{m})(z_{g.m}) := \alpha_{m}(g^{-1}.z_{g.m})$ for 
$\alpha_m \in \bV^*_{m}$, so that 
$$ \Psi(g.s)(\alpha_m) 
= \alpha_m(g.s(g^{-1}.m)) 
= (g^{-1}.\alpha_m)(s(g^{-1}.m)) 
= \Psi(s)(g^{-1}.\alpha_m) $$
implies that $\Psi$ is equivariant with respect to the natural 
$G$-actions on $\Gamma(\bV)$ and $\cO(\bV^*)$. 
Therefore the reproducing kernel $K$ of $\cH_K := \Psi(\cH) 
\subeq \cO(\bV^*)$ is 
$G$-invariant, and we are thus in the situation of 
Theorem~\ref{thm:2.7}. In addition, the fiberwise linearity of 
the functions in $\cH_K$ leads to 
$K_{z\alpha} = \oline z K_\alpha$ for $\alpha \in \bV^*$ and 
$z \in \C^\times$. 
If $\cH_K \not=\{0\}$, then the homogeneity of 
the bundle $\bV$ implies that $K_\alpha\not=0$ for every 
$0\not= \alpha \in \bV^*$. 
Writing $\bP(\bV^*)$ for the projective bundle associated to 
$\bV^*$ whose fibers are the projective spaces of the fibers of 
$\bV^*$, we therefore obtain a well-defined map 
$$ \Phi \: \bP(\bV^*) \to \g', 
\quad [\alpha_m] \mapsto \Phi_{\pi_K}([K_{\alpha_m}])\quad \mbox{ for } 
\quad \alpha_m \in \bV^*_m \setminus \{0\}. $$
Since this map is $G$-equivariant and $G$ acts transitively 
on $M$, we obtain for each $m_0 \in M$: 
\begin{equation}
  \label{eq:imom1}
\Phi(\bP(\bV^*)) = \Ad^*(G)\Phi(\bP(\bV^*_{m_0})).
\end{equation}

If, in addition, the requirements of Theorem~\ref{thm:2.7} are 
satisfied, i.e., the action of the one-parameter-groups 
generated by $-x \in B(I_{\pi_K})^0$ on $\bV$, resp., $\bV^*$, 
 extend holomorphically to $\C_+$, we obtain 
 \begin{equation}
   \label{eq:imom2}
I_{\pi_K} = \oline\conv(\im(\Phi)). 
 \end{equation}

If, in particular, $\bV$ is a line bundle, i.e., the fibers are one-dimensional, then $\bP(\bV^*) \cong M$, and we obtain a $G$-equivariant map 
$\Phi \: M \to \g'$ whose image is a single coadjoint orbit~$\cO_\pi$. 
\end{rem}

\begin{ex} \mlabel{ex:5.9} If $\cH$ is a complex Hilbert space, 
then its projective space $\bP(\cH)$ is a complex Hilbert manifold. 
Moreover, there exists a holomorphic line bundle 
$q \: \bL_\cH \to \bP(\cH)$ with the property that for 
every non-zero continuous linear functional 
$\alpha \in \cH'$ we have on the open subset
$U_\alpha := \{ [v] \in \bP(\cH) \: \alpha(v) \not=0\}$ 
a bundle chart 
$$ \phi_\alpha \: (\bL_\cH)\res_{U_\alpha} \to U_\alpha \times \C $$
such that the transition functions are given by 
$$ \phi_\beta \circ \phi_\alpha^{-1}([v],z) 
= \Big([v], \frac{\alpha(v)}{\beta(v)}\Big)\quad \mbox{ for } \quad 
0\not=\alpha, \beta \in \cH'. $$ 
This implies that each $0\not=v \in \cH$ defines a linear functional 
on the fiber $(\bL_\cH)_{[v]}$ by 
$$ \phi_\alpha^{-1}([v], z) \mapsto \alpha(v) z, $$
which further implies that $(\bL_\cH)_{[v]}^* = [v]$, 
i.e., $\bL_\cH^*$ is the tautological bundle over $\bP(\cH)$. 

The complement of the zero-section 
of $\bL_\cH$ is equivalent, as a $\C^\times$-bundle,  
to the projection $\cH\setminus \{0\} \to \bP(\cH)$ by 
the map $\phi_\alpha([v],z) \mapsto \frac{1}{z \alpha(v)}v$. 
This identification can be used to show that the natural map 
$$ \Psi \: \cH' \to \Gamma(\bL_\cH), \quad 
\Psi(\alpha)([v]) = \phi_\beta^{-1}\Big([v], \frac{\alpha(v)}{\beta(v)}\Big) 
\quad \mbox{for } \quad 
\beta(v)\not=0 $$
defines a linear isomorphism (see \cite[Thm.~V.4]{Ne01b} for details). 

As the group $\U(\cH)$ acts smoothly by holomorphic bundle 
isomorphisms on $\bL_\cH$, this construction shows that the 
unitary representation $\pi^* \: \U(\cH) \to \U(\cH')$, 
given by $\pi^*(g)\alpha = \alpha \circ \pi(g)^*$ 
can be realized in the space $\Gamma(\bL_\cH)$ of holomorphic 
sections of $\bL_\cH$. 

To realize the identical representation on $\cH$ itself 
by holomorphic sections, we simply exchange the role 
of $\cH$ and $\cH'$, which leads to a $\U(\cH)$-equivariant isomorphism 
$\cH \to \Gamma(\bL_{\cH'})$. 
\end{ex}

To simplify the applications of Theorem~\ref{thm:2.7}, 
we need a criterion for its applicability. Here the main 
idea is that, since  every Hilbert space $\cH$ can be realized as a 
space of holomorphic sections of the bundle $\bL_{\cH'}$, we obtain 
similar realizations from cyclic $G$-orbits in $\bP(\cH)$ which 
are complex manifolds. 

\begin{thm} {\rm(Complex Realization Theorem)} \mlabel{thm:5.7}
{\rm(a)} Let $G$ be a Fr\'echet--Lie group with Lie algebra $\g$ and 
$H \subeq G$ be a closed subgroup for which the coset space 
$G/H$ carries the structure of a complex manifold such that 
the projection $q \: G\to G/H$ is a smooth $H$-principal bundle and 
$G$ acts on $G/H$ by holomorphic maps. 
Let $x_0 = \1 H \in G/H$ be the canonical base point and 
$\fp \subeq \g_\C$ be the kernel of the complex linear extension 
of the map $\g \to T_{x_0}(G/H)$ to $\g_\C$, so that $\fp$ is a closed 
subalgebra of~$\g_\C$.  

Let $(\pi, \cH)$ be a unitary representation 
of $G$ and $\dd \pi \: \g_\C \to \End(\cH^\infty)$ be the complex 
linear extension of the derived representation. 
Suppose that $0\not=v \in \cH^\infty$ is 
an eigenvector for $H$ and of the subalgebra 
$\oline\fp := \{ \oline{x + iy} = x - i y \: x + i y\in \fp\}$ of $\g_\C$. 
Then the map 
$$ \eta \: G/H \to \bP(\cH'), 
\quad gH \mapsto [\pi^*(g)\alpha_v] = [\alpha_v \circ \pi(g)^{-1}], 
\quad \alpha_v(w) = \la w, v\ra, $$ 
is holomorphic and $G$-equivariant. If, in addition, 
$v$ is cyclic, then we obtain a 
$G$-equivariant injection $\cH \into \Gamma(\eta^*\bL_{\cH'})$, 
where $\eta^*\bL_{\cH'}$ is a $G$-equivariant holomorphic 
line bundle over $G/H$. If, in addition, $G$ is connected, 
then $\pi$ is irreducible. 

{\rm(c)} Suppose that, for each $x \in W_\pi$, the flow 
$(t,m) \mapsto \exp_G(-tx)H$ 
extends holomorphically to $\C_+$. Then the momentum set is given by 
$$ I_\pi 
= \oline{\conv}(\Phi_\pi(G[v])) = \oline{\conv}(\cO_{\Phi_\pi([v])}). $$
\end{thm}

\begin{prf} 
(a) First we observe that the map 
$\Upsilon \: \cH \to \cH', v \mapsto \alpha_v$ is an antilinear 
isometry and that the contragredient representation 
$\pi^*(g)\alpha := \alpha \circ \pi(g)^{-1}$ is also unitary. 

Next we recall that the smooth action of $G$ on 
$M := G/H$ defines a homomorphism 
$\g \to \cV_{\cal O}(M)$, the Lie algebra of holomorphic vector fields. 
Since $M$ is complex, $\cV_{\cal O}(M)$ is a complex 
Lie algebra, and, for each $p \in M$, the subspace 
$\{ X \in \cV_{\cal O}(M) \: X(p) = 0\}$ is a complex subalgebra. 
This proves that $\fp$ is a Lie subalgebra of $\g_\C$, and its closedness 
follows from the continuity of the map $\g_\C \to T_{x_0}(M) \cong \g/\fh$. 

(b) Since $v$ is a $\oline\fp$-eigenvector, there exists a continuous 
linear functional $\lambda \: \fp \to \C$ with 
$\dd \pi(\oline z) v = \oline{\lambda(z)} v$ for 
$z \in \fp$, so that we have for 
$w \in \cH$ and $z \in \fp$ the relation 
$$ \alpha_v(z.w) 
= \la \dd\pi(z)w, v \ra 
= - \la w, \dd\pi(\oline z)v \ra 
= - \lambda(z) \alpha_v(w). $$
We conclude that $\alpha_v$ is an $\fp$-eigenvector. This implies that 
the tangent map 
$$ T_{x_0}(\eta) \: T_{x_0}(M) \to T_{[\alpha_v]}(\bP(\cH')) $$
is complex linear, i.e., compatible with the respective complex 
structures. Since $M$ is homogeneous, $T(\eta)$ is complex 
linear on each tangent space, and this means that $\eta$ is holomorphic. 
Therefore $\eta^*\bL_{\cH'}$ is a holomorphic line bundle over 
$M$ and we obtain a $G$-equivariant pullback map 
$\cH \cong \Gamma(\bL_{\cH'}) \to \Gamma(\eta^*\bL_{\cH'})$. 
Its kernel consists of all those sections 
vanishing on $\eta(M)$, which corresponds to the elements 
$w \in (\pi(G)v)^\bot$. In particular, this map is injective if 
$v$ is a cyclic vector. 

That $(\pi, \cH)$ is irreducible if $G$ is connected follows from 
\cite[Prop.~XV.2.7]{Ne00}. 

(c) If $x \in W_\pi$, then $i\dd\pi(x)$ is bounded from above, 
so that, so that  $(z,\alpha) \mapsto \alpha \circ e^{z \dd\pi(x)}$ 
defines a continuous action of $\C_+$ on $\cH'$ which is holomorphic 
on $\C_+^0$ and satisfies $(t,\alpha) \mapsto \pi^*(\exp_G(-tx))\alpha$ 
(cf.\ Example~\ref{ex:5.7}). 
As $\eta$ is $G$-equivariant and holomorphic, 
it is also equivariant with respect to the $\C_+$-actions 
on $M$ and the holomorphic action on $\bP(\cH')$ by 
$(z,[\alpha]) \mapsto [\alpha \circ e^{z \dd\pi(x)}]$. 
Therefore the $\C_+$ actions on $M$ and on $\cH'$ combine to a holomorphic 
action of $\C_+$ on $\eta^*\bL_{\cH'}$. Now we combine 
Theorem~\ref{thm:2.7}  with 
Remark~\ref{rem:5.8} to obtain (c). Here we only have to 
observe that the realization of $\cH$ by holomorphic sections 
of $\bL_{\cH'}$ leads on the dual bundle, whose complement 
of the zero-section can be identified with $\cH'\setminus \{0\}$,  
to the evaluation functional $K_{\alpha_v} = v \in \cH$ 
(cf.\ Example~\ref{ex:5.9}). 
\end{prf}

\begin{ex} \mlabel{ex:5.11}
An important special case where the 
requirements of Theorem~\ref{thm:5.7} are satisfied 
occurs if the $G$-action on the complex manifold 
$M$ extends to a holomorphic action of a complex 
Lie group $G_\C\supeq G$ with holomorphic exponential function. 

If $(\pi, \cH)$ is an irreducible representation of a compact 
Lie group $G$, then it is finite dimensional and in particular
bounded, so that $\pi$ extends to a holomorphic representation 
$\hat\pi \: G_\C \to \GL(\cH)$. Since the highest weight orbit 
$G[v_\lambda]$ is a compact complex manifold, it is also invariant 
under the $G_\C$-action on $\bP(\cH)$, and Theorem~\ref{thm:5.7}, 
applied to the maximal torus $H \subeq G$, 
implies that $I_\pi = \conv(\cO_\pi)$ for 
$\cO_\pi = \Phi_\pi(G[v_\lambda])$. 
\end{ex}

\begin{rem} (a) If $G$ is finite dimensional, then Theorem~\ref{thm:2.7} applies 
to all irreducible semibounded representations 
(\cite[Prop.~XII.3.6]{Ne00}), where $G/H$ is the 
highest weight orbit $G[v_\lambda] \subeq 
\bP(\cH_\lambda)$, and this eventually leads to 
$I_{\pi_\lambda} = \oline\conv(\cO_{-i\lambda})$. 

(b) We have already seen in \eqref{eq:extorbit} that, 
for finite dimensional groups, 
the highest weight orbit $G[v_\lambda] \subeq \bP(\cH_\lambda)$ 
can be characterized as the inverse image of the set 
$\Ext(I_{\pi_\lambda})$ of extreme points of the momentum set 
(Remark~\ref{rem:findim}). 
For the special case where $G$ is a compact group, one finds in 
\cite{Ha82} a different characterization of the highest weight 
orbit as the set of all those elements for which 
the subrepresentation of $\Sym^2(\cH_\lambda)$ generated by 
$v \otimes v$ is irreducible. It would be interesting to see if a 
similar result holds in other situations. 

If $G$ is finite dimensional, for any 
$[v] \in G[v_\lambda]$ the orbit of $[v \otimes v]$ is 
a complex manifold (actually a holomorphic image of 
$G[v_\lambda])$, and from that one can derive that the cyclic 
representation it generates is irreducible (\cite[Prop.~XV.2.7]{Ne00}). 

Another interesting characterization of the highest weight orbit 
of an irreducible representation of a compact Lie group 
is given in \cite{DF77}. Starting with an orthonormal basis 
$i F_1,\ldots, i F_n$ 
of the Lie algebra $\g$ with respect to an invariant scalar 
product, one defines the {\it invariant dispersion} of a state $[v]$ 
by 
$$ (\Delta F)^2 = \Big\la \sum_r (F_r - \la F_r \ra)^2 \Big\ra, 
\quad \mbox{ where } \quad 
\la A \ra = \frac{\la Av,v\ra}{\la v,v\ra} \quad \mbox{ for } \quad 
A = A^*. $$
If $C := \sum_r F_r^2$ is the corresponding Casimir operator, 
which acts on $\cH$ as a multiple $c \1$ 
of the identity, one easily finds that 
$$  (\Delta F)^2 = c - \|\Phi_\pi([v])\|^2, $$
which is minimal if $\|\Phi_\pi([v])\|$ is maximal. 
As $I_\pi = \conv(\cO_\pi)$ holds for a coadjoint orbit $\cO_\pi$, 
and the scalar product on 
$\g$, resp., $\g'$ is invariant, the orbit $\cO_\pi$ is contained 
in a sphere. Therefore the invariant dispersion $(\Delta F)^2$ 
is minimal in a state $[v]$ 
if and only if $\Phi_\pi([v]) \in \cO_\pi = \Ext(I_\pi)$. 
\end{rem}

\section{Invariant cones in Lie algebras} \mlabel{sec:6}

In this section we take a closer look at important examples 
of invariant convex cones in Lie algebras. 

\subsection{Invariant cones in unitary Lie algebras} 

\begin{exs} \mlabel{ex:liealgcon}
If $\cH$ is a complex Hilbert space, then the 
Lie algebra $\fu(\cH)$ of skew-hermitian bounded operators 
contains the open invariant cone 
$$ C_{\fu(\cH)} := \{ x \in \fu(\cH) \: i x \< 0 \} $$ 
(cf.\ Example~\ref{ex:4.9}). 
\glossary{name={$C_{\fu(\cH)}$},description={$=\{ x \in \fu(\cH) \: i x \< 0 \}$}}

More generally, for any unital $C^*$-algebra $\cA$, the 
Banach--Lie algebra $\fu(\cA) = \{ x \in \cA \: x^* = -x \}$ contains 
the open invariant cone 
$$ C_{\fu(\cA)} := \{ x \in \fu(\cA) \: i x \< 0 \}. $$
\end{exs} 

\subsection{Invariant cones in symplectic Lie algebras} 

\begin{defn} \mlabel{def:sympcon} 
If $V$ is a Banach space and 
$q \: V\to \R$ a continuous quadratic form, then we 
say that $q$ is {\it strongly positive definite}, written 
$q \> 0$, if $\sqrt{q}$ defines a Hilbert norm on $V$. 
In particular we are asking for $V$ to be complete with respect 
to this norm. 
\glossary{name={$q \> 0$},description={$q$ is positive definite, for quad. form}}

Now let $(V,\omega)$ be a strongly symplectic Banach space 
(cf.~Remark~\ref{rem:b.2}) 
and $\sp(V,\omega) \subeq \gl(V)$ be the corresponding symplectic 
Lie algebra. We associate to each $X \in \sp(V,\omega)$ the 
quadratic Hamiltonian $H_X(v) := \frac{1}{2}\omega(Xv,v)$ and obtain 
an open invariant cone by  \glossary{name={$H_X(v)$},description={$=\frac{1}{2}\omega(Xv,v)$, hamilt. funct., $X \in \sp(V,\omega)$}} 
$$ W_{\sp(V,\omega)} := \{ X \in \sp(V,\omega) \: H_X \> 0 \} $$ 
which \glossary{name={$W_{\sp(V,\omega)}$},description={$=\{ X \in \sp(V,\omega) \: H_X \> 0 \}$}}
is non-empty if and only if $V$ is topologically isomorphic 
to a real Hilbert space which carries a complex Hilbert space 
structure $\la \cdot,\cdot\ra$ with 
$$ \omega(v,w) =\Im \la v, w\ra\quad\mbox{ for } \quad v,w \in V $$
(cf.\ Proposition~\ref{prop:b.3}, \cite[Thm.~3.1.19]{AM78}).
\end{defn}

\begin{defn} \mlabel{def:hspcon} 
Actually $\sp(V,\omega)$ is a Lie subalgebra of the semidirect product 
$$ \hsp(V,\omega) := \heis(V,\omega) \rtimes \sp(V,\omega), $$ 
\glossary{name={$\hsp(V,\omega)$},description={$= \heis(V,\omega) \rtimes \sp(V,\omega)$, 
Jacobi--Lie algebra of $(V,\omega$)}} where $\heis(V,\omega) = \R \oplus_\omega V$ is the Heisenberg algebra 
associated to $(V,\omega)$, with the bracket 
$$ [(z,v), (z',v')] := (\omega(v,v'), 0). $$
Since every continuous linear functional on $V$ is of the form 
$i_x\omega$, the discussion in Remark~\ref{rem:b.2} 
implies that 
$\hsp(V,\omega)$ can be identified with the space of 
continuous polynomials of degree $\leq 2$ on $V$, endowed with the 
Poisson bracket
$$ \{f,g\} = \omega(X_g, X_f)\quad \mbox{ where }\quad 
i_{X_f} \omega = \dd f. $$
Let $\Heis(V,\omega) := \R \times V$ be the {\it Heisenberg group 
of $(V,\omega)$} with the multiplication 
$$ (t,v) (t',v') := (t + t' + \shalf\omega(v,v'), v + v'). $$
Then the Jacobi group 
\glossary{name={$\HSp(V,\omega)$},description={$=\Heis(V,\omega) \rtimes \Sp(V,\omega)$, Jacobi group}} 
$\HSp(V,\omega) := \Heis(V,\omega) \rtimes 
\Sp(V,\omega)$ acts by 
$ \sigma\big((c,w,g), v\big) := w + gv$ on $V$ and the 
corresponding derived action is given by 
$\dot\sigma(z,w,x)(v) = - w - Xv$, so that 
$$ (i_{\dot\sigma(z,w,x)}\omega)_v = -\omega(w + Xv, \cdot) 
= - i_{w + Xv}\omega. $$
Therefore this action is Hamiltonian with equivariant momentum map 
$$ \Phi \: \cH \to \hsp(V,\omega)', \quad 
\Phi(v)(c,w,A) := - c - \omega(w,v) - \shalf \omega(Av,v) $$
(cf.\ \cite[Prop.~A.IV.15]{Ne00}, where we use a different sign 
convention). 
\end{defn} 

\begin{lem} \mlabel{lem:6.4} The convex cone 
\begin{align*}
 W_{\hsp(V,\omega)}^+ 
&:= \{ (c,v,A) \in \hsp(V,\omega) \: (\forall v \in V) 
c + \omega(x,v) + H_A(v) > 0, H_A \> 0 \} 
\end{align*}
is open and invariant in the Banach--Lie algebra $\hsp(V,\omega)$. 
It is contained in the larger open invariant cone 
\begin{align*}
 W_{\hsp(V,\omega)} 
&:= \{ (c,v,A) \in \hsp(V,\omega) \: H_A \> 0 \}.
\end{align*}
\end{lem}
\glossary{name={$W_{\hsp(V,\omega)}$},description={$=\{ (c,v,A) \in \hsp(V,\omega) \: H_A \> 0\}$}}

\begin{prf} It is clear that $W :=  W_{\hsp(V,\omega)}^+$ is an invariant 
convex cone. It remains to show that it is open. 
If $f(v) = c + \omega(x,v) + \frac{1}{2}\omega(Av,v)$ is such that 
$H_A$ is strictly positive, then 
$\dd f(v) = i_x\omega + i_{Av}\omega,$
which vanishes if and only if $v = - A^{-1}x$. It follows in particular, 
that each such function has a unique minimal value which is given by  
$$ f(-A^{-1}x) = c - \omega(x,A^{-1}x) + \frac{1}{2} \omega(x,A^{-1}x) 
= c - \frac{1}{2} \omega(x,A^{-1}x). $$
Therefore the condition $f > 0$ is equivalent to 
$c > \frac{1}{2} \omega(x,A^{-1}x),$
showing that $W$ is indeed open in $\hsp(V,\omega)$. 
\end{prf}

From now on we assume that $\cH$ is a complex Hilbert space, 
$V = \cH_\R$ is the underlying real Banach space, and 
$\omega(v,w) := \Im \la v,w \ra$ is the corresponding symplectic form. 
Then $Iv = iv$ is a complex structure on $\cH_\R$ leaving $\omega$ 
invariant. It satisfies 
\begin{equation}
  \label{eq:iposcom}
\omega(Iv,v) = \Im \la iv,v\ra = \|v\|^2. 
\end{equation}
Formalizing this property leads to: 

\begin{defn}
We call a real linear complex structure $J \: \cH \to \cH$ 
{\it $\omega$-positive} if $\omega(Jv,w)$ is symmetric and 
positive definite and write $\cI_\omega$ for the set of $\omega$-positive 
complex structures on~$\cH$. 
\end{defn}

The following lemma is well known for the finite dimensional 
case, but it carries over to infinite dimensional Hilbert 
spaces. It implies in particular that 
$\cI_\omega = \cO_I$ is an adjoint orbit of the Lie algebra $\sp(\cH)$.

\begin{lem}\mlabel{lem:sp.1} The following assertions hold: 
  \begin{description}
  \item[\rm(i)] $\Sp(\cH) \cap \sp(\cH) 
= \{ g \in \GL(\cH_\R) \: I g^\top I = - g^{-1} = g \}$ 
is the set of complex structures $J$ on $\cH$ for which 
$\omega(Jv,w)$ is symmetric. 
  \item[\rm(ii)] $\cI_\omega = I e^{\fp}$ for 
$\fp := \{ x \in \sp(\cH) \: Ix = - xI\} 
= \{ x \in \sp(\cH) \: x^\top = x\}.$
  \item[\rm(iii)] The conjugation action of $\Sp(\cH)$ on 
$\cI_\omega$ leads to a diffeomorphism $\cI_\omega \cong \Sp(\cH)/\U(\cH)$. 
\item[\rm(iv)] 
If $A \in \sp(\cH)$ is such that $H_A \> 0$, then there exists a 
unique $\omega$-positive complex structure $J$ on $\cH$ commuting with $A$. 
\item[\rm(v)] $\cI_\omega^{\rm res} := \Ad(\Sp_{\rm res}(\cH))I 
=  \{ J \in \cI_\omega \: \|I - J \|_2 < \infty\} 
= \cI_\omega \cap \sp_{\rm res}(\cH)$. 
  \end{description}
\end{lem}

\begin{prf} (i) That $\omega(Jv,w)$ is symmetric is equivalent to 
$J \in \sp(\cH)$, and as $J^{-1} = -J$ characterizes complex structures, 
(i) follows. 

(ii) We have seen in \eqref{eq:iposcom} 
that $I$ is $\omega$-positive. 
Let $J$ be another $\omega$-positive complex structure. 
Writing $J = u e^x$ according to the polar 
decomposition of $\Sp(\cH)$ with 
$u \in \U(\cH)$ and $x\in  \fp$ (\cite[Thm.~I.6(iv)]{Ne02a}), 
we see that $J^2 = -\1$ is equivalent to $u^2 = -\1$ 
($u$ is a complex structure) and $ux=-xu$ ($x$ is antilinear 
with respect to $u$). If this is the case, then 
$$\omega(Jv,v) = \omega(u e^x v, v) 
= \omega(e^{-x/2}u e^{x/2}v,v)
= \omega(u e^{x/2}v,e^{x/2}v) $$
shows that $J$ is $\omega$-positive if and only if 
$u$ has this property. Since $u$ is complex linear, 
$\cH$ decomposes into $\pm i$-eigenspaces 
$\cH_\pm$ of~$u$. Then $\omega(uv,v)$ is positive definite on 
$\cH_+$ and negative definite on $\cH_-$, so that the 
$\omega$-positivity implies $u = I$. This proves (ii). 

(iii) For $g = u e^x \in \Sp(\cH)$, the relation 
$g^{-1}Ig = e^{-x} I e^x = I e^{2x}$
shows that $\Sp(\cH)$ acts transitively on 
$\cI_\omega$. As $\U(\cH)$ is the stabilizer of $I$, 
(iii) follows from the smoothness of the map 
$J = I e^x \mapsto x = \shalf\log(J^\top J)$. 

(iv) Let $(v,w)_A := \omega(Av,w)$ denote the real Hilbert 
structure on $\cH$ defined by $A$. Then 
$$ \omega(x,y) = \omega(A(A^{-1}x), y) = (A^{-1}x,y)_A $$ 
implies that $A^{-1}$ is skew-symmetric with respect to 
$(\cdot,\cdot)_A$, and the same holds for $A$ itself. 
Therefore $-A^2 \geq 0$ and 
$J :=  (-A^2)^{-1/2}A$ is a complex structure leaving 
$(\cdot,\cdot)_A$ invariant (cf.\ \cite[Thm.~3.1.19]{AM78}). 
Then 
$$ \omega(Jv,v) 
= \omega((-A^2)^{-1/2}Av,v)
= ((-A^2)^{-1/2}v,v)_A $$
is positive definite, so that $J$ is $\omega$-positive.  

If $J'$ is another $\omega$-positive complex structure commuting 
with $A$, then the construction of $J$ shows that it also commutes 
with $J$. Therefore $JJ'$ is an involution and since 
$J$ and $J'$ are $\omega$-positive, the $-1$-eigenspace of this 
involution is trivial, which leads to $J' = J$. 

(v) From the polar decomposition
$$\Sp_{\rm res}(\cH) = \U(\cH) \exp(\fp_2) 
\quad \mbox{ with } \quad 
\fp_2 := \fp \cap \sp_{\rm res}(\cH) 
= \{ x \in \fp \: \|x\|_2 < \infty\} $$
we derive that 
$\cI_\omega^{\rm res} = I e^{\fp_2}.$
For the entire function $F(z) := \frac{e^z - 1}{z}$, we then have 
$$ Ie^x - I = I F(x) x, $$
Since $F(x)$ is invertible for the real symmetric operator 
$x$ by the Spectral Mapping Theorem, it follows that 
$Ie^x - I$ is Hilbert--Schmidt if and only if $x$ 
has this property.

For $J = I e^x$ we also observe that 
$$[I,J] = -e^x - Ie^x I = e^{-x} - e^{x} 
= e^{-x}(\1 - e^{2x}) 
= - e^{-x} F(2x) 2 x, $$
so that $[I,J]$ is Hilbert--Schmidt if and only if $x$ is. 
This leads to $\cI_\omega^{\rm res} \break = \cI_\omega \cap \sp_{\rm res}(\cH)$. 
\end{prf}

The following theorem is a useful tool when dealing with invariant 
cones in symplectic Lie algebras. 

\begin{thm} \mlabel{thm:coneconj} For the canonical open invariant cones 
in $\sp(\cH)$, $\sp_{\rm res}(\cH)$ and $\hsp(\cH)$ we have the following 
conjugacy results: 
\begin{description}
\item[\rm(i)] $W_{\sp(\cH)} = \Ad(\Sp(\cH)) C_{\fu(\cH)}$ and 
$C_{\fu(\cH)} = W_{\sp(\cH)} \cap \fu(\cH)$. 
\item[\rm(ii)] $W_{\sp_{\rm res}(\cH)} :=  W_{\sp\cH)} \cap 
\sp_{\rm res}(\cH) = \Ad(\Sp_{\rm res}(\cH)) C_{\fu(\cH)}$. 
\glossary{name={$W_{\sp_{\rm res}(\cH)}$},description={$=  W_{\sp\cH)} \cap \sp_{\rm res}(\cH)$}} 
\item[\rm(iii)] $W_{\hsp(\cH)} = \Ad(\HSp(\cH))(\R \times \{0\} \times 
C_{\fu(\cH)})$ 
for the corresponding Lie group \linebreak 
$\HSp(\cH) = \Heis(\cH) \rtimes \Sp(\cH)$. 
\item[\rm(iv)] $W_{\hsp_{\rm res}(\cH)} := 
W_{\hsp(\cH)} \cap \hsp_{\rm res}(\cH) 
= \Ad(\HSp_{\rm res}(\cH))(\R \times \{0\} \times C_{\fu(\cH)})$ 
for $\HSp_{\rm res}(\cH) := \Heis(\cH) \rtimes \Sp_{\rm res}(\cH)$. 
\end{description}
\end{thm}

\begin{prf} (i) If $A \in \sp(\cH)$ is complex linear, then 
$A^* =- A$, so that $\omega(Av,v) = \Im \la Av,v\ra =  \la -iAv,v\ra$. 
Therefore $A \in W_{\sp(\cH)}$ is equivalent to $iA \< 0$, i.e., 
$W_{\sp(\cH)} \cap \fu(\cH) = C_{\fu(\cH)}$.  

For $A \in W_{\sp(\cH)}$ we find with 
Lemma~\ref{lem:sp.1}(iv) a $J \in \cI_\omega$ commuting with $A$. 
For any $g \in \Sp(\cH)$ with $J = \Ad(g)I$, 
whose existence follows from Lemma~\ref{lem:sp.1}(iii), 
we conclude that 
$\Ad(g)^{-1}A$ commutes with $I$, hence is contained in
$W_{\sp(\cH)} \cap \fu(\cH) = C_{\fu(\cH)}$. This proves (i). 

(ii) In view of Lemma~\ref{lem:sp.1}(v), it suffices to show that 
for any $A \in W_{\sp_{\rm res}(\cH)}$ the corresponding $\omega$-positive 
complex structure $J =  (-A^2)^{-1/2}A$ from Lemma~\ref{lem:sp.1}(iv) 
is contained in 
$\cI_\omega^{\rm res}$, i.e., $[I,J]$ is Hilbert--Schmidt. 
Since $A^2$ commutes with $I$, we have 
$[I,J] = (-A^2)^{-1/2}[I,A]$, so that the invertibility of 
$(-A^2)^{-1/2}$ implies that $A \in \sp_{\rm res}(\cH)$ is equivalent 
to $J \in \sp_{\rm res}(\cH)$. 

(iii) As we have seen in the proof of Lemma~\ref{lem:6.4} above, 
for each element $(c,x,A) \in W_{\hsp(\cH)}$, the Hamiltonian function
$f(v) = c + \omega(x,v) + H_A(v)$
has a unique minimum in $- A^{-1}x$. Since the adjoint action of the 
Heisenberg group $\Heis(\cH) \subeq \HSp(\cH)$ corresponds to 
a translation action on $\cH$, each adjoint orbit $\cO_X$ in 
$W_{\hsp(V,\omega)}$ contains an element $Y$ whose corresponding 
Hamiltonian function is minimal in $0$, so that
$Y \in \R \times \{0\} \times \sp(\cH)$. 
In view of (i), each orbit in $W_{\hsp(V,\omega)}$ meets the subalgebra 
$\R \times \{0\} \times \fu(\cH)$, and this proves (iii).

(iv) follows by combining the argument under (iii) with the 
proof of (ii). 
\end{prf}

\subsection{Invariant Lorentzian cones} \mlabel{subsec:6.3}

If $(\g,\beta)$ is a Lorentzian Lie algebra, 
i.e., $\beta$ is an invariant Lorentzian form (which is negative definite 
on a closed hyperplane), then 
each half of the open double cone $\{ x \in \g \: \beta(x,x) > 0 \}$ 
is an open invariant cone in $\g$. 

\begin{ex} \mlabel{ex:6.5} 
A particularly important example is $\g = \fsl_2(\R)$ with 
$\beta(x,y) = -\tr(xy)$. Indeed, the basis 
$$
h =
\begin{pmatrix}
1 &  \phantom{-}0 \\
0 & -1
\end{pmatrix},
\quad
u =
\begin{pmatrix}
\phantom{-}0 & 1 \\
-1 & 0
\end{pmatrix},
\quad
t =
\begin{pmatrix}
0 & 1 \\
1 & 0
\end{pmatrix}
$$
is orthogonal with 
$$ \beta(x h + y u + z t, x h + y u + z t) = -2 x^2 + 2 y^2 - 2 z^2. $$
For $\g = \fsl_2(\R)$, the adjoint group 
coincides with the identity component \break $\SO(\g,\beta)_0 
\cong \SO_{1,2}(\R)_0$, which implies that the adjoint and 
coadjoint orbits are the connected components of the level surfaces of the 
associated quadratic form, and the $0$-level surface of isotropic 
vectors decomposes into the $0$-orbit and two isotropic orbits 
lying in the boundary of the double cone. 
This description of the adjoint orbits implies in 
particular that there are precisely two non-trivial open invariant 
cones, namely the connected components of the set 
$\{ x \: \beta(x,x) > 0\}$. 
\end{ex}

\begin{rem}
Other examples of Lorentzian Lie algebras 
arise as double extensions 
from Lie algebras $\g^0$, endowed with an invariant scalar product 
$\kappa^0$: 
If $D \in \der(\g^0,\kappa^0)$ is a continuous skew-symmetric derivation, 
then 
$\g := \R \times \g^0 \times \R$ is a Lie algebra with respect to the 
bracket 
$$[(z,x,t), (z',x',t')] = (\kappa^0(Dx,x'), tDx' - t'Dx + [x,x'],0)$$
and the continuous symmetric bilinear form 
$$ \kappa((z,x,t), (z',x',t')) := zt' + z't + \kappa^0(x,x'), $$ 
is easily seen to be invariant. 
The pair $(\g,\kappa)$ is called a {\it double 
extension} of $(\g^0,\kappa^0)$ (cf.\ \cite{MR85}).  
\end{rem}

\begin{ex} \mlabel{ex:lorentz} (a) If $\fk$ is a compact Lie algebra and 
$\g^0 := C^\infty(\bS^1,\fk)$ is the corresponding loop algebra, 
then we identify its elements with $2\pi$-periodic functions on $\R$. 
With an invariant scalar product $\kappa_\fk$ on $\fk$, we obtain 
the invariant scalar product 
$$ \kappa^0(\xi,\eta) := \int_0^{2\pi} \kappa_\fk(\xi(t), \eta(t))\, dt $$
on $\g^0$, and the derivation $D\xi := \xi'$ is skew-symmetric. 
The corresponding double extension produces the (unitary forms) of the 
untwisted affine Kac--Moody Lie algebras. 

(b) If $\g^0 := \fu_2(\cH)$ is the Lie algebra of skew-hermitian 
Hilbert--Schmidt operators on the Hilbert space $\cH$ and 
$A = - A^* \in \fu(\cH)$, then 
$\kappa^0(x,y) := \tr(xy^*) = -\tr(xy)$ is an invariant scalar product on 
$\g^0$ and one obtains a double extension for the derivation 
$D(x) := [A,x]$. It is non-trivial if and only if 
$A \not\in \R i \1 + \fu_2(\cH)$ (cf.\ \cite{Ne02a}). 
\end{ex}

\begin{prob} Classify infinite dimensional Lorentzian Lie algebras 
$\g$ which are complete in the sense that for $x \in \g$ 
with $\beta(x,x) > 0$ the orthogonal space is a Hilbert space with 
respect to $-\beta$. The construction in Example~\ref{ex:lorentz}(b) 
produces interesting examples. 

Finite dimensional indecomposable Lorentzian Lie algebras 
have been classified by Hilgert and Hofmann in \cite{HiHo85}. 
The simple result is that an indecomposable finite dimensional 
Lorentzian Lie algebra is either $\fsl_2(\R)$, endowed with the 
negative of its Cartan--Killing form, or a double extension 
of an abelian Lie algebra, defined by an invertible skew-symmetric 
derivation $D$ (\cite[Thm.~II.6.14]{HHL89}). 
\end{prob}

\subsection{Invariant cones of vector fields} 

\begin{ex}
  \mlabel{ex:vectcirc}
Let $\g = \cV(\bS^1) = C^\infty(\bS^1) \partial_\theta$ 
be the Lie algebra of smooth vector fields 
on the circle $\bS^1 \cong \R/\Z$, where $\partial_\theta := 
\frac{\partial}
{\partial \theta}$ denotes the generator of the right rotations. Then 
$$ W_{\cV(\bS^1)} := \Big\{ f \partial_\theta \: f > 0 \Big\} $$
is an open invariant cone (cf.\ Section~\ref{sec:7}). 
\end{ex}
\glossary{name={$W_{\cV(\bS^1)}$},description={$=\Big\{ f \partial_\theta \: f > 0 \Big\}$, 
invariant cone in $\cV(\bS^1)$}}

\begin{ex}
  \mlabel{ex:conform} 
The preceding example has a natural higher dimensional generalization. 
Let $(M,g)$ be a compact Lorentzian manifold possessing a timelike 
vector field $T$, i.e., $g_m(T(m), T(m)) > 0$ for every $m \in M$. 
Then 
$$ W := \{ X \in \cV(M) \: (\forall m \in M)\, 
g(X,X) > 0, g(X,T) > 0 \} $$
is an open convex cone in the Fr\'echet space $\cV(M)$ and its
intersection with the subalgebra $\conf(M,g)$ of conformal vector fields 
is an open invariant cone. 
For $M = \bS^1$ all vector fields are conformal and we thus obtain 
Example~\ref{ex:vectcirc}. 
\end{ex}

\subsection{Invariant cones and symmetric Hilbert domains} 

\begin{ex}
  \mlabel{ex:boundom} (a) 
We recall from Example~\ref{ex:liegrp}(g) the concept 
of a symmetric Hilbert domain and write 
$G := \Aut(\cD)_0$ for the identity component of its automorphism 
group. We assume w.l.o.g.\ that $\cD$ is the open unit ball of a 
complex Banach space $V$ (cf.\ \cite{Ka83}, \cite{Ka97}). 
Let $K = G \cap \GL(V)$ be the subgroup of linear 
automorphisms of $\cD$, i.e., the group of complex linear isometries of 
$V$ and let $\fk = \L(K)$ be its Lie algebra. Then 
$$ W_\fk := \{ x \in \fk \: \|e^{i x}\| < 1 \} $$ 
is a pointed open convex cone in the subalgebra $\fk$, and 
$$ W := \Ad(G)W_\fk $$ 
is a pointed open invariant convex cone in $\g$ 
(cf.\ \cite[Thm.~V.9]{Ne01a}, and \cite[Thm.~5]{Vin80} for the finite 
dimensional case). 
By definition, this open cone has the interesting 
property that every orbit in $W$ meets~$\fk$. 

(b) Because they will also show up in the following, we take a closer 
look at some characteristic examples. In Examples~\ref{ex:liegrp}(g), 
we have seen that the group 
$G = \U_{\rm res}(\cH_+, \cH_-)$ acts naturally on the 
circular symmetric Hilbert 
domain $\cD := \{ z \in B_2(\cH_+, \cH_-) \: \|z\| < 1\}$. 
Here the stabilizer of $0$ is 
$$ K := \U(\cH_-) \times \U(\cH_+), $$
which acts by $(a,d)z = azd^{-1}$. From the Hilbert space 
isomorphism \break $B_2(\cH_+, \cH_-) 
\cong \cH_- \hat\otimes \cH_+^*$, one now derives with 
\begin{align*}
\|e^{ix_-} \otimes e^{-ix_+}\| 
&= \|e^{ix_-}\|\|e^{-ix_+}\| 
= \sup\Spec(i x_-) + \sup\Spec(-i x_+)\\
&= \sup\Spec(i x_-) - \inf\Spec(i x_+) 
\end{align*}
that 
$$ W_\fk = \{ (x_-, x_+) \in \fk \cong 
\fu(\cH_-) \times \fu(\cH_+) \: 
\sup\Spec(i x_-) < \inf\Spec(i x_+) \}, $$
and $W = \Ad(G)W_\fk$ is the corresponding open invariant cone in 
$\g = \fu(\cH_+, \cH_-)$. Since the action of $G$ on $\cD$ is not 
faithful, this cone has a non-trivial edge $H(W) = \R i\1 
= \L(\T \1)$. 

(c) For $\cH_+= \cH_- = \cH$ and the subgroup 
$G = \Sp_{\rm res}(\cH) \subeq \U_{\rm res}(\cH_+, \cH_-)$,  
we have $K \cong \U(\cH)$, corresponding to the pairs of the  form 
$(a,a^{-\top}) \in \U(\cH) \times \U(\cH)$, and, accordingly, 
$$ W_\fk = \{ x \in \fu(\cH)\ i x \< 0 \} = C_{\fu(\cH)}$$
(Examples~\ref{ex:liegrp}(h) and \ref{ex:liealgcon}). 
From Theorem~\ref{thm:coneconj}(ii) it now follows that 
\begin{equation}
  \label{eq:sprescone}
W = \Ad(G)W_\fk =   W_{\sp_{\rm res}(\cH)}.
\end{equation}
\end{ex}

\subsection{A general lemma} 

The following lemma captures the spirit of some arguments 
in the previous constructions of invariant convex cones 
in Lie algebras in a quite natural way. 
F.i., it applies to finite dimensional simple 
algebras as well as $\cV(\bS^1)$. 

\begin{lem}\mlabel{lem:doubcone} 
Suppose that the element $d \in \g$ has the following 
properties: 
\begin{description}
\item[\rm(a)] The interior $W_{\rm min}$ of the invariant convex cone 
generated by $\cO_d = \Ad(G)d$ is non-empty and different from $\g$. 
\item[\rm(b)] There exists a continuous linear  projection 
$p \: \g \to \R d$ which preserves every open and closed convex subset. 
\item[\rm(c)] There exists an element $x \in \g$ for which 
$p(\cO_x)$ is unbounded. 
\end{description}
Then each non-empty open invariant cone $W$ contains $W_{\rm min}$ or 
$-W_{\rm min}$, and for $\lambda \in \g'$ the following are equivalent
\begin{description}
\item[\rm(i)] $\lambda \in \g'_{\rm seq}$, i.e., 
$\cO_\lambda$ is semi-equicontinuous. 
\item[\rm(ii)] $\cO_\lambda(d)$ is bounded from below or above. 
\item[\rm(iii)] $\lambda \in W_{\rm min}^\star \cup - W_{\rm min}^\star$. 
\end{description}
\end{lem}

\begin{prf} If $W \subeq \g$ is an open invariant cone, 
then $p(W) \subeq \R d$ contains an open half line, so that (b) 
implies that $W$ contains either $d$ or $-d$. Accordingly, we 
then have $W_{\rm min} \subeq W$ or $-W_{\rm min} \subeq W$. 

(iii) $\Rarrow$ (i): As $W_{\rm min}$ is open, 
$W_{\rm min}^\star \cup - W_{\rm min}^\star$ consists 
of semi-equicontinuous coadjoint orbits (Example~\ref{ex:2.3}(b)). 

(i) $\Rarrow$ (ii): Let $\lambda \in \g'_{\rm seq}$. 
Since $B(\cO_\lambda)^0$ is an open invariant cone, (a) implies that 
it either contains $d$ or $-d$, which is (ii). 

(ii) $\Rarrow$ (iii): Assume w.l.o.g.\ that 
$\inf \cO_\lambda(d) > - \infty$. We claim that 
$\cO_\lambda \subeq W_{\rm min}^\star$, which is equivalent to 
$\cO_\lambda(d) \geq 0$. 

To this end, we consider $\mu \in \g'$ defined by $p(x) = \mu(x)d$. 
Then $-d \not\in W_{\rm min}$ and (b) imply that 
$\mu \in W_{\rm min}^\star$, so that 
$\mu(\cO_d) \geq 0$. If $\mu(\cO_d)$ is bounded from above, 
then $\pm d \in B(\cO_\mu)$, so that the invariance of the convex 
cone $B(\cO_\mu)$ leads to $\pm W_{\rm min}\subeq B(\cO_\mu)$, 
and hence to $\g = B(\cO_\mu)$, contradicting (c). We conclude that 
$\mu(\cO_d)$ is unbounded. 
Applying (b) to the closed convex subset $C := \oline{\conv}(\cO_d)$, 
it follows that $[1,\infty[ \cdot d \subeq C$, and hence that 
$d \in \lim(C)$ (Lemma~\ref{lem:limcone}(iii)). 
As $\lim(C)$ is an invariant cone, (a) entails 
$W_{\rm min} \subeq \lim(C)$. We finally conclude that 
$\lambda \in B(\cO_d) = B(C) \subeq \lim(C)^\star 
\subeq W_{\rm min}^\star$ (Lemma~\ref{lem:limcone})(vi)).
\end{prf}

The preceding lemma applies in particular to finite dimensional 
simple non-compact Lie algebras (\cite{Vin80}): 
\begin{prop} \mlabel{prop:simpfin} 
For a finite dimensional simple non-compact  
Lie algebra $\g$, the following assertions hold: 
\begin{description}
\item[\rm(i)] Every non-empty invariant convex subset $C\not=\{0\}$ 
has interior points. 
\item[\rm(ii)] If $\g$ contains a proper open invariant convex cone, 
then there exist minimal and maximal open invariant cones
$W_{\rm min} \subeq W_{\rm max}$ such that for any other open invariant 
cone $W$ we either have 
$$W_{\rm min} \subeq W \subeq W_{\rm max} \quad \mbox{ or } \quad 
W_{\rm min} \subeq -W \subeq W_{\rm max}. $$
In this case the set of semi-equicontinuous coadjoint orbit in 
$\g'$ coincides with $W_{\rm min}^\star \cup - W_{\rm min}^\star$.
\end{description}
\end{prop}

\begin{prf} (i) If $\eset \not= C \not= \{0\}$ is invariant and convex, 
then $\Spann(C)$ is a non-zero ideal of $\g$, hence equal to $\g$. 
If the affine subspace generated by $C$ is proper, then its translation 
space is a hyperplane ideal of $\g$, contradicting the simplicity of 
$\g$. Therefore $C$ generates $\g$ as an affine space, hence 
has interior points because $\dim \g < \infty$ and 
$C$ contains a simplex of maximal dimension. 

(ii) Let $\g = \fk \oplus \fp$ be a Cartan decomposition of 
$\g$ and $\fz := \fz(\fk)$ be the center of $\fk$. 
Then the existence of invariant cones implies that 
$\fz = \fz_\g(\fk) = \R d$ is one-dimensional (\cite{Vin80}). 
Therefore we have a fixed point projection 
$p_\fz \: \g \to \fz$ with respect to the action of the 
compact subgroup $e^{\ad \fk}$, and this projection 
preserves all open and closed invariant convex subsets 
(Proposition~\ref{prop:project}(a)). 
If $p_\fz(\cO_x)$ is bounded for every $x \in \g$, then 
$p_\fz$,  considered as a linear functional on $\g$, 
has a bounded orbit. As $\g$ is simple, this implies that all 
coadjoint orbits are bounded, and this contradicts the non-compactness 
of $\g$. Therefore (a)-(c) in Lemma~\ref{lem:doubcone} are satisfied, 
which implies that the set of semi-equicontinuous coadjoint 
orbits is $W_{\rm min}^\star \cup - W_{\rm min}^\star$. 
\end{prf}

\section{Connections to $C^*$-algebras} \mlabel{sec:7}

In this section we discuss two aspects of semibounded representation 
theory in the context of $C^*$-algebras. The first one concerns 
the momentum sets of restrictions of 
representations of a unital $C^*$-algebra $\cA$ to its unitary group 
$\U(\cA)$, and the second one concerns covariant representations 
for $C^*$-dynamical systems defined by a Banach--Lie group $H$ acting 
on a $C^*$-algebra $\cA$. Using the results from Appendix~\ref{app:a}, 
it follows that covariant representations lead to 
smooth unitary representations of the Lie group 
$\U(\cA^\infty) \rtimes H$, so that spectral conditions for covariant 
representations can 
be interpreted in terms of semibounded representations and the 
theory of invariant cones in Lie algebras becomes available. 

In this subsection $\cA$ denotes a unital $C^*$-algebra and 
$\U(\cA)$ its unitary group, considered as a Banach--Lie group 
(Example~\ref{ex:liegrp}(a)). We identify the state space 
$$S(\cA) := \{ \phi \in \cA' \: \phi(\1) = \|\phi\| = 1\} 
\subeq \{ \phi \in \cA' \: \phi(\fu(\cA)) \subeq i \R\} $$ 
\glossary{name={$S(\cA)$},description={states of the $C^*$-algebra $\cA$}}

of $\cA$ with a subset of $\fu(\cA)'$ by mapping 
$\phi \in S(\cA)$ to the real-valued functional $-i \phi\res_{\fu(\cA)} 
\in \fu(\cA)'$ (cf.\ \cite[Sect.~X.5]{Ne00} for more details). 

\subsection{Representations of the unitary group} 

If $(\pi, \cH)$ is a unitary representation of $\U(\cA)$ 
obtained by restricting an algebra representation, then 
its momentum set is simply given by 
$$ I_\pi = \{ \phi \in S(\cA) \: \phi(\ker \pi) =\{0\}\} 
= S(\cA) \cap (\ker \pi)^\bot \cong S(\pi(\cA)) $$
(cf.\ \cite[Thm.~X.5.13]{Ne00}). In particular, the momentum set 
is completely determined by the kernel of the representation $\pi$, resp., 
the $C^*$-algebra $\pi(\cA)$. This is why representations with the 
same kernel are called {\it physically equivalent} in 
algebraic Quantum Field Theory (\cite{HK64}). 

In \cite[Thm.~9.1]{Dix64} one finds a characterization 
of the separable $C^*$-algebras $\cA$ of type $I$ as those for 
which irreducible representations are determined by their 
kernels, hence by their momentum sets. 
More generally, postliminal $C^*$-algebras 
have this property (\cite[Thm.~4.3.7]{Dix64}). 
We conclude in particular that the existence 
of separable $C^*$-algebras $\cA$ which are not of type $I$ implies the 
existence of non-equivalent irreducible representations 
$(\pi_1, \cH_1)$ and $(\pi_2, \cH_2)$ with 
$\ker \pi_1 = \ker \pi_2$, and hence with 
$I_{\pi_1} = I_{\pi_2}$. We thus observe: 

\begin{thm} Bounded irreducible unitary representations 
of Banach--Lie groups 
are in general not determined up to equivalence by their 
momentum sets.   
\end{thm}

\begin{prop} {\rm(Momentum sets of irreducible representations)} \mlabel{prop:7.2}
For irreducible representations $(\pi, \cH)$ of a $C^*$-algebra 
$\cA$, the following assertions hold: 
\begin{description}
\item[\rm(i)] $\U(\cA)$ acts transitively on $\bP(\cH)$, so that 
$\cO_\pi := \im(\Phi_\pi)$ is a single coadjoint orbit with 
$I_\pi = \oline{\conv}(\cO_\pi)$.  
\item[\rm(ii)] $\cO_\pi$ consists  of pure states, i.e., 
$\cO_\pi \subeq \Ext(S(\cA))$. 
\item[\rm(iii)] Two irreducible representations 
$\pi_1$ and $\pi_2$ are equivalent if and only if 
$\cO_{\pi_1} = \cO_{\pi_2}$. 
\end{description}
\end{prop}

\begin{prf} (i) \cite[Thm.~2.8.3]{Dix64} 

(ii) \cite[Prop.~2.5.4]{Dix64} 

(iii) (cf.~\cite[Cor.~2.8.6]{Dix64})  
From the naturality of the momentum map it follows that 
equivalent representations have the same 
orbits. The converse follows from the fact that an irreducible 
representation can be recovered from any of its pure states $\phi$ by 
the GNS construction, and states in the same $\U(\cA)$-orbit 
lead to equivalent representations. 
\end{prf}

In \cite{BN10} these results are generalized to irreducible 
representations of $\U(\cA)$ occurring in tensor products 
of algebra representations and their duals. 

\begin{exs} \mlabel{ex:7.3} 
(a) If $\pi$ is the identical representation 
of the $C^*$-algebra $\cA = B(\cH)$ on $\cH$, then the corresponding 
momentum map is given by 
$$ \Phi_\pi([v])(x) = \frac{1}{i} \frac{\la xv, v\ra}{\la v,v\ra} 
= - i \tr(xP_v), $$
where $P_v(w) = \frac{\la w,v\ra}{\la v,v\ra} v$ is the orthogonal 
projection onto $[v]$. Clearly, this representation is faithful 
and irreducible, so that $I_\pi = S(\cA)$. On 
the other hand $\cO_\pi = \im(\Phi_\pi)$ can be identified 
with the set of rank-one projections, as elements of the 
dual space $\cA'$. With the trace pairing, we can embed 
the subspace $\Herm_1(\cH)$ of hermitian trace class operators 
in to the dual $\cA'$. Then 
$$ \Herm_1(\cH) \cap S(\cA) = \{ S \in \Herm_1(\cH) \: 
S \geq 0, \tr S = 1\}, $$
and the Spectral Theorem for compact hermitian operators 
implies that 
$$ \Ext(S(\cA)) \cap \Herm_1(\cH) = \Phi_\pi(\bP(\cH)) $$
is a single coadjoint orbit. However, since the 
{\it Calkin algebra} $\cA/K(\cH)$ is non-trivial if $\dim \cH = \infty$, 
$\cA$ also has pure states 
vanishing on the ideal $K(\cH)$ of compact operators and 
$\U(\cH)$ does not act transitively on $\Ext(S(\cA))$. 

(b) For the natural representation of $G = \U(\cH)$ on 
the symmetric powers $S^n(\cH)$, the elements 
$[v^n]$, $0\not=v \in \cH$, form a single $G$-orbit on 
which the complex group $\GL(\cH) = \cA^\times$ acts holomorphically  
(Proposition~\ref{prop:7.2}(i)). Since 
$$\Phi_{S^n(\pi)}([v^n]) = n \Phi_\pi([v]) \quad \mbox{ for } \quad 
0 \not= v\in \cH, $$
it thus follows from Example~\ref{ex:5.11} that 
\begin{equation}
  \label{eq:sym-mom}
I_{S^n(\pi)} = n I_\pi. 
\end{equation}
For a generalization to more general subrepresentations of 
$\cH^{\otimes n}$ we refer to \cite{BN10}. 
\end{exs}

\subsection{$C^*$-dynamical systems} 

\begin{defn} \mlabel{def:7.3} 
Let $G$ be a topological group and $\cA$ a $C^*$-algebra. 
A {\it $C^*$-dynamical system} is a triple 
$(\cA, G, \alpha)$, where $\alpha \: G \to \Aut(\cA), g \mapsto \alpha_g$, 
is a homomorphism defining a continuous action of $G$ on $\cA$. 
\end{defn}

\begin{thm}\mlabel{thm:a.8}
  \mlabel{cor:cstarsmooth} {\rm(\cite{Ne10})} 
If $G$ is a Banach--Lie group and 
$(\cA, G, \alpha)$ a $C^*$-dynamical system, then 
the space $\cA^\infty$ of smooth vectors is a 
Fr\'echet algebra with respect to the locally convex 
topology defined by the  seminorms 
$$ p_n(a) := \sup \{ \|\dd\alpha(x_1) \cdots \dd\alpha(x_n)a\| \: 
x_i \in \g, \|x_i\|\leq  1\}, \quad n \in \N_0, $$
and the action of $G$ on $\cA^\infty$ is smooth. 
If, in addition, $\cA$ is unital, then $\cA^\infty$ is a continuous 
inverse algebra. 
\end{thm}

Now let $(\cA,H,\alpha)$ be a $C^*$-dynamical system, 
where $H$ is a Banach-Lie group and $\cA$ is a unital $C^*$-algebra. 
Then $\U(\cA^\infty)$ 
carries a natural Fr\'echet--Lie group structure 
(cf.\ Example~\ref{ex:liegrp}(a)), and we can form the 
semidirect product Lie group $G := \U(\cA^\infty) \rtimes H$. 

For the proof of the following theorem, we record 
a general observation on invariant cones. 

\begin{lem} \mlabel{lem:2.1} 
If $W \subeq \g$ is an  open invariant 
convex cone, then we have for each element $x \in \g$ satisfying 
$(\ad x)^2 =0$ the relation $[x,\g] \subeq H(W)$. 
In particular, $x$ is central if $W$ is pointed. 
\end{lem}

\begin{prf} Let $w \in W$. 
Then, for each $t \in \R$, we have 
$w + \R [x,w] = e^{\R \ad x}w \subeq W,$
so that $[x,w] \in H(W)$ (Lemma~\ref{lem:limcone}). 
Now $\g = W - W$ leads to $[x,\g] \subeq H(W)$. 
If, in addition, $W$ is pointed, this implies that 
$[x,\g] = \{0\}$,  i.e., $x \in \z(\g)$. 
\end{prf}

\begin{thm} \mlabel{thm:7.4} 
Let $(\pi,\rho, \cH)$ be a covariant representation of 
$(\cA,H,\alpha)$, i.e., $\pi$ is a non-degenerate representation of 
$\cA$ on $\cH$ and $\rho$ is a unitary representation of $H$ satisfying 
$$ \pi(\alpha_g A) = \rho(g)\pi(A) \rho(g)^{-1} \quad \mbox{ for } 
\quad A \in \cA, g \in H. $$
Then the following assertions hold: 
\begin{description}
\item[\rm(i)] The corresponding representation 
$\hat\pi(a,h) := \pi(a) \rho(h)$ of $G = \U(\cA^\infty) \rtimes H$ 
is smooth if and only if 
$\rho$ has this property. 
\item[\rm(ii)] $W_{\hat\pi} = \fu(\cA^\infty) \times W_\rho$ and 
$\hat\pi$ is semibounded if and only if $\rho$ has this property 
if and only if $C_{\hat\pi} \not=\eset$. 
\item[\rm(iii)] If $\cA$ is commutative and 
$\hat\pi$ is semibounded, then the identity component $H_0$ of $H$ 
acts trivially on $\pi(\cA)$. 
\end{description}
\end{thm}

\begin{prf} (i) For every $H$-smooth vector $v \in \cH$, the smoothness of the 
map $\U(\cA) \times H \to \cH, (a,g) \mapsto \pi(a)\rho(g)v$ 
follows from the smoothness of the action of the Banach--Lie group 
$\U(\cA)$ on $\cH$. Since the inclusion $\U(\cA^\infty) \to \U(\cA)$ 
is smooth, it follows that every $H$-smooth vector is smooth for 
$G$, so that the corresponding unitary representation 
$\hat\pi \: G \to \U(\cH), (a,g) \mapsto \pi(a)\rho(g)$ 
of $G$ is smooth whenever $\rho$ has this property. 

(ii) Clearly $i\1 = \dd\pi(i\1) \in \dd\hat\pi(\g)$, so that 
$\hat\pi$ is semibounded if and only if 
$C_{\hat\pi} \not=\eset$ (Proposition~\ref{prop:repcone}). 
Since $\cA$ acts by bounded operators, we have 
$\fu(\cA^\infty) \subeq H(W_{\hat\pi})$, and thus 
$$ W_{\hat\pi} 
= \fu(\cA^\infty) \times (W_{\hat\pi} \cap \fh) 
= \fu(\cA^\infty) \times W_\rho. $$
Therefore $\hat\pi$ is semibounded if and only if $\rho$ is semibounded. 

(iii) As $\fu(\cA^\infty)$ is an abelian 
ideal of $\g$, Lemma~\ref{lem:2.1} implies that 
$[\fu(\cA^\infty),\fh] \subeq H(C_\pi) = \ker \dd\hat\pi$, i.e., 
that $H_0$ acts trivially on the $C^*$-algebra $\pi(\cA)$. 
\end{prf}

\begin{rem}
A closely related fact is well-known in the 
context of $C^*$-dynamical systems with the group 
$H = \R^d$. 
To explain the connection, let $C \subeq \fh'$ be a closed convex cone. 
Then we say that $(\pi, \rho, \cH)$ satisfies the 
{\it $C$-spectrum condition} if $-i\dd\rho(x) \geq 0$ for 
$x \in C^\star$, i.e., $I_\rho \subeq C$. 
If $C^\star$ has interior points, we have just seen that 
the $C$-spectrum condition implies that $\hat\pi$ is semibounded, 
and if $\pi$ is faithful and $\cA$ is commutative, 
this can only happen if $H$ acts trivially on $\cA$ 
(cf.\ \cite[Thm.~IV.6.2]{Bo96}).
To obtain non-trivial situations, one has to consider non-commutative 
algebras. Typical examples arise for any 
semibounded representation $(\rho, \cH)$ of $H$ 
for $\cA := K(\cH)$ (compact operators on $\cH$) and 
$\alpha_g(A) := \rho(g)A\rho(g)^{-1}$. 
\end{rem}

\begin{ex} \mlabel{ex:7.5} 
(a) If $(\rho, \cH)$ is a smooth representation of $G$, 
then the corresponding action of $G$ on $K(\cH)$ defined by 
$\alpha_g(A) := \rho(g)A\rho(g)^{-1}$ also has a dense space $K(\cH)^\infty$ 
of smooth vectors because for every pair $(v,w)$ of smooth vectors 
the corresponding rank-one operator $P_{v,w}$, defined by 
$P_{v,w}(x) := \la x,w \ra v$ satisfies 
$$ \alpha_g P_{v,w} = P_{\rho(g)v, \rho(g)w}, $$
which easily implies that its orbit map is smooth. 

(b) If $\cH$ is a complex Hilbert space and 
$\cA := \Car(\cH)$, then the canonical action of the orthogonal 
group $\rO(\cH)$ (of the underlying real Hilbert space) 
on $\cA$ is continuous and the subalgebra 
$\cA^\infty$ of smooth vectors is dense because it contains 
$a(\cH)$ and hence the $*$-subalgebra generated by this subset 
(cf.\ Section~\ref{sec:10}). 

(c) To find a similar situation for the CCR is not so obvious because 
the Weyl algebra $\CCR(\cH)$, i.e., the $C^*$-algebra 
defined by the generators 
$W(f)$, $f \in \cH$, and the {\it Weyl relations} 
$$ W(f)^* = W(-f), \quad 
W(f) W(h) = e^{\frac{i}{2}\Im \la f, h \ra} W(f+h) $$
is a very singular object. The map $W \: \cH \to \CCR(\cH)$ 
is discontinuous, even on every ray in $\cH$ 
(\cite[Thm.~5.2.8]{BR97}). 
Since the action of the symplectic 
group $\Sp(\cH)$ preserves these relations, it acts by 
$\alpha_g(W(f)) = W(gf)$ on $\CCR(\cH)$, but this action 
is highly discontinuous (cf.\ Section~\ref{sec:9}). 

It seems that one possible way out of this dilemma is to find 
suitable $C^*$-algebras consisting of operators with a more regular 
behavior than $\CCR(\cH)$. For interesting recent results 
in this direction we refer to Georgescu's work 
\cite{Ge07}. Another step in this direction is the construction 
of a $C^*$-algebra $\cA$ for each countably dimensional symplectic space 
whose representations correspond to those representations of the 
corresponding Weyl relations which are continuous on each one-parameter 
group (\cite{GrNe09}). 
\end{ex}

\section{The Virasoro algebra and vector fields on $\bS^1$} \mlabel{sec:8}

In this section we discuss invariant cones 
in the Lie algebra $\cV(\bS^1)$ of smooth vector fields on the circle 
and its (up to isomorphy unique) non-trivial central extension 
$\vir$, the Virasoro algebra. For $\cV(\bS^1)$ we show that, up to sign, 
there is only one open invariant convex cone given by vector fields 
of the form $f \frac{\partial}{\partial\theta}$ with $f > 0$.  
As is well-known on the Lie algebra level, all unitary highest weight 
representations of $\cV(\bS^1)$, resp., the subalgebra of vector fields 
for which $f$ is a finite Fourier polynomial, are trivial. 
On the group level we show the closely related result that 
all semibounded unitary representations of $\Diff(\bS^1)_+$ are trivial. 
This is the main reason for the Virasoro algebra and the corresponding 
simply connected group $\Vir$ playing a more important role 
in mathematical physics than $\Diff(\bS^1)_+$ itself 
(cf.~\cite{SeG81}, \cite{Mick89}, \cite{Ot95}). 
For $\vir$ we prove a convexity theorem for 
adjoint and coadjoint orbits which provides complete 
information on invariant cones and permits us to 
determine the momentum sets of the unitary highest weight representations 
of $\Vir$ and to show that they are semibounded. 

\subsection{The invariant cones in $\cV(\bS^1)$}

We consider $\bS^1$ as the quotient $\R/2\pi\Z$ and identify 
smooth functions on $\bS^1$ with the corresponding 
$2\pi$-periodic smooth functions on $\R$, where the coordinate is denoted 
by $\theta$. Accordingly $\partial_\theta := \frac{d}{d\theta}$ 
is the vector field generating the rigid rotations of $\bS^1$. 
We write $G := \Diff(\bS^1)_+^{\rm op}$ 
for the group of orientation preserving 
diffeomorphisms endowed with the group structure defined by 
$\phi \cdot \psi := \psi \circ \phi$, so that 
$\g = \cV(\bS^1) = C^\infty(\bS^1) \partial_\theta$ 
is its Lie algebra (cf. Examples~\ref{ex:liegrp}(e)). 
 We represent orientation preserving 
diffeomorphisms of $\bS^1$ by smooth functions 
$\phi \: \R \to \R$ satisfying 
$\phi(\theta+2\pi) = \phi(\theta) + 2\pi$ for $\theta \in \R$ and $\phi' > 0$. 

In the following it will be convenient to consider 
the spaces 
$$\cF_s(\bS^1) := C^\infty(\bS^1) (d\theta)^s, \quad s \in \R, $$
of $s$-densities on $\bS^1$. Here $(d\theta)^s$ denotes the canonical 
section of the $s$-density bundle of $\bS^1$ and $G$ acts on 
$\cF_s(\bS^1)$ by pullbacks 
\begin{equation}
  \label{eq:diffdens}
\phi^*(u (d\theta)^s) 
= (u \circ \phi) (\phi^*d\theta)^s
= (u \circ \phi) (\phi')^s (d\theta)^s. 
\end{equation}
The corresponding derived action of the Lie algebra $\cV(\bS^1)$ 
is given by the Lie derivative 
\begin{equation}
  \label{eq:lieder}
\cL_{f\partial_\theta}(u(d\theta)^s) 
= (fu' + s f'u)(d\theta)^s. 
\end{equation}
The space $\cF_1$ is the space of $1$-forms and $\cF_{-1} \cong 
\cV(\bS^1)$ is the space of vector fields on which 
\eqref{eq:diffdens} describes the adjoint action. We have 
equivariant multiplication maps $\cF_s \times \cF_t \to \cF_{s+t}$, 
and an invariant integration map 
$$ I \: \cF_1 \to \R, \quad f d\theta \mapsto \int_0^{2\pi} f(\theta)\, d\theta,$$
which leads to an invariant pairing 
$\cF_{-1} \times \cF_2 \to \R,$ and hence to an equivariant embedding 
$$ \cF_2 \into \cV(\bS^1)' = \g'. $$
Its image is called the {\it smooth dual of $\g$}. Identifying it 
with $\cF_2$, the coadjoint action of $G$ on $\g'$ 
corresponds to the natural action on $\cF_2$. 

In view of \eqref{eq:diffdens}, the adjoint action 
clearly preserves the open cone 
$$ W_{\cV(\bS^1)} := \{ f \partial_\theta \: f > 0, f \in C^\infty(\bS^1) \} $$
of all vector fields corresponding to positive functions. 
Since every positive function $f$ with $\int_0^{2\pi} f\, d\theta = 2\pi$ 
arises as $\phi'$ for some $\phi \in G$, each $G$-orbit 
in $W_{\cV(\bS^1)}$ intersects the (maximal) abelian subalgebra 
$\ft := \R \partial_\theta$. We also have a projection map 
$$ p_\ft \: \g \to \ft, \quad 
f \partial_\theta \mapsto \frac{1}{2\pi}\int_0^{2\pi} f(s)\, ds \cdot  
\partial_\theta = \int_T \Ad(\phi)(f\partial_\theta)\, d\mu_T(\phi),$$
where $T := \exp(\ft)\cong \T$ is the group of rigid rotations 
and $\mu_T$ the normalized Haar measure on $T$. 

The following lemma provides some fine information on the convex geometry 
of adjoint orbits in $W_{\cV(\bS^1)}$. 

\begin{lem} \mlabel{lem:9.1} The function 
$\chi \: W_{\cV(\bS^1)} \to \R, f \partial_\theta \mapsto 
\chi(f) :=  \frac{1}{2\pi}\int_0^{2\pi} \frac{1}{f}\, d\theta$ 
has the following properties: 
\begin{description}
\item[\rm(i)] It is smooth, $G$-invariant and strictly convex. 
\item[\rm(ii)]  If the sequence $(f_n)$  in $W_{\cV(\bS^1)}$ converges to 
$f \in \partial W_{\cV(\bS^1)}$, then $\chi(f_n) \to \infty$. 
\item[\rm(iii)]  $\chi(f + g) \leq \chi(f)$ for $f,g \in W_{\cV(\bS^1)}$. 
\item[\rm(iv)] For each $c > 0$, the set 
$I_c := \{ x \in W_{\cV(\bS^1)} \: \chi(x) \leq \frac{1}{c}\}$ 
is an invariant closed convex subset of $\g$ with  
$\lim(I_c) = \oline{W_{\cV(\bS^1)}}$. Its boundary is a single 
orbit 
$$ \cO_{c\partial_\theta} = \partial I_c = \Ext(I_c) 
= \Big\{ f \in W_{\cV(\bS^1)} \: \chi(f) = \frac{1}{c}\Big\}$$ 
which coincides with its set of extreme points and satisfies 
$I_c = \oline{\conv}(\cO_{c\partial_\theta})$. 
\item[\rm(v)] $p_\ft(\cO_{c\partial_\theta}) 
= p_{\ft}(I_c) = [c,\infty[ \cdot \partial_\theta$ and 
the only inverse image of the ``minimal value'' $c\partial_\theta$ is the 
element $c\partial_\theta$ itself. 
\end{description}
\end{lem}

\begin{prf} (i)  The $G$-invariance of $\chi$ follows immediately 
from the Substitution Rule. 
The function $\chi$ is convex because inversion is a convex function on 
$\R_+^\times$ and integrals of convex functions are convex. 
It is smooth because on $W_{\cV(\bS^1)}$ pointwise inversion is a smooth 
operation (since $C^\infty(\bS^1,\R)$ is a real continuous inverse 
algebra; \cite{Gl02}), and integration 
is a continuous linear functional. 

To verify that $\chi$ is strictly convex, we observe that 
\begin{equation}
  \label{eq:secderiv}
(\partial_h \chi)(f) = -\frac{1}{2\pi}\int_0^{2\pi} \frac{h}{f^2}\, d\theta
\quad \mbox{ and } \quad 
(\partial_h^2 \chi)(f) = \frac{1}{\pi}\int_0^{2\pi} \frac{h^2}{f^3}\, d\theta,
\end{equation}
which is positive definite for each $f > 0$. 

(ii) Now we turn to the boundary behavior of $\chi$. Suppose that 
the sequence $(f_n)$ in $W_{\cV(\bS^1)}$ tends to a boundary point $f \in \partial W_{\cV(\bS^1)}$. 
Then $f(\theta_0) = 0$ for some $\theta_0 \in [0,2\pi[$, and 
$f \geq 0$ implies that we also have $f'(\theta_0) = 0$, hence 
$f(\theta) \leq C (\theta-\theta_0)^2$ in a compact $\delta$-neighborhood $U$ of $\theta_0$. 
Given $\eps > 0$, we eventually have 
$f_n \leq 2C(\theta-\theta_0)^2 + \eps$ on $U$ (here we use $C^2$-convergence), 
and therefore 
$$2\pi \chi(f_n) 
\geq \int_{\theta_0-\delta}^{\theta_0+\delta} \frac{1}{f_n(\theta)}\, d\theta 
\geq \int_{-\delta}^{\delta} \frac{1}{2C \theta^2 + \eps}\, d\theta. $$
Since  
$\int_{-\delta}^{\delta} \frac{1}{\theta^2}\, d\theta = \infty,$
the Monotone Convergence Theorem implies that 
$$ \lim_{\eps \to 0} \int_{-\delta}^{\delta} \frac{1}{2C\theta^2+\eps}\, d\theta 
= \infty, $$
and therefore that $\chi(f_n)\to \infty$. It also follows that 
$\int_0^{2\pi} \frac{1}{f(\theta)}\, d\theta = \infty$. 

(iii) If $f,g \in W_{\cV(\bS^1)}$, then $\frac{1}{f+g} \leq \frac{1}{f}$ 
implies the assertion. 

(iv) From  (ii) we derive that $I_c$ is closed in $\g$. 
Its invariance follows from the invariance of $\chi$ and 
its convexity from the convexity of $\chi$. 
The boundary of $I_c$ is a level set of $\chi$, and 
since every orbit in $W_{\cV(\bS^1)}$ meets $\R_+ \partial_\theta$ in a unique point 
and $\chi \res_{\R_+^\times \partial_\theta}$ is injective, it follows that 
$\partial I_c = \cO_{c \partial_\theta}$. 
The fact that $\chi$ is strictly convex further implies that 
$\partial I_c \subeq \Ext(I_c)$, and since the converse inclusion is trivial, 
equality follows. 

From (iii) we derive that $W_{\cV(\bS^1)} \subeq \lim(I_c)$, 
so that the closedness of $\lim(I_c)$ 
(Lemma~\ref{lem:limcone}(i)) implies that 
$\oline {W_{\cV(\bS^1)}} \subeq\lim(I_c)$. Now equality follows from 
$\lim(I_c) \subeq \lim(W_{\cV(\bS^1)}) = \oline{W_{\cV(\bS^1)}}$ 
(Lemma~\ref{lem:limcone}(ii)).

To see that $I_c$ coincides with the closed convex hull 
$D := \oline\conv(\cO_{c \partial_\theta})$, 
we first observe that we trivially 
have $D \subeq I_c$. Next we note that 
$\partial_\theta$ is contained in the $3$-dimensional subalgebra 
$$\fs := \Spann \{ 1, \cos(\theta), \sin(\theta) \} \partial_\theta 
\cong \fsl_2(\R) $$   
corresponding to the action of $\SL_2(\R)$ on 
$\bP_1(\R) \cong \bS^1$. In $\fs$ 
the element $\partial_\theta$ corresponds to the matrix 
$$u =
\begin{pmatrix}
\phantom{-}0 & 1 \\
-1 & 0
\end{pmatrix}, $$
so that Example~\ref{ex:6.5} implies that the corresponding 
group $S := \la \exp \fs \ra \cong \PSL_2(\R)$ satisfies 
$p_\ft(\Ad(S)\partial_\theta) = [1,\infty[ \cdot \partial_\theta$ 
(cf.\ Example~\ref{ex:vectcirc}) and from that we derive in 
particular that 
\begin{equation}
  \label{eq:orb-proj}
[c,\infty[ \cdot \partial_\theta \subeq p_\ft(\cO_{c \partial_\theta})
\subeq D, 
\end{equation}
so that $I_c = \Ad(G)([c,\infty[\cdot \partial_\theta)$ leads to 
$I_c \subeq D$. 

(v) We have already seen in \eqref{eq:secderiv} above that 
$ (\partial_h \chi)(c) = - \frac{1}{2\pi c^2} \int_0^{2\pi} h\, d\theta$, so that 
$p_\ft^{-1}(c\partial_\theta)$ is a tangent hyperplane of the strictly convex 
set $I_c$. This implies that $c \partial_\theta$ is the unique minimum of the 
linear functional $-\dd\chi(c)$ on $I_c$ and hence that 
$p_{\ft}^{-1}(c\partial_\theta) \cap I_c = \{c\partial_\theta\}$. 

We also conclude that 
$p_\ft(I_c) \subeq [c,\infty[ \cdot \partial_\theta$, 
so that (v) follows from \eqref{eq:orb-proj}. 
\end{prf} 

\begin{rem} \mlabel{rem:9.13} 
The topological dual $\cV(\bS^1)'$ of $\cV(\bS^1) = C^\infty(\bS^1) 
\partial_\theta$ naturally identifies with the space 
of distributions on $\bS^1$. 
Then each $\lambda \in W_{\cV(\bS^1)}^\star$ is a distribution satisfying 
$\lambda(f) \geq 0$ for $f \geq 0,$
and this implies that $\lambda$ extends continuously from 
$C^\infty(\bS^1)$ to the Banach space $C(\bS^1)$ and thus defines a 
(finite) positive Radon measure 
on $\bS^1$ (cf.\ \cite[Ch.~I, \S 4, Thm.~V]{Sw73}). 
This shows that the functionals in $W_{\cV(\bS^1)}^\star$ satisfy a strong 
regularity condition. 
\end{rem}

\begin{thm}{\rm(Classification of open invariant cones in 
$\cV(\bS^1)$)} \mlabel{thm:classcon}
The two open cones $\pm W_{\cV(\bS^1)}$ are the only 
non-empty proper open invariant cones  in $\cV(\bS^1)$.
\end{thm}

\begin{prf} Let 
$C \subeq \cV(\bS^1)$ be a non-emptry open invariant cone. Then $C$  
is in particular invariant under the adjoint action 
of the rotation group 
$T \cong \T$, generated by the vector field 
$d := \partial_\theta$. Averaging over $T$, we see that 
either $d$ or $-d$ is contained in $C$ 
(Proposition~\ref{prop:project}), which leads to 
$W_{\cV(\bS^1)} \subeq C$ or $-W_{\cV(\bS^1)} \subeq C$. 

Let us assume that $W_{\cV(\bS^1)} \subeq C$. 
Now we apply the same argument to the dual cone $C^\star \subeq W_{\cV(\bS^1)}^\star$ (Proposition~\ref{prop:project}(b)). 
If $C$ is proper, then $C^\star$ contains a non-zero functional 
$\lambda$, and then $d \in C$ leads to $\lambda(d) > 0$. 
Averaging over $T$ now leads to a $T$-invariant 
functional $\lambda^* \in C^\star \subeq W_{\cV(\bS^1)}^\star$ satisfying 
$\lambda^*(d) = \lambda(d) > 0$. We conclude that $\lambda^*$ is a positive 
multiple of the invariant measure $\mu = d \theta$ 
on 
$\bS^1 \cong \R/2\pi\Z$ (Remark~\ref{rem:9.13}), 
which corresponds to the constant function $1 \in \cF_2$, 
so that
$$ \Ad^*(\phi) \mu =  (\phi')^2 \mu $$
follows from \eqref{eq:diffdens}. 
For any $f \partial_\theta \in C$ we now find 
$\int_0^{2\pi} \psi(\theta)^2 f(\theta)\, d\theta \geq 0$
for each positive function $\psi$ with $\int_0^{2\pi}  \psi(\theta)\, d\theta 
=2\pi$  (these are precisely the functions occurring as $\phi'$ for an 
orientation preserving diffeomorphism), and this implies the 
corresponding relation for all non-negative functions 
$\psi$ with the integral $2\pi$. 
If $f(\theta_0) < 0$, then we may choose $\psi$ supported by 
the set $\{ f < 0\}$ and obtain a contradiction. 
This proves that $C \subeq \oline{W_{\cV(\bS^1)}}$ and hence that $C = W_{\cV(\bS^1)}$ 
because $W_{\cV(\bS^1)} = (\oline{W_{\cV(\bS^1)}})^0$ 
(Lemma~\ref{lem:bou}). 
\end{prf}

\begin{prop} \mlabel{prop:semieq-vect} 
 For $\lambda \in \cV(\bS^1)'$ the following are equivalent
  \begin{description}
  \item[\rm(i)] $\cO_\lambda$ is semi-equicontinuous. 
  \item[\rm(ii)] $\inf \cO_\lambda(\partial_\theta) > -\infty$ or 
$\sup \cO_\lambda(\partial_\theta) < \infty$. 
  \item[\rm(iii)] $\lambda \in W_{\cV(\bS^1)}^\star \cup - W_{\cV(\bS^1)}^\star$. 
  \end{description}
In particular, $\inf \cO_\lambda(\partial_\theta) > -\infty$ implies 
$\cO_\lambda(\partial_\theta) \geq 0$. 
\end{prop}

\begin{prf} We have to verify the assumptions of Lemma~\ref{lem:doubcone} 
for $d := \partial_\theta$. 
First, $W_{\cV(\bS^1)}$ is the open convex cone generated by 
the orbit $\cO_d$. Second, the projection 
$p_\ft \: \g \to \R d = \ft$ is the fixed point projection for the 
adjoint action of the circle group $T$, hence preserves 
open and closed convex subsets (Proposition~\ref{prop:project}). 
Finally, we recall from Lemma~\ref{lem:9.1}(v) 
that $p_\ft(\cO_d) = [1,\infty[d$ is unbounded. Now the assertion 
follows from Lemma~\ref{lem:doubcone}. 
\end{prf}

\begin{thm} \mlabel{thm:8.6} 
The quotient $M := G/T$ of the group 
$G = \Diff(\bS^1)_+$ by the subgroup $T$ of rigid rotations 
carries the structure of a complex Fr\'echet manifold on which 
$G$ acts smoothly by holomorphic maps. 
Here the tangent space $T_{x_0}(G/T)$ in the base point 
$x_0 = \1 T$ is canonically identified with 
$\g_\C/\fp$, where 
$$ \fp = \Big\{ f \partial_\theta \in C^\infty(\bS^1,\C) \partial_\theta \: 
f = \sum_{n \leq 0}^\infty a_n e^{in\theta}, a_n\in \C \Big\}. $$ 
For each element $x \in W_{\cV(\bS^1)}$, the flow on $M$ defined by 
$(t,gT) \mapsto \exp(tx)g T$ extends to a smooth flow on the upper 
half plane $\C_+$ which is holomorphic on 
$\C_+^0 \times M$. 
\end{thm}

\begin{prf} The complex structure on $M$ has been discovered 
by Kirillov and Yuriev (cf.\ \cite{Ki87}, \cite{KY87}, \cite{Ov01}; 
see in particular \cite{Le95} for rigorous arguments concerning the 
complex Fr\'echet manifold structure). 
Complex structures always 
come in pairs, and we therefore write 
$\oline M$ for the same smooth manifold, endowed with the 
opposite complex structure. 
  
According to \cite{KY87}, the complex manifold structure on $\oline M$ 
can be obtained by identifying it with the space 
$\cF_{\rm reg}$ of {\it normalized regular univalent functions} 
$f \: \cD \to \C$, where $\cD$ is the open unit disc in $\C$. 
These are the injective holomorphic maps $f  \: \cD \to \C$ 
extending smoothly to the 
closure of $\cD$ and satisfying $f(0) = 0$ and $f'(0) = 1$ 
(cf.\ \cite[Sect.~6.5.6]{GR07} for a detailed discussion). 
Here the complex structure is determined by its action on 
$T_{x_0}(G/T) \cong T_{\id_\cD}(\cF_{\rm reg}) \cong \g/\ft$: 
$$ I[(e^{in\theta}+ e^{-in\theta})\partial_\theta]
= [(-i e^{in\theta}+ i e^{-in\theta})\partial_\theta]
\quad \mbox{ for } \quad n > 0. $$

Since each $x \in W_{\cV(\bS^1)}$ is conjugate to a multiple of $\partial_\theta$, 
it suffices to assume that $x = \partial_\theta$ is the generator 
of the rigid rotations of $\bS^1$. In this case the action 
of the one-parameter group $T = \exp_G(\R x)$ on $\cF_{\rm reg}$ 
is given explicitly by 
$$ R_\alpha f = R_\alpha \circ f \circ R_{-\alpha}, \quad \mbox{where} 
\quad R_\alpha z = e^{i\alpha}z, \alpha \in \R $$
(cf.\ \cite[Prop.~6.5.14]{GR07}). 
If $\Im \alpha \leq 0$, then $|e^{-i\alpha}| \leq 1$, so that 
$R_\alpha(\cD) \subeq \cD$ implies that 
$R_\alpha f$ can still be defined as above, is continuous on 
$\C_- := -\C_+$ and depends holomorphically on $\alpha$ for $\Im \alpha < 0$. 
This implies the holomorphic extension to $\C_-$ for the 
complex manifold $\oline M$. For the manifold $M$ we therefore 
obtain an extension to $\C_+$. 
\end{prf}

The holomorphic extension of actions of one-parameter groups 
on $G/T$ can be carried much further. As shown by Neretin in \cite{Ner90}, 
one even has a ``complex semigroup'' containing $G$ in its boundary 
which acts on $G/T$. 

Below we shall use the preceding theorem to identify the momentum 
sets for the unitary highest weight representations of $\vir_\C$. 

\begin{defn} For the following we recall some algebraic 
aspects of $\cV(\bS^1)$. 
In the complexification 
$\cV(\bS^1)_\C$, we consider the elements 
\begin{equation}
  \label{eq:wittgen}
d_n :=  i e^{i n\theta} \partial_\theta, \quad 
n \in \Z, 
\end{equation}
satisfying the commutation relations 
$$ [d_n, d_m] = (n-m) d_{n+m}. $$
The standard involution on this Lie algebra 
is given by $(f \partial_\theta)^* = -\oline f \partial_\theta$, 
so that $x^* = - x$ describes the elements of $\cV(\bS^1)$. 
Note that $d_n^* = d_{-n}$ and in particular 
$d_0^* = d_0$, so that $d_0 = i \partial_\theta$ 
is a hermitian element (cf.\ \cite[p.~9]{KR87}). 
\end{defn}

\begin{thm} \mlabel{thm:sembotriv} 
All semibounded unitary representations 
of the group $\Diff(\bS^1)_+$ are trivial. 
\end{thm} 

\begin{prf} Let $(\pi, \cH)$ be a semibounded unitary representation 
of $G = \Diff(\bS^1)_+^{\rm op}$. Then $W_\pi \not=\eset$, and 
in view of Theorem~\ref{thm:classcon}, we may w.l.o.g.\ 
assume that $- \partial_\theta \in -W_{\cV(\bS^1)} \subeq W_\pi$, 
so that the spectrum of the image of 
$d_0 := i \partial_\theta \in \g_\C$ under the derived representation 
is bounded from below. In view of $\exp_G(2\pi i d_0) = \1$, it is contained 
in $\Z$ and Proposition~\ref{prop:maxsmoothvect} implies the 
existence of a smooth unit vector $v \in\cH^\infty$ which is an 
eigenvector for the minimal  eigenvalue $h$ of $d_0$. 

Now the relation $[d_0, d_n] = - n d_n$ implies 
that $d_n.v = 0$ for each $n > 0$. 
For $n > 0$ we then obtain
\begin{align*}
\la d_{-n}v, d_{-n}v\ra 
&=  \la d_{-n}^* d_{-n}v, v\ra 
=  \la d_n d_{-n}v, v\ra 
=  \la [d_n, d_{-n}]v, v\ra \\
&=  2n \la d_0 v, v\ra 
=  2n h.  
\end{align*}
This implies in particular that $h \geq 0$ and that 
$h = 0$ implies $d_n v= 0$ for each $n \in \Z$. 
Now an easy direct calculation leads to 
$$ 0 \leq \det\pmat{ 
\la d_{-2n}v, d_{-2n}v \ra & 
\la d_{-n}^2v, d_{-2n}v \ra \\ 
\la d_{-2n}v, d_{-n}^2v \ra & 
\la d_{-n}^2v, d_{-n}^2v \ra} 
= 4 n^3 h^2 (8h - 5n)  $$
(\cite[p.~90]{KR87}; see also \cite{GO86}).
If $h \not=0$, this expression is negative for sufficiently large $n$, 
so that we must have $h = 0$. This means that 
$\g.v = \{0\}$, and hence that $v \in \cH^G$ is a fixed vector 
(cf.~\cite[Rem.~II.3.7]{Ne06}). 

The preceding argument implies that each semibounded 
unitary representation $(\pi, \cH)$ of $G$ on a non-zero 
Hilbert space satisfies $\cH^G \not=\{0\}$. 
Applying this to the representation on the invariant subspace 
$(\cH^G)^\bot$, which is also semibounded, we find that this 
space is trivial, and hence that $\cH = \cH^G$, i.e., 
the representation is trivial. 
\end{prf}

\begin{rem} A smooth unitary representation  
$(\pi, \cH)$ of $\Diff(\bS^1)_+$ is said to be a {\it positive 
energy representation} if the operator 
$-i\dd\pi(\partial_\theta)$ has non-negative spectrum. 
This means that $\partial_\theta \in I_\pi^\star$, so that 
$W_{\cV(\bS^1)} \subeq I_\pi^\star$ leads to $W_{\cV(\bS^1)} \subeq C_\pi$, 
and therefore $\pi$ is semibounded. 
Hence the preceding theorem implies in particular that 
all positive energy representations of $\Diff(\bS^1)_+$ are 
trivial. 
\end{rem}

\begin{prob} It would be nice to have an analog 
of Theorem~\ref{thm:sembotriv} for the universal covering 
group $\tilde G$ of $G = \Diff(\bS^1)_+$, which has the 
fundamental group $\pi_1(G) = \Z$. 
Then $\partial_\theta$ generates a subgroup 
$\tilde T \subeq \tilde G$ isomorphic to $\R$. Since this group 
is non-compact, we cannot expect it to have eigenvectors, so that 
the argument in the proof of Theorem~\ref{thm:sembotriv} does 
not apply. 

What we would need in this context is a suitable direct 
integral decomposition with respect to the subgroup 
$Z := Z(\tilde G) \cong \Z$. If $(\pi, \cH)$ is a semibounded 
representation of $\tilde G$ with $\pi(Z) \subeq \T \1$, i.e., 
$\pi$ has a central character (which is a consequence of 
Schur's Lemma if $\pi$ is irreducible), then 
$\Spec(\dd\pi(d_0)) \subeq \lambda + \Z$ for some 
$\lambda \in \R$, and the argument from above applies. 
This proves that all irreducible semibounded representations 
of $\tilde G$ are trivial. 

We expect that a general semibounded representation of $\tilde G$ 
has a direct integral decomposition 
$\cH \cong \int^\oplus_{\hat Z} \cH_\chi\, d\mu(\chi)$ 
with respect to some measure $\mu$ on $\hat Z \cong \T$ 
and that the semiboundedness of $\pi$ implies that all representations 
$\pi_\chi$ on the spaces $\cH_\chi$ with central character 
$\chi$ are semibounded, hence trivial, and this would imply that 
all semibounded representations of $\tilde G$ are trivial. 
\end{prob}

\subsection{Invariant cones in the Virasoro algebra} 

In the analytic context, the Virasoro algebra 
is usually defined as the 
central extension $\vir = \R \oplus_{\omega_{GF}} \cV(\bS^1)$ 
defined by the {\it Gelfand--Fuchs} cocycle 
\begin{equation}
  \label{eq:vircocyc}
\omega_{GF}(f\partial_\theta, g\partial_\theta) 
:= \int_0^{2\pi} f'g''\, d\theta 
= \frac{1}{2}\int_0^{2\pi} f'g'' - f''g'\, d\theta 
= \int_0^{2\pi} f'''g\, d\theta. 
\end{equation}
In many situations, the cohomologous cocycle 
\begin{eqnarray}
  \label{eq:vircocyc2}
\omega(f\partial_\theta, g\partial_\theta) 
&:=& \int_0^{2\pi} (f'''+f')g\, d\theta 
= \omega_{GF}(f\partial_\theta, g\partial_\theta)  
- \shalf \int_0^{2\pi} fg' -f'g\, d\theta \notag\\
&=& \omega_{GF}(f\partial_\theta, g\partial_\theta)  
- \shalf \lambda([f\partial_\theta, g\partial_\theta]),  
\end{eqnarray}
with $\lambda(f\partial_\theta) = \int_0^{2\pi}f\, d\theta$, 
turns out to be more convenient. 

\begin{rem} \mlabel{rem:vir-rel}
On the generators $d_n =  i e^{i n\theta} \partial_\theta 
\in \cV(\bS^1)_\C$ from \eqref{eq:wittgen} we have 
\begin{equation}
  \label{eq:vircocgen}
\omega(d_n, d_{-n}) = 2\pi i (n^3-n). 
\end{equation}
With the central element $\hat c := (24 \pi i,0) \in i \vir \subeq \vir_\C$, 
we thus obtain the relation 
$$ [d_n, d_m] = (n-m) d_{n+m} + \delta_{n,-m} \frac{n^3 -n}{12}\hat c $$ 
if we identify $d_n$ with the corresponding element 
$(0,d_n) \in \vir_\C$ (\cite[p.~9]{KR87}).

Since we shall need them in the following, we record some related formulas. 
First we observe that 
\begin{equation}
  \label{eq:vir-cartan}
\ft := \R c + \R d \quad \mbox{ with } \quad c := (1,0), 
d := (0,\partial_\theta) 
\end{equation}
is a maximal abelian subalgebra of $\vir$. 
The relation $[d_0, d_n] = - n d_n$ implies that 
$d_n \in \vir_\C$ is a root vector 
for the root $\alpha_n \in \ft_\C^*$ defined by 
$$ \alpha_n(\hat c) = 0 \quad \mbox{ and } \quad 
\alpha_n(d_0) = -n. $$
In particular, 
\begin{equation}
  \label{eq:rootpos}
\alpha_n([d_n, d_n^*]) = \alpha_n([d_n, d_{-n}]) 
= 2n \alpha_n(d_0) = - 2n^2 < 0 \quad \mbox{ for } \quad n \not=0.
\end{equation}
We also observe that 
\begin{eqnarray} \label{eq:rootbrack}
[d_n^*, d_n] = [d_{-n}, d_n] 
&= (\omega(d_{-n}, d_n), -2n d_0)
= -i (2\pi(n^3 - n), 2n \partial_\theta)\notag \\
= - 2i n (\pi(n^2 - 1), \partial_\theta). 
\end{eqnarray}
\end{rem} 

In view of the $G$-invariant pairing of the space $\cF_2$ of $2$-densities 
with vector fields, the $1$-cocycle $\g \to \g', x \mapsto i_x\omega$ 
corresponds to the $1$-cocycle 
$$ \g \to \cF_2, \quad f \partial_\theta \mapsto 
 (f''' + f') (d\theta)^2. $$
To obtain a formula for the adjoint action of $G = \Diff(\bS^1)_+^{\rm op}$ 
on $\vir$, we therefore need a group cocycle $G \to \cF_2$ integrating 
this $1$-cocycle (Remark~\ref{rem:centext}). 

\begin{defn}
The {\it Schwarzian derivative} 
$$ S(\phi) 
:= \frac{\phi' \phi''' - \frac{3}{2} (\phi'')^2}{(\phi')^2} 
= \frac{\phi'''}{\phi'} - \frac{3}{2} \Big(\frac{\phi''}{\phi'}\Big)^2 $$ 
assigns to $\phi \in G$ a 
$2\pi$-periodic smooth function. It satisfies the cocycle identity 
$$ S(\phi \circ \psi) = \big(S(\phi) \circ \psi\big)\cdot
 (\psi')^2 + S(\psi) $$
(cf.\ \cite{Ov01}), 
which means that it defines an $\cF_2$-valued  $1$-cocycle 
on $G$. 
\end{defn}
We easily derive that 
$T_{\id}(S)(f) = f''',$ and therefore the {\it modified Schwarzian 
derivative} 
$$ \tilde S(\phi) := S(\phi) + \shalf((\phi')^2 - 1)
= S(\phi) + \phi.\shalf - \shalf $$
is a cohomologous $1$-cocycle with 
$T_{\id}(\tilde S)(f) = f''' + f'$ (cf.\ \cite[Sect.~7]{SeG81}). Therefore 
Remark~\ref{rem:centext} implies that the coadjoint action 
of $G$ on the smooth dual 
$\R \times \cF_2 \cong \R \times C^\infty(\bS^1)$ 
of $\vir$ is given by  
\begin{equation}
  \label{eq:centcoadvir}
\Ad^*_{\phi}(a, u) 
= \big(a, (u\circ \phi) (\phi')^2 - a \tilde S(\phi)\big), 
\end{equation}
whereas, in view of \eqref{eq:diffdens}, the adjoint action on $\vir 
= \R \oplus_\omega \cV(\bS^1) \cong \R \times C^\infty(\bS^1)$ 
is given by 
$$ \Ad_\phi(z,f) 
= \Big(z - \int_0^{2\pi} f \tilde S(\phi^{-1})\, d\theta, (f \circ \phi) 
\cdot (\phi')^{-1}\Big). $$

We are especially interested in the adjoint action 
on the open convex cone 
$$ W_{\rm max} := \{ (z,f) \in \vir \: f > 0\}, $$
which is the inverse image of the positive invariant cone $W_{\cV(\bS^1)} 
\subeq \cV(\bS^1)$ 
under the quotient map $\vir \to \cV(\bS^1), (z,f) \mapsto f$. 
From the corresponding results for $W_{\cV(\bS^1)}$, we derive immediately 
that each orbit in $W_{\rm max}$ intersects $\ft$. 

\begin{prop} \mlabel{prop:8.11} 
For each $(z,f) \in W_{\rm max}$, the $G$-orbit of $(z,f)$ 
meets $\ft$ in a unique element $(\beta(z,f),\alpha(z,f))$, given by 
$$ \alpha(z,f) := \frac{1}{\chi(f)}  \quad \mbox{ and } \quad 
\beta(z,f) := z  - \int_0^{2\pi} \frac{(f')^2}{2f}\, d\theta  
+ \shalf \int_0^{2\pi} f\, d\theta - \frac{\pi}{\chi(f)}. $$  
For $G = \Diff(\bS^1)_+$ and $T = \exp(\R \partial_\theta)$, the 
orbit map induces a diffeomorphism 
$$ \Gamma\: G/T \times (W_{\rm max} \cap \ft) \to W_{\rm max}, \quad 
(\phi T,(\beta, \alpha)) \mapsto \Ad_{\phi}(\beta, \alpha). $$
\end{prop}

\begin{prf} Let $(z,f) \in W_{\rm max}$ and assume that 
$\Ad_{\phi^{-1}}(z,f) = (\beta, \alpha) \in \ft$, i.e., 
$$ (\beta, \alpha) 
= \Big(z - \int_0^{2\pi} f 
\tilde S(\phi)\, d\theta, (f\cdot \phi')\circ \phi^{-1}
\Big). $$
Then $\phi' = \frac{\alpha}{f}$ leads to $\phi'' = - \alpha \frac{f'}{f^2}$ 
and 
$$ \phi''' 
= - \alpha \frac{f''}{f^2} + 2\alpha \frac{(f')^2}{f^3} 
= \frac{\alpha}{f^3}(2 (f')^2 - ff''). $$
This leads to 
\begin{align*}
S(\phi) = \frac{\phi'''}{\phi'} - \frac{3}{2}\Big(\frac{\phi''}{\phi'}\Big)^2 
&= \frac{1}{f^2}(2(f')^2 - ff'') - \frac{3}{2}\frac{(f')^2}{f^2} 
= \frac{1}{2f^2}((f')^2 - 2ff'')\\
&= \frac{(f')^2}{2f^2}- \frac{f''}{f}, 
\end{align*}
so that we obtain with $\int_0^{2\pi} f''\, d\theta 
= f'(2\pi) - f'(0) = 0$ the relation 
$$ \int_0^{2\pi} f S(\phi)\, d\theta = 
\int_0^{2\pi} \frac{(f')^2}{2f}\, d\theta. $$
Next we use $\phi' = \frac{\alpha}{f}$ to obtain 
$$ \shalf \int_0^{2\pi} f ((\phi')^2 -1)\, d\theta 
= \shalf \int_0^{2\pi} \alpha \phi' - f \, d\theta 
= \pi \alpha - \shalf \int_0^{2\pi} f \, d\theta. $$
Combining all that, we get 
\begin{align*}
\beta 
&= z - \int_0^{2\pi} \frac{(f')^2}{2f}\, d\theta 
+ \shalf \int_0^{2\pi} f \, d\theta - \pi \alpha.
\end{align*}
We also obtain from  $\phi' = \frac{\alpha}{f}$ the relation 
$\chi(f) 
= \frac{1}{2\pi}\int_0^{2\pi} \frac{1}{f}\, d\theta 
= \frac{1}{\alpha},$
so that $\alpha = \frac{1}{\chi(f)}$. 
This proves the first assertion. 

Since $T$ fixes the subalgebra $\ft$ pointwise, 
$\Gamma$ is a well-defined smooth map. As a manifold, we may identify 
$G/T$ with the set 
$$\{ \phi' \: \phi \in G \} = \Big\{ h \in C^\infty(\bS^1) \: h > 0, 
\int_0^{2\pi} h\, d\theta = 2 \pi\Big\}. $$
As we have seen above, the inverse of $\Gamma$ is given by 
$$ \Gamma^{-1}(z,f) = \big(\phi, (\beta(z,f), \alpha(z,f))\big), $$
where $\phi' = \frac{\alpha}{f} = \frac{1}{f \chi(f)}$, 
and this map is also smooth. 
Therefore $\Gamma$ is a diffeomorphism.   
\end{prf}

The function $\beta$ is rather complicated, but it can 
be analyzed to some extent as follows. First we observe that 
for the probability measure $\mu = \frac{1}{2\pi} \, d\theta$ 
Jensen's inequality and the convexity of the function $\frac{1}{x}$ 
on the positive half line imply that 
$$ \frac{1}{\int_0^{2\pi} f\, d\mu} \leq \int_0^{2\pi} \frac{1}{f}\, d\mu 
= \chi(f) = \frac{1}{\alpha}, $$
which implies that 
$\int_0^{2\pi} f\, d\theta \geq 2\pi \alpha.$
Therefore the sign of $\beta(0,f)$ is not clear at all, 
but we shall see below that $\beta(0,f) \leq 0$. 

\begin{lem}
  \mlabel{lem:beta-concave}
 The function $\beta$ is concave. 
\end{lem}

\begin{prf} To show that $\beta$ is concave, we have 
to verify that $\partial_h^2\beta \leq 0$ in each point of $W_{\rm max}$. 
Since $\beta$ is $\Ad(G)$-invariant and each orbit meets $\ft$, 
it suffices to verify this in points $(z,f)$, where 
$f$ is  constant, so that 
$$ (\partial_h^2 \beta)(z,f) 
= - f^{-1} \int_0^{2\pi}(h')^2\, d\theta - \pi \partial_h^2(\chi^{-1})(f). $$
Further, 
$$ \partial_h(\chi^{-1}) = - \chi^{-2} \partial_h \chi 
\quad \mbox{ and } \quad 
\partial_h^2(\chi^{-1}) = - \chi^{-2} (\partial_h^2 \chi) 
+ 2 \chi^{-3} (\partial_h \chi)^2, $$
so that we obtain with $\chi(f) = f^{-1}$ ($f$ is constant) and 
the formulas \eqref{eq:secderiv} in the 
proof of Lemma~\ref{lem:9.1} the relations 
$$ (\partial_h \chi)(f) = -\frac{1}{2\pi f^2}\int_0^{2\pi} h\, d\theta
\quad \mbox{ and } \quad 
 (\partial_h^2 \chi)(f) = \frac{1}{\pi f^3}\int_0^{2\pi} h^2\, d\theta. $$ 
This leads further to 
\begin{align*}
\partial_h^2(\chi^{-1})(f)
&= - \frac{f^2}{\pi f^3}\int_0^{2\pi} h^2\, d\theta 
+ 2 f^3 \frac{1}{4\pi^2f^4}\Big(\int_0^{2\pi} h\, d\theta\Big)^2 \\
&= \frac{1}{\pi f}\Big(-\int_0^{2\pi} h^2\, d\theta 
+ \frac{1}{2\pi}\Big(\int_0^{2\pi} h\, d\theta\Big)^2\Big). 
\end{align*}
Putting everything together, we arrive at 
$$ f \cdot (\partial_h^2 \beta)(z,f) 
= - \int_0^{2\pi} (h')^2\, d\theta  
+ \int_0^{2\pi} h^2\, d\theta 
- \frac{1}{2\pi}\Big(\int_0^{2\pi} h\, d\theta\Big)^2. $$
We thus obtain a rotation invariant 
quadratic form on $C^\infty(\bS^1)$, so that it 
is diagonal with respect to 
Fourier expansion. Evaluating it in the basis functions 
$\cos(n\theta)$ and $\sin(n\theta)$ immediately shows that it is 
negative semidefinite. 
\end{prf}

\begin{thm} {\rm(Convexity Theorem for adjoint orbits of $\vir$)} 
\mlabel{thm:concoad-vir}
For each $x \in \R c + \R_+ \partial_\theta \subeq \ft,$ we have 
$$ p_\ft(\cO_x) \subeq x + C_+ \quad \mbox{ for } \quad 
C_+ := \R_+ c + \R_+ \partial_\theta. $$
If $x$ is not central, then we even have the equality 
$$ p_\ft(\oline{\conv}(\cO_x)) = x + C_+. $$
\end{thm}

\begin{prf} By continuity of the projection $p_\ft$, it suffices 
to assume that $x = (\beta_0, \alpha_0)$ with a constant $\alpha_0 > 0$, so that 
$x \in W_{\rm max} \cap \ft$. Then $\alpha(x) =\alpha_0$ and 
$\beta(x) = \beta_0$. Further, 
$p_\ft(\cO_x) \subeq \oline{\conv}(\cO_x)$ 
(Proposition~\ref{prop:project}), so that the 
convexity of the functions $\chi = \alpha^{-1}$ and $-\beta$ 
(Lemma~\ref{lem:beta-concave}) implies that 
for $(\tilde\beta, \tilde\alpha) \in p_\ft(\cO_x)$, we have 
$$ \tilde\alpha \geq \alpha_0 \quad \mbox{ and } \quad 
\tilde\beta \geq \beta_0. $$
This means that $p_\ft(\cO_x) \subeq x + C_+$. 

Now we assume that $x = \beta c + \alpha \partial_\theta$ with 
$\alpha > 0$. In view of \eqref{eq:rootpos} 
in Remark~\ref{rem:vir-rel}, 
Proposition~\ref{prop:orb-pro} implies that 
$$ p_\ft(\cO_x) 
\supeq x + \R^+ \alpha_n(x) [d_n^*, d_n]
=  x + \R^+ \alpha_n(-i x) i [d_n^*, d_n]. $$
Further $d_0 = i \partial_\theta$ leads to 
$\alpha_n(-ix) > 0$, so that 
$$ p_\ft(\cO_x) \supeq  x + \R^+ i [d_n^*, d_n], \quad n \in \N. $$
Next we recall from \eqref{eq:rootbrack} that 
$$ i [d_n^*, d_n] = 2n (\pi(n^2 - 1), \partial_\theta) 
= 2n (\pi(n^2-1)c + d). $$
For $n = 1$ we obtain $2(0,\partial_\theta)$, and for $n \to \infty$ 
we have positive multiples of 
$(1, \frac{1}{\pi(n^2-1)}\partial_\theta) \to c$, 
so that the closed convex cone generated by the elements 
$i[d_n^*, d_n]$ is $C_+ = \R_+ c + \R_+ d$. 
This proves that 
$$ p_\ft(\oline{\conv}(\cO_x))
= \oline{\conv}(\cO_x) \cap \ft \supeq x + C_+, $$
and our proof is complete. 
\end{prf} 

Note that the following theorem can not be derived from 
the ``general'' Lemma~\ref{lem:doubcone} because 
$\ft$ is $2$-dimensional. 

\begin{thm} {\rm(Classification of open invariant cones in $\vir$)} 
\mlabel{thm:classinvcon-vir}
The following statements classify the open invariant convex cones 
in $\vir$: 
\begin{description}
\item[\rm(i)] Each proper open invariant convex cone in $\vir$ 
is either contained in $W_{\rm max}$ or $-W_{\rm max}$. 
\item[\rm(ii)] Each proper open invariant convex cone $W \subeq \vir$ 
is uniquely determined by 
$C := W \cap \ft$ via $W = \Ad(G)C$. 
\item[\rm(iii)] Let $C_{\rm max} := W_{\rm max} \cap \ft$. Then an 
open convex cone $C \subeq C_{\rm max}$ 
is the trace of an invariant open convex cone if and only if $C_+^0 \subeq C$. 
\item[\rm(iv)] If $W_{\rm min}$ is the open invariant cone corresponding 
to $C_{\rm min} := C_+^0$, then each open invariant convex 
cone $W \subeq W_{\rm max}$ contains $W_{\rm min}$. 
\item[\rm(v)] An open invariant convex cone $W$ is pointed if and only if $C$ is 
pointed. In particular, $W_{\rm min}$ is pointed. 
\end{description}
\end{thm}

\begin{prf} (i) Let $W \subeq \vir$ be an open invariant convex cone. 
If $\fz(\vir) \cap W = \eset$, then 
$\fz(\vir) +  W$ is a proper cone, 
and therefore its image in $\cV(\bS^1)$ is contained in $W_{\cV(\bS^1)}$ or 
$-W_{\cV(\bS^1)}$, 
so that $W \subeq W_{\rm max}$ or $W \subeq - W_{\rm max}$ 
(Theorem~\ref{thm:classcon}). 
To verify (i), we therefore have to show that 
$c \not\in C := W \cap \ft$.  Suppose the converse. 
Then there exists an $\eps > 0$ with 
$x_\pm := c \pm \eps \partial_\theta \in C$. 
From Theorem~\ref{thm:concoad-vir} 
we derive that 
$$ p_\ft(\oline{\conv}(\cO_{x_+})) = x_+ + C_+, $$
and, applying it also to $-x_-$, we find that 
$$ p_\ft(\oline{\conv}(\cO_{x_-})) = x_- - C_+. $$
Since both sets are contained in $\oline C$, we see that 
$\pm C_+ \subeq \lim(\oline C) = \lim(C)$ 
(Lemma~\ref{lem:limcone}), so that $\lim(C) = \ft$, and thus 
$0 \in C = \ft$. As $W$ is open, $0 \in W$ leads to 
$W = \vir$. 

(ii) follows from (i) and Proposition~\ref{prop:8.11}. 

(iii) If $C = W \cap \ft$ for an invariant open convex cone $W$, 
then $C \subeq C_{\rm max}$ implies that $C \cap \z(\vir) = \{0\}$. 
Therefore we have for $x \in C$ the relation 
$$ p_\ft(\oline{\conv}(\cO_x)) = x + C_+ \subeq \oline{C}, $$
and thus $C_+ \subeq \lim(\oline{C}) = \oline{C}$, 
which in turn yields $C_+^0 \subeq C$. 

If, conversely, $C_+^0 \subeq C$ holds for an open convex cone 
$C \subeq C_{\rm max}$, 
then $C_+ \subeq \lim(C) = \oline C$ leads for each $x \in C$ to 
$$ p_\ft(\cO_x) \subeq x + C_+ \subeq C. $$
Therefore 
$$ W_C 
:= \{ x \in \vir \: p_\ft(\cO_x) \subeq C \}  
= \bigcap_{\phi \in G} \phi.p_\ft^{-1}(C) $$ 
is a convex invariant cone containing the subset 
$\Ad(G)C$ which is open by Proposition~\ref{prop:8.11}. 
Hence $W_C^0$ is an open invariant convex cone satisfying \break 
$W_C^0 \cap \ft =~C$. 

(iv) follows immediately from (iii). 

(v) If $W$ is pointed, i.e., it contains no affine lines, 
then the same holds for $C := W \cap \ft$. If, conversely, 
$C$ contains no affine lines, then 
$H(W)$ is a closed ideal of $\vir$ intersecting 
$\ft$ trivially. Hence it is contained in 
$[\ft, \vir]$ and $H(W)_\C$ is adapted to the root 
decomposition with respect to $\ft_\C$. If it contains 
$d_n$, then its $*$-invariance implies that it also contains $d_n^* = d_{-n}$, 
which leads to the contradiction $[d_n, d_{-n}] \in H(W)_\C$. 
This implies that $H(W) = \{0\}$, so that $W$ is pointed. 
\end{prf} 

As a consequence of the preceding theorem, the cones 
$W_{\rm min}$, resp., $W_{\rm max}$ play the role of a minimal, resp., 
maximal open invariant cone in $\vir$. The existence of minimal 
and maximal invariant cones is a well known phenomenon for 
finite dimensional hermitian Lie algebras (cf.\ \cite{Vin80} and 
Proposition~\ref{prop:simpfin}). 

\begin{cor} \mlabel{cor:d-gen} 
The smallest closed convex invariant cone in $\vir$ containing 
$\partial_\theta$ is the closure of $W_{\rm min}$. 
\end{cor}

\begin{prf} If $D \subeq \vir$ is a closed convex invariant cone, 
then $\partial_\theta \in D$ implies that 
$\partial_\theta + C_+ \subeq D$
(Theorem~\ref{thm:concoad-vir}), so that 
$C_+ \subeq D$ (Lemma~\ref{lem:limcone}), 
and thus $W_{\rm min} \subeq D$ by invariance. 
\end{prf}

\begin{rem} In view of the preceding theorem, the open invariant convex 
cones in $\vir$ can be classified as follows. Since the closure of the 
cone $C$ contains 
$C_+ = \R_+ c + \R_+ \partial_\theta$ and is contained in 
$\oline{C_{\rm max}} = \R c + \R_+ \partial_\theta$, 
we have $\oline C = \R_+ c + \R_+ (\partial_\theta - \alpha c)$ 
for some $\alpha > 0$ whenever $C \not= C_{\rm min}, C_{\rm max}$. 
\end{rem}

\subsection{Semi-equicontinuity of coadjoint orbits of $\vir'$} 

In this final section on the Virasoro algebra we apply the 
detailed results on invariant cones to semibounded representations 
and semi-equicontinuous coadjoint orbits. 
In particular, we show that the set $\g'_{\rm seq}$ of semi-equicontinuous 
coadjoint orbits coincides with 
the double cone $W_{\rm min}^\star \cup -W_{\rm min}^\star$. 
This in turn is used to show that the unitary highest weight 
representations of the Virasoro group are precisely the 
irreducible semibounded representations 
and to determine their momentum sets. 

\begin{prop} For $\lambda \in \vir'$ and $d = (0,\partial_\theta)$, 
the following 
are equivalent:
\begin{description}
\item[\rm(i)] $\cO_\lambda$ is semi-equicontinuous. 
\item[\rm(ii)] The convex cone $B(\cO_\lambda)$ contains 
$W_{\rm max}$ or $-W_{\rm max}$. 
\item[\rm(iii)] $\cO_\lambda(d) = \lambda(\cO_d)$ is bounded from 
below or above. 
\item[\rm(iv)] $\lambda \in W_{\rm min}^\star \cup -W_{\rm min}^\star$. 
\item[\rm(v)] $\cO_\lambda(d) \geq 0$ or $\lambda(\cO_d) \leq 0$.
\end{description}
\end{prop} 

\begin{prf} (i) $\Rarrow$ (ii): If $\cO_\lambda$ is semi-equicontinuous, 
then the invariant convex cone $B(\cO_\lambda)$ has interior points 
and contains $\fz(\vir)$. Therefore 
$B(\cO_\lambda)^0/\fz(\vir)$ is an open invariant convex cone 
in $\cV(\bS^1)$, hence contains $W_{\cV(\bS^1)}$ or $-W_{\cV(\bS^1)}$ 
(Theorem~\ref{thm:classcon}). 
This in turn implies that $B(\cO_\lambda)$ contains either 
$W_{\rm max}$ or $-W_{\rm max}$. 
  
(ii) $\Rarrow$ (i) follows from Proposition~\ref{prop:seqcrit} 
because $\vir$ is a Fr\'echet space. 

(ii) $\Rarrow$ (iii) follows from $d \in W_{\rm max}$. 

(iii) $\Rarrow$ (ii) follows from $\R c \subeq B(\cO_\lambda)$ 
and $W_{\rm max} = \Ad(G)(\R c + \R_+^\times d)$ 
(Proposition~\ref{prop:8.11}). 

(i) $\Rarrow$ (iv): If $\lambda(c) = 0$, then $\cO_\lambda$ can be 
identified with a semi-equicontinuous coadjoint orbit of $\cV(\bS^1)$, 
so that Proposition~\ref{prop:semieq-vect} implies 
that $\lambda \in \break W_{\cV(\bS^1)}^\star \cup - W_{\cV(\bS^1)}^\star$ 
in $\cV(\bS^1)'$, 
which means that $\lambda \in W_{\rm max}^\star \cup - W_{\rm max}^\star$ in 
$\vir'$. 

If $\lambda(c) \not=0$, then $\cO_\lambda$ is contained 
in the closed invariant hyperplane $\lambda + c^\bot$, so that 
the construction in Remark~\ref{rem:conered} implies that the cone 
$\cO_\lambda^\star$ has interior points. 
Clearly, this cone is proper, so that Theorem~\ref{thm:classinvcon-vir} 
implies that it either contains $W_{\rm min}$ or $-W_{\rm min}$, 
which in turn leads to $\lambda \in W_{\rm min}^\star \cup - W_{\rm min}^\star$. 

(iv) $\Rarrow$ (i): Since $\pm W_{\rm min}$ are open invariant cones, 
their duals are semi-equicontinuous sets (Example~\ref{ex:2.3}(b)).

(iv) $\Leftrightarrow$ (v): Since the closed convex invariant cone generated 
by $d$ is $W_{\rm min}$ (Corollary~\ref{cor:d-gen}), 
$\cO_\lambda(d) \geq 0$ is equivalent to 
$\lambda \in W_{\rm min}^\star.$
\end{prf}

For any $\lambda \in \ft^* \cong [\ft,\g]^\bot \subeq \g'$, the fact that 
$\cO_\lambda$ is constant on the central element 
$c= (1,0)$ implies that
$p_{\ft^*}(\cO_\lambda) \subeq \ft^*$ is a connected 
subset of the affine line 
$$\{ \mu \in \ft^* \: \mu(c) = \lambda(c)\} 
= \lambda + (\ft^* \cap c^\bot). $$
In particular, this set is convex. 

\begin{prop} \mlabel{prop:semeq-vir}
If $\lambda \in \ft^*$, then 
\begin{description}
\item[\rm(a)] $\cO_\lambda(d)$ is bounded from below 
if and only if $\lambda(d)\geq 0$ and $\lambda(c) \geq 0$. 
If this is the case and $\lambda \not=0$, then  
$$ B(\cO_\lambda)^0 = W_{\rm max} \quad \mbox{ and } \quad 
W_{\rm min} \subeq \cO_\lambda^\star. $$
\item[\rm(b)] $p_{\ft^*}(\cO_\lambda)$ is contained in an affine 
half-line if and only if $\cO_\lambda$ is semi-equi\-con\-ti\-nuous 
if and only 
if $\lambda(c)\lambda(d) \geq 0$.  
\end{description}
\end{prop}

\begin{prf} (a)  We recall from Theorem~\ref{thm:concoad-vir} that 
$p_{\ft}(\oline{\conv}(\cO_d)) = d + C_+.$
This implies that 
$\lambda(\cO_d)$ is bounded from below if and only if 
$\lambda \in (C_+)^\star$, i.e., 
$\lambda(c), \lambda(d) \geq 0$. 

Suppose that these conditions are satisfied and that $\lambda \not=0$. 
Then $B(\cO_\lambda)^0$ is a proper open invariant cone, hence 
determined by its intersection with $\ft$ (Theorem~\ref{thm:classinvcon-vir}).
 As this intersection 
contains $d$ and is invariant under translation with $\R c$, 
Theorem~\ref{thm:classinvcon-vir} implies that it coincides with 
$C_{\rm max}$. This proves that $W_{\rm max} = B(\cO_\lambda)^0$. 
We have already seen above that $\lambda \in C_+^\star$ and 
since $p_{\ft^*}(\cO_\lambda)$ is a half-line constant 
on $c$ and bounded below on $d$, 
it follows that $C_+ \subeq \cO_\lambda^\star$, which 
leads to $W_{\rm min} \subeq \cO_\lambda^\star$. 

(b) We use Proposition~\ref{prop:orb-pro} to obtain with the notation 
of Remark~\ref{rem:vir-rel} 
$$p_{\ft^*}(\cO_\lambda) \supeq 
\lambda + \R^+ \lambda([d_n^*, d_n]) \alpha_n. $$ 
We also know that 
$i [d_n^*, d_n] = 2n (\pi(n^2 - 1), \partial_\theta),$
so that, for each $n \in \N$, 
$$p_{\ft^*}(\cO_\lambda) \supeq 
\lambda -  \R^+ \lambda(\pi(n^2-1), \partial_\theta) i \alpha_n 
= \lambda -  \R^+ \lambda(\pi(n^2-1), \partial_\theta) i \alpha_1. $$ 
If $p_{\ft^*}(\cO_\lambda)$ is contained in a half-line, the signs 
of the numbers 
$$ \lambda(\pi(n^2-1), \partial_\theta), \quad n \in \N, $$ 
have to coincide, which is equivalent to 
$\lambda(d) \lambda(c) \geq 0.$
If, conversely, this condition is satisfied, 
then (a) implies that $\cO_\lambda(d)$ is semibounded, 
so that $p_{\ft^*}(\cO_\lambda)$ is contained 
in an affine half-line. 
\end{prf}

\begin{defn} We write $\Vir$ for the (up to isomorphism unique) 
simply connected Lie group with Lie algeba $\vir$. 
\glossary{name={$\Vir$},description={Virasoro group, simply connected with 
$\L(\Vir) = \vir$}}

A unitary representation $(\pi, \cH)$ of $\Vir$ 
is called a {\it highest weight representation} 
if there exists a smooth cyclic vector $0 \not=v \in \cH^\infty$ 
which is a $\ft$-eigenvector annihilated by each $d_n$, $n > 0$. 
Then the corresponding eigenfunctional $\lambda \in i\ft^*$ is called 
the {\it highest weight} and $v$ a highest weight vector. 
\end{defn}

\begin{thm} \mlabel{thm:hiwei-vir} 
If $(\pi_\lambda, \cH_\lambda)$ is a unitary 
highest weight representation of the simply connected
 Lie group $\Vir$ with Lie algebra $\vir$ of highest weight 
$\lambda \in i \ft^*$, then 
\begin{equation}
  \label{eq:hiweisign}
\lambda(d_0) \geq 0 \quad \mbox{ and } \quad 
\lambda(\hat c) \geq 0, 
\end{equation}
i.e., $i \lambda(\partial_\theta) \geq 0$ 
and $i\lambda(c) \geq 0.$ 
The representation $(\pi_\lambda, \cH_\lambda)$ 
is semibounded and its momentum set is given by 
$$ I_{\pi_\lambda} = \oline{\conv}(\cO_{-i\lambda}). $$
For $\lambda\not=0$ we have 
$$ -W_{\rm max} = W_{\pi_\lambda} \quad \mbox{ and } \quad 
-W_{\rm min} \subeq C_{\pi_\lambda}. $$
\end{thm} 

\begin{prf} First we use \cite[Prop.~3.5]{KR87} to see that the 
unitarity of the irreducible highest weight module $L(\lambda)$ of $\vir_\C$ 
with highest weight 
$\lambda \in i\ft^*$ implies \eqref{eq:hiweisign}. 
Actually, this follows from the simple observation that 
if $v_\lambda$ is a highest weight vector of unit length, then 
\begin{align*}
0 &\leq 
\la d_{-n}v_\lambda, d_{-n}v_\lambda \ra 
= \la d_{-n}^*d_{-n}v_\lambda, v_\lambda \ra 
= \la [d_{-n}^*, d_{-n}]v_\lambda, v_\lambda \ra \\
&= \lambda([d_n, d_{-n}]) 
= \lambda\Big(2n d_0 + \frac{n^3 - n}{12}\hat c\Big) \geq 0 
\end{align*}
holds for each $n \in \N$. 

The existence of a corresponding continuous unitary highest weight 
representation $(\pi_\lambda, \cH_\lambda)$ of the simply connected 
Lie group $\Vir$ with Lie algebra $\vir$ has been shown by Goodman 
and Wallach \cite{GW85}. A more general method of integration 
which applies in particular to highest weight modules of $\Vir$ 
has been developed by Toledano Laredo (\cite[Thm.~6.1.1]{TL99b}). 
It is based on techniques related to regular Lie groups, 
and \cite[Cor.~4.2.2]{TL99b} implies in particular that the 
highest weight vector $v_\lambda$ is smooth. 
This vector is an eigenvector for the closed subalgebra 
of $\vir_\C$ generated by $\ft_\C$ and the 
$d_n$, $n > 0$. 

Next we recall from Theorem~\ref{thm:8.6} the complex manifold 
$$ M = \Diff(\bS^1)_+/\exp(\R \partial_\theta)  
\cong \Vir/T, $$
where $T := \exp\ft \subeq \Vir$ is the subgroup corresponding to 
$\ft$. Since the tangent space in the base point 
$x_0 = \1 T$ can be identified with 
$\vir_\C/\fp$, $\fp := \oline{\ft_\C + \sum_{n < 0} \C d_n}$, 
Theorem~\ref{thm:5.7} implies that the 
map 
$$\eta \: M \to \bP(\cH_\lambda'),\quad gT \mapsto 
[\pi^*_\lambda(g)\alpha_{v_\lambda}]  $$
is holomorphic. As $v_\lambda$ is cyclic, we thus obtain a 
realization of the unitary representation $(\pi_\lambda, \cH_\lambda)$ 
in the space of holomorphic sections of the holomorphic line bundle 
$\eta^*\bL_{\cH_\lambda'}$ over 
$M$ (cf.~Theorem~\ref{thm:5.7})\begin{footnote}
{An infinitesimal version of this construction 
can already be found in \cite{KY88} and \cite{Ki98} contains 
various formal aspects of the realization of the highest weight 
representations in spaces of holomorphic functions, resp., 
sections of holomorphic line bundles on $M$.}  
\end{footnote} 

From the highest weight structure of $\cH_\lambda$ it follows that the 
set of $\ft_\C$-weights on $\cH_\lambda$ is given by 
$\lambda - \N_0 \alpha_1$, so that $i \dd\pi_\lambda(\partial_\theta) 
= \dd\pi_\lambda(d_0)$ is bounded from below. Therefore the 
cone $C_{\rm max}$ from the Classification Theorem~\ref{thm:classinvcon-vir} 
satisfies $- C_{\rm max} \subeq W_{\pi_\lambda}$, 
which immediately leads to $- W_{\rm max} \subeq W_{\pi_\lambda}$. 
Since  $\partial_\theta\not\in W_{\pi_\lambda}$, 
this cone is proper, and the Classification Theorem thus implies 
the equality $W_{\pi_\lambda} = -W_{\rm max}$. 
With this information we now apply Theorem~\ref{thm:2.7}  
to determine the precise momentum set. 

From Theorem~\ref{thm:8.6} it now follows that,  
for each $x \in W_{\rm max} = - W_{\pi_\lambda} = W_{\pi_\lambda^*}$, 
the smooth action of the corresponding one-parameter group 
on $M$ extends to a holomorphic action of $\C_+$. 
Therefore Theorem~\ref{thm:5.7}(c) implies that 
$I_{\pi_\lambda}$ is the closed convex hull of the coadjoint 
orbit $\Phi_{\pi_\lambda}(G[v_\lambda])$. In view of 
$\Phi_{\pi_\lambda}([v_\lambda]) = -i\lambda \in \ft^* \subeq \vir'$, 
this proves that $I_{\pi_\lambda} = \oline{\conv}(\cO_{-i\lambda})$. 
Now 
$$C_{\pi_\lambda} = (I_{\pi_\lambda}^\star)^0 =  
(\cO_{-i\lambda}^\star)^0 \supeq -W_{\rm min} $$ 
follows directly from Proposition~\ref{prop:semeq-vir}, 
and since $W_{\rm min}$ has interior points, 
$(\pi_\lambda, \cH_\lambda)$ is semibounded. 
\end{prf}

For more details on the classification of unitary highest 
weight modules of $\vir$, we refer to \cite{KR87}.  
We conclude this section with the following converse of 
Theorem~\ref{thm:hiwei-vir}: 

\begin{thm} \mlabel{thm:8.22}  
Every irreducible semibounded representation 
$(\pi, \cH)$ of $\Vir$ is either a highest weight representation 
or the dual of a highest weight representation. 
\end{thm} 

\begin{prf} We assume that the representation 
$\pi$ is non-trivial, so that $I_\pi \not= \{0\}$. 
Let $p_{\ft^*} \: \vir' \to \ft^*$ be the restriction map. 
If $p_{\ft^*}(I_\pi)= \{0\}$, then 
$\ft \subeq I_\pi^\bot = \ker \dd\pi$ and thus 
$\vir= \ft + [\ft,\vir] \subeq \ker \dd\pi$, contradicting the 
non-triviality of the representation. 
Next we use Proposition~\ref{prop:project} 
to conclude that $\{0\} \not= p_{\ft^*}(I_\pi)\subeq I_\pi$. 
Then each non-zero $\lambda \in \ft^* \cap I_\pi$ has a 
semi-equicontinuous orbit, so that 
$\cO_\lambda(d)$ is bounded from below or above, and 
in this case $B(\cO_\lambda)^0 = W_{\rm max}$ or $-W_{\rm max}$ 
(Proposition~\ref{prop:semeq-vir}). This implies in particular 
that $W_\pi \not=\vir$, i.e., $\pi$ is not bounded. 

As the open invariant cone $W_\pi$ is proper, 
Theorem~\ref{thm:classinvcon-vir}(i) 
implies that $W_\pi$ is either 
contained in $W_{\rm max}$ or $-W_{\rm max}$. We assume the latter 
and claim that $\pi$ is a highest weight representation. 
In the other case, 
$W_{\pi^*} = - W_\pi \subeq - W_{\rm max}$, so that $\pi$ is the 
dual of a highest weight representation. 

First we note that $W_\pi \subeq - W_{\rm max}$ implies 
$-W_{\rm min} \subeq W_\pi$. 
As $c \in H(W_\pi)$ follows from $\dd\pi(c) \in i \R\1$ (Schur's Lemma), 
we obtain the relation $\R c - W_{\min} \subeq W_\pi$, and therefore 
$W_{\rm max} = W_{\rm min} + \R c$ leads to 
$W_\pi = - W_{\rm max}\ni - d$. 
Hence the spectrum of $i\dd\pi(d)= \dd\pi(d_0)$ is bounded from below. 
In view of $\exp(2\pi d) \in Z(\Vir)$, 
Proposition~\ref{prop:maxsmoothvect} implies the 
existence of a smooth unit vector $v \in\cH^\infty$ which is an 
eigenvector for the minimal  eigenvalue $h$ of $i\dd\pi(d)$. 
Then $\ft = \R c + \R d$, $v$ is a $\ft$-eigenvector and 
$[d_0, d_n] = - n d_n$ implies 
that $d_n.v = 0$ for each $n > 0$. 
Therefore $(\pi, \cH)$ is a highest weight representation. 
\end{prf}

\section{Symplectic group and metaplectic representation} \mlabel{sec:9}

In this section we study the metaplectic representation 
$(\pi_s, S(\cH))$ of the 
central extension $\hat\Sp_{\rm res}(\cH)$ of $\Sp_{\rm res}(\cH)$ 
on the symmetric Fock space 
$S(\cH)$. This representation arises from self-intertwining 
operators of the Fock representation of the Heisenberg group 
$\Heis(\cH)$. We show that it is semibounded and 
determine the corresponding cone $W_{\pi_s}$. 
For the larger central extension 
$\hat\HSp_{\rm res}(\cH)$ of $\cH \rtimes \Sp_{\rm res}(\cH)$ 
the  representation on $S(\cH)$ is irreducible 
and semibounded and we show that its momentum set 
is the weak-$*$-closed convex hull of a single coadjoint orbit. 

\subsection{The metaplectic representation} \mlabel{sec:9.1}

On the dense subspace $S(\cH)_0 = \sum_{n = 0}^\infty S^n(\cH)$ of the 
symmetric Fock space $S(\cH)$ (cf.\ Appendix~\ref{app:c}) 
we have for each $f \in \cH$ the {\it creation operator} 
$$ a^*(f)(f_1 \vee \cdots \vee f_n) := f \vee f_1 \vee \cdots \vee f_n. $$
This operator has an adjoint $a(f)$ on $S(\cH)_0$, given by 
$$ a(f)\Omega = 0, \quad 
a(f)(f_1 \vee \cdots \vee f_n) 
= \sum_{j = 0}^n \la f_j, f \ra f_1 \vee \cdots \vee \hat{f_j} \vee 
\cdots \vee f_n, $$
where $\hat{f_j}$ means omitting the factor $f_j$. Note that 
$a(f)$ defines a derivation on the algebra $S(\cH)_0$. 
One easily verifies 
that these operators satisfy the 
canonical commutation relations (CCR): 
\begin{equation}
  \label{eq:ccr}
[a(f), a(g)] = 0, \quad [a(f), a^*(g)] = \la g, f\ra \1. 
\end{equation}
For each $f \in \cH$, the operator $a(f) + a^*(f)$ on $S(\cH)_0$ 
is essentially self-adjoint (\cite[p.~70]{Ot95}), so that 
$$ W(f) := e^{\frac{i}{\sqrt{2}}\oline{a(f) + a^*(f)}} \in \U(\cH) $$
is a unitary operator. These operators satisfy the {\it Weyl 
relations} 
$$ W(f) W(f') = W(f + f') e^{\frac{i}{2}\Im \la f, f'\ra},\quad 
f,f' \in \cH. $$
For the Heisenberg group 
$$ \Heis(\cH) := \R \times \cH, $$
with the multiplication 
$$ (t,v) (t',v') := (t + t' + \shalf\omega(v,v'), v + v'), \quad 
\omega(v,v') := \Im \la v,v'\ra $$
we thus obtain by $W(t, f) := e^{it} W(f)$
a unitary representation on $S(\cH)$, called the {\it 
Fock representation}. 
It is a continuous irreducible representation 
(\cite[Cor.~3.11]{Ot95}). That it actually is smooth 
follows from the smoothness of 
\begin{equation}
  \label{eq:matcoeff}
\la W(t,f)\Omega, \Omega \ra 
= e^{it - \frac{1}{4}\|f\|^2} 
\end{equation}
and Theorem~\ref{thm:unitautsmooth}. 

\begin{defn}  The symplectic group $\Sp(\cH)$ acts via 
$g.(t,v) := (t,gv)$ by automorphism on the Heisenberg group, 
and since the unitary representation 
$(W,S(\cH))$ is irreducible, there exists for each 
$g \in \Sp(\cH)$, up to multiplication with $\T$, at most one unitary 
operator $\pi_s(g) \in \U(S(\cH))$ with 
\begin{equation}
  \label{eq:autrel}
\pi_s(g) W(t,f) \pi_s(g)^* = W(t,gf) \quad \mbox{ for } \quad 
t\in \R, f \in \cH. 
\end{equation}
According to \cite{Sh62}, such an operator exists 
if and only if $g \in \Sp_{\rm res}(\cH)$ (cf.\ Example~\ref{ex:liegrp}(d)), 
which immediately 
leads to a projective unitary representation 
$$ \oline\pi_s \: \Sp_{\rm res}(\cH) \to \PU(S(\cH)),  $$
determined by \eqref{eq:autrel} for any lift 
$\pi_s(g)\in \U(S(\cH))$ of $\oline\pi_s(g)$. 
Writing $\oline u$ for the image of $u \in \U(\cH)$ in the projective 
unitary group $\PU(\cH)$, the corresponding pull-back 
$$ \hat\Sp_{\rm res}(\cH)  := \oline\pi_s^*\U(S(\cH)) 
= \{ (g,u) \in \Sp_{\rm res}(\cH) \times \U(S(\cH)) \: 
\oline u = \oline\pi_s(g)\}$$ 
is called the {\it metaplectic group}. 
\glossary{name={$\hat\Sp_{\rm res}(\cH)$},description={metaplectic group}}
\glossary{name={$(\pi_s, S(\cH))$},description={metaplectic representation of metapl. group}}
It is a central 
extension of $\Sp_{\rm res}(\cH)$ by $\T$. We shall 
see below that this group is a Lie group and that 
its canonical representation $\pi_s(g,u) = u$ 
on $S(\cH)$, the {\it metaplectic 
representation}, is smooth and semibounded. Our strategy 
is to use Theorem~\ref{thm:unirep-cenext}, 
which requires a suitable lift of $\oline\pi_s$. 
\end{defn}

\begin{rem} \mlabel{rem:9.2} 
From the relation $\dd W(f) = \frac{i}{\sqrt 2}(a(f) + a^*(f))$, 
we recover the antilinear and linear part of $\dd W$ by 
$$ a(f) = \frac{1}{i\sqrt 2} (\dd W(f) + i \dd W(If)), \quad 
a^*(f) =  \frac{1}{i\sqrt 2}(\dd W(f) - i \dd W(If)), $$ 
so that by \eqref{eq:autrel} 
$$ \pi_s(g) a(f) \pi_s(g)^{-1} 
= \frac{1}{i \sqrt{2}}\big(\dd W(gf) + i \dd W(gIf)\big) 
= a(g_1 f) + a^*(g_2 f) =: a_g(f), $$ 
where $g = g_1 + g_2$ is the decomposition into linear an antilinear 
part.   
\end{rem}

\begin{thm}\mlabel{thm:metasmooth} 
The topological group $\hat\Sp_{\rm res}(\cH)$ is a Lie group and 
the metaplectic representation is smooth. 
A Lie algebra cocycle $\eta$ defining $\hat\sp_{\rm res}(\cH)$ 
as an extension of $\sp_{\rm res}(\cH)$ by $\R$ is given by 
$$ \eta(x,y) = \frac{1}{2i} \tr([x_2, y_2]), $$
where $x_2$ denotes the antilinear component of $x$.
\end{thm}

\begin{prf} 
{\bf Step 1:} We consider the unbounded operators $a(f)$ and $a^*(f)$ on 
$S(\cH)_0 \subeq S(\cH)$ and observe that $\C\Omega$ is the common kernel
of the annihilation operators $a(f)$. For $g\in \Sp_{\rm res}(\cH)$, 
we observe that $\pi_s(g)\Omega$ lies in the common kernel of the 
operators $a_g(f) = \pi_s(g) a(f) \pi_s(g)^{-1}$ 
(cf.\ Remark~\ref{rem:9.2}). 
If $g_1$ is invertible, which is the case in some open 
$\1$-neighborhood in $\Sp_{\rm res}(\cH)$ (actually on the whole group), 
we consider the antilinear operator $T(g) := g_2 g_1^{-1}$ 
for which $a_g(g_1^{-1} f) = a(f) + a^*(T(g)f)$. 
Therefore $F = \pi_s(g)\Omega$ is a solution of the 
following system of equations: 
\begin{equation}
  \label{eq:absdiffeq2}
a(f)F = - a^*(T(g)f)F \quad \mbox{ for } \quad f \in \cH. 
\end{equation}

{\bf Step 2:} For each $n \in \N$, the subset $a^*(\cH)S^n(\cH)$ is total
in $S^{n+1}(\cH)$, which implies that 
\begin{equation}
  \label{eq:homker}
\{ T \in S^{n+1}(\cH) \: (\forall f \in \cH)\, a(f)T = 0\} 
= \{0\}. 
\end{equation}
If an element $F= \sum_{n = 0}^\infty F_n \in S(\cH)$ 
with $F_n \in S^n(\cH)$ satisfies 
\eqref{eq:absdiffeq2}, then 
\begin{equation}
  \label{eq:absdiffeq3}
a(f)F_1 = 0, \quad 
a(f)F_{n+1} = - a^*(T(g)f)F_{n-1} \quad \mbox{ for } \quad f \in \cH, 
n \in \N.  
\end{equation}
This implies $F_1 = 0$, and inductively we obtain with 
\eqref{eq:homker} $F_{2k+1} = 0$ for $k \in \N_0$. 
We also derive from \eqref{eq:homker} that $F$ is completely determined 
by $F_0$, hence that the solution space of \eqref{eq:absdiffeq2} 
is at most one-dimensional. 

If $F$ is a solution, we may w.l.o.g.\ assume that $F_0 = \Omega$. 
Then $F_2$ satisfies 
$$ a(f) F_2 = - a^*(T(g)f) \Omega = - T(g)f, \quad f \in \cH, $$
i.e., $F_2 = - \hat{T(g)}$ (Lemma~\ref{lem:hs2}). 
This observation implies Shale's 
result that only for $g \in \Sp_{\rm res}(\cH)$, i.e., 
$\|g_2 \|_2 < \infty$, equation 
\eqref{eq:autrel} has a solution $\pi_s(g)\in \U(S(\cH))$. 

{\bf Step 3:} Combining Lemma~\ref{lem:hs2}(i) with 
Remark~\ref{rem:fock-alg}, we conclude that the exponential series 
$e^{-\hat{T(g)}}$ converges  in $S(\cH)$ 
for $\|T(g)\|_2 < 1$, which holds on an open $\1$-neighborhood in 
$\Sp_{\rm res}(\cH)$. 
Since the operators $a(f)$ act as derivations on $S(\cH)_0$, it follows that 
$$ a(f) e^{-\hat{T(g)}} 
=  -  a(f) \hat{T(g)}  \vee e^{-\hat{T(g)}} 
=  -  T(g)f  \vee e^{- \hat{T(g)}} 
=  - a^*(T(g)f) e^{-\hat{T(g)}}, $$ 
so that $e^{-\hat{T(g)}}$ satisfies  \eqref{eq:absdiffeq2}. 
We conclude that 
\begin{equation}
  \label{eq:omega_g}
\pi_s(g) \Omega = c(g) e^{-\hat{T(g)}}, \quad c(g) \in \C^\times. 
\end{equation}
Choosing the operators $\pi_s(g)$, $g \in \Sp_{\rm res}(\cH)$,  
in such a way that 
$$c(g) = \la \Omega_g, \Omega\ra = \la \pi_s(g)\Omega, \Omega \ra > 0,$$  
it follows that 
$c(g) = \|e^{-\hat{T(g)}}\|^{-1}$ 
(cf.\ \cite[p.~97]{Ot95}). 

{\bf Step 4:} Since the map 
$g \mapsto T(g) = g_2 g_1^{-1}$ 
is smooth in an identity neighborhood and the map 
$$ \Big\{ A \in \fp_2 \: \|A\|_2 < 1\Big\} \to S(\cH), 
\quad A \mapsto e^{\hat A} $$
is analytic (cf.~Remark~\ref{rem:fock-alg} and  
the proof of Lemma~\ref{lem:sp.1}(v) for $\fp_2$), 
hence in particular smooth, it follows that 
$e^{-\hat{T(g)}}$ and hence also $c(g)$ are smooth in an identity neighborhood 
of $\Sp_{\rm res}(\cH)$. 

From $\la \pi_s(g)^*\Omega, \Omega\ra = \la \Omega, \pi_s(g)\Omega\ra > 0$ 
we further obtain $\pi_s(g)^* = \pi_s(g^{-1})$. 
This implies that the map
$$ (g,h) \mapsto \la \pi_s(g)\pi_s(h)\Omega,\Omega\ra 
= \la \pi_s(h)\Omega,\pi_s(g)^*\Omega\ra 
= \la \pi_s(h)\Omega,\pi_s(g^{-1})\Omega\ra $$
is smooth in an identity neighborhood. Now 
Theorem~\ref{thm:unirep-cenext} implies that 
$\hat\Sp_{\rm res}(\cH)$ is a Lie group and that 
$\Omega$ is a smooth vector for the corresponding unitary representation, 
also denoted~$\pi_s$. 
Since $\Sp_{\rm res}(\cH)$ acts smoothly on $\Heis(\cH)$, the 
space $S(\cH)^\infty$ of smooth vectors is invariant under 
$\Heis(\cH)$, so that the irreducibility of the Fock 
representation of $\Heis(\cH)$ on $S(\cH)$ (\cite[Cor.~3.11]{Ot95}) 
implies the smoothness of~$\pi_s$. 

{\bf Step 5:} With Remark~\ref{rem:a.11} we can now calculate a 
suitable cocycle $\eta$ with $F(g) := \pi_s(g)\Omega = 
c(g) e^{-\hat{T(g)}}$ by 
$$  \eta(x,y) 
= 2\Im \la \dd F(\1)x, \dd F(\1)y\ra  + i \la \dd F(\1)[x,y],\Omega \ra.$$
As $F(g)$ only depends on $T(g) = g_2 g_1^{-1}$, we have
$F(gu) = F(g)$ for $u \in \U(\cH)$, which leads to 
$\dd F(\1)x = 0$ for $x \in \fu(\cH)$. Hence 
$\eta(x,\cdot) = 0$ for $x \in \fu(\cH)$. 
For $x \in \fp_2$ we find with 
$T(\exp x) = \cosh(x)\sinh(x)^{-1} = \tanh(x)$ 
the relation $\dd T(\1)x = x$, and hence 
$$ \dd F(\1)x = \dd c(\1)(x) \Omega -\hat x. $$
Since $c$ is real-valued, this leads with Lemma~\ref{lem:hs2}(ii) 
for $x, y \in \fp_2$ to 
\begin{align*}
 \eta(x,y) 
&= 2\Im \la \dd c(\1)(x) \Omega - \hat x, 
\dd c(\1)(y) \Omega - \hat y \ra 
= 2 \Im \la \hat x, \hat y \ra\\
&= \Im \tr(xy) = \shalf\Im \tr([x,y]) = \frac{1}{2i} \tr([x,y]) 
\end{align*}
because the trace of the symmetric operator $xy+yx$ is real 
and $\tr([x,y]) \in i \R$. 
\end{prf}

The content of the preceding theorem is essentially known 
(cf.\ \cite{Ot95}, \cite{Ve77}, \cite{Sh62}), 
although all references known to the author only discuss 
the metaplectic representation 
as a representation of a topological group 
and not as a Lie group. In \cite[Thm.~3.19]{Ot95} 
and \cite{Lm94} one finds quite explicit formulas for group 
cocycles describing $\hat\Sp_{\rm res}(\cH)$ as a central extension. 

\begin{rem} \mlabel{rem:9.4} 
The metaplectic representation of 
$G := \hat\Sp_{\rm res}(\cH)$ is not irreducible. 
Since the representation of $\U(\cH)$ on each  $S^n(\cH)$ 
is irreducible (Example~\ref{ex:7.3}), every 
$G$-invariant subspace is the direct sum  of some $S^n(\cH)$, 
$n \in \N_0$. 

The preceding proof immediately shows that 
$$\pi(G)\Omega \subeq 
S^{\rm even}(\cH) := \hat\oplus_{n \in \N_0} S^{2n}(\cH)$$ 
and that all projections onto the subspaces $S^{2n}(\cH)$ 
are non-zero. 

Using the fact that the operators in 
$\dd\pi_s(\g)$ contain all multiplication operators 
with elements $\hat A$, $A \in \fp_2$ (cf.\ Lemma~\ref{lem:hs2}), 
and their adjoints (\cite{Ot95}), 
it easily follows that the representations of $G$ on the two 
subspaces 
$S^{\rm even}(\cH)$ and $S^{\rm odd}(\cH)$ are irreducible. 
\end{rem}

\subsection{The metaplectic group} \mlabel{sec:9.2}

In this subsection we describe a convenient description 
of the metaplectic group $\hat\Sp_{\rm res}(\cH)$ as a Banach--Lie group. 
For further details we refer to \cite{Ne02a}.
First we recall from \cite[Def.~III.3]{Ne02a}  
that 
\begin{align*}
\Sp_{1,2}(\cH) 
&:= \Big\{ g = \pmat{ a & b \\ c & d} \in \Sp(\cH) \subeq \U(\cH,\cH) \: \\
&\qquad\qquad \qquad \|b\|_2, \|c\|_2 < \infty, 
\|a - \1\|_1, \|d - \1\|_1 < \infty\Big\} 
\end{align*}
carries the structure of a Banach--Lie group 
with polar decomposition 
$$ \U_1(\cH) \times \fp_2 \to \Sp_{1,2}(\cH), \quad 
\fp_2 = \Big\{ \pmat{ 0 & b \\ b^* & 0} \: b = b^\top \in B_2(\cH)\Big\}$$ 
(cf.\ \cite[Def.~IV.7, Lemma~IV.13]{Ne02a}).   
The full unitary group $\U(\cH)$ acts smoothly by conjugation 
on $\Sp_{1,2}(\cH)$, so that we can form the semidirect product
\break $\Sp_{1,2}(\cH) \rtimes \U(\cH)$, and the multiplication map 
$$ \mu_1 \: \Sp_{1,2}(\cH) \rtimes \U(\cH) \to \Sp_{\rm res}(\cH), 
\quad (g,u) \mapsto gu $$
is a quotient morphism of Banach--Lie groups with kernel 
$$ N := \{ (g,g^{-1}) \: g \in \U_1(\cH) \}, $$
so that 
$$ \Sp_{\rm res}(\cH) \cong (\Sp_{1,2}(\cH) \rtimes \U(\cH))/N $$ 
(cf.\ \cite[Def.~IV.7]{Ne02a}). 

From the polar decomposition of 
$\Sp_{1,2}(\cH)$ we derive that 
$\pi_1(\Sp_{1,2}(\cH)) \cong \pi_1(\U_1(\cH)) \cong \Z$, hence the 
existence of a unique $2$-fold covering group 
\break $q \: \Mp_{1,2}(\cH)\to \Sp_{1,2}(\cH)$. On the inverse 
image $\hat\U_1(\cH)$ of $\U_1(\cH)$ in $\Sp_{1,2}(\cH)$ we then 
have a unique character 
\begin{equation}
  \label{eq:sqdet}
\sqrt{\det} \: \hat\U_1(\cH) \to \T \quad \mbox{ with} \quad 
\L(\sqrt{\det}) = \frac{1}{2i} \tr. 
\end{equation}
Next we observe that the smooth action of $\U(\cH)$ on 
$\Sp_{1,2}(\cH)$  
lifts to a smooth action on $\Mp_{1,2}(\cH)$. We also note that, 
for $\SU_1(\cH) := \ker(\det)$, we have 
$\U_1(\cH) \cong \SU_1(\cH) \rtimes \T$, where the determinant 
is the projection onto the second factor. Accordingly, 
$\hat\U_1(\cH) \cong \SU_1(\cH) \rtimes \T$ with 
$q(g,t) = (g,t^2)$ and $\sqrt{\det}(g,t) = t$. 

Writing $x = x_1 + x_2$ for the decomposition of $x \in \sp(\cH)$  
into linear and antilinear component, 
we recall from Theorem~\ref{thm:metasmooth} that 
$$ \hat\sp_{\rm res}(\cH) \cong \R \oplus_\eta \sp_{\rm res}(\cH) 
\quad \mbox{ with } \quad 
\eta(x,y) = \frac{1}{2i} \tr([x_2, y_2]). $$
This implies that 
$$ \sigma \: \sp_{1,2}(\cH) \to \hat\sp_{\rm res}(\cH) 
\cong \R \oplus_\eta \sp_{\rm res}(\cH), \quad 
\sigma(x) := \Big(\frac{1}{2i} \tr(x_1), x\Big) $$
is a homomorphism of Banach--Lie algebras. Here we use that 
$$ [x,y]_1 = [x_1, y_1] + [x_2, y_2] 
\quad \mbox{ and } \quad \tr([x_1, y_1])=0. $$
On the subgroup $\hat\U_1(\cH) \subeq \Mp_{1,2}(\cH)$, 
$\sigma$ integrates to the group homomorphism 
$$\sigma_G \: \hat\U_1(\cH) \to \hat\Sp_{\rm res}(\cH) 
\subeq \Sp_{\rm res}(\cH) \times \U(S(\cH)), \quad 
\sigma(g) := (q(g), \sqrt{\det}(g)\pi_s(g)), $$
so that the 
polar decomposition of $\Sp_{1,2}(\cH)$ implies that 
$\sigma$ integrates to a morphism of Banach--Lie groups 
$\sigma_G \: \Mp_{1,2}(\cH) \to \hat\Sp_{\rm res}(\cH).$
Combining this map with the canonical inclusion 
$\U(\cH) \into \hat\Sp_{\rm res}(\cH)$, the equivariance of 
$\sigma_G$ under conjugation with unitary operators implies the 
existence of a homomorphism 
$$ \mu \: \Mp_{1,2}(\cH) \rtimes \U(\cH) \to \hat\Sp_{\rm res}(\cH), 
\quad (g,u) \mapsto \sigma_G(g)u.$$ 

\begin{prop} \mlabel{prop:olinemu} The homomorphism $\mu$ 
factors through an isomorphism 
$$ \oline\mu \: (\Mp_{1,2}(\cH) \rtimes \U(\cH))/\ker \mu 
\to \hat\Sp_{\rm res}(\cH)$$ 
of connected Banach--Lie groups with 
$\ker \mu \cong \SU_1(\cH).$ 
\end{prop}

\begin{prf} If $\mu(g,u) = \1$, then $\sigma_G(g) = u^{-1}$ implies that 
$g \in \hat\U_1(\cH)$ with $q(g) = u^{-1}$. 
Now 
$$ \Omega = \pi_s(\mu(g,u)) \Omega = \sqrt{\det}(g)\Omega $$
implies $\sqrt{\det}(g) = 1$. This shows that 
$\ker \mu = \{(g,q(g)^{-1}) \: g \in \SU_1(\cH)\},$ 
and the assertion follows. 
\end{prf}

The inverse image of 
the center $Z \cong \T$ of $\hat\Sp_{\rm res}(\cH)$ is 
$$\hat N := \{ (g,q(g)^{-1}) \: g \in \hat\U_1(\cH)\}.$$ 
The inverse image $\hat\U(\cH)$ of the subgroup $\U(\cH)$ in 
$\hat\Sp_{\rm res}(\cH)$ is of the form 
\begin{equation}
  \label{eq:hatuh}
\hat\U(\cH) \cong 
(\hat \U_1(\cH) \rtimes \U(\cH))/\ker \mu \cong Z \times \U(\cH). 
\end{equation}
In particular, it splits as a central extension of Lie groups.

\subsection{Semiboundedness of the metaplectic representation} 
\mlabel{sec:9.3}

In this subsection we show that the metaplectic representation 
$(\pi_s, S(\cH))$ of $\hat\Sp_{\rm res}(\cH)$ is semibounded and 
determine the open invariant cone $W_{\pi_s}$ as the inverse 
image of the canonical cone $W_{\sp_{\rm res}(\cH)}$ in 
$\sp_{\rm res}(\cH)$. 

We start with a closer look at adjoint orbits of the 
symplectic Lie algebra. 

\begin{lem} \mlabel{lem:convall} If $(V,\omega)$ is a finite dimensional 
symplectic space and $X \in \sp(V,\omega)$ is such that 
the Hamiltonian function $H_X(v) = \shalf\omega(Xv,v)$ is indefinite, then 
$\conv(\cO_X) = \sp(V,\omega)$. 
\end{lem} 

\begin{prf} Let $\g := \sp(V,\omega)$. Then Proposition~\ref{prop:simpfin}, 
combined with the uniqueness of the open invariant cone in 
$\sp(V,\omega)$ up to sign 
(cf.\ \cite{Vin80}), 
implies that the set of semi-equicontinuous coadjoint 
orbits coincides with the double cone 
$W_{\sp(V,\omega)}^\star \cup - W_{\sp(V,\omega)}^\star.$

Using the Cartan--Killing form to identify $\g$ with its dual, 
we accordingly see that any semi-equicontinuous 
adjoint orbit is contained in the double cone 
$$ \oline{W_{\sp(V,\omega)}} \cup - \oline{W_{\sp(V,\omega)}}, $$ 
i.e., the corresponding Hamiltonian function is either 
positive or negative. 

Finally, we observe that if $x \in \g$ satisfies 
$\conv(\cO_x) \not=\g$, then $B(\cO_x)$ is a non-zero invariant 
cone and Proposition~\ref{prop:simpfin}(i) 
implies that it has interior points. 
Therefore $\cO_x$ is semi-equicontinuous because every finite dimensional 
space is barrelled. 
\end{prf}

\begin{lem} \mlabel{lem:ideal} The subspace 
$\fp := \{ X \in \sp(\cH) \: IX = - XI\}$
contains no non-trivial Lie algebra ideals, and the same holds 
for the subspace 
$$\fp_2 := \{ X \in \sp_{\rm res}(\cH) \: IX = - XI\} 
= \fp \cap \sp_{\rm res}(\cH).$$ 
\end{lem}

\begin{prf} First we show that $X \in \fp$ and $[X,[I,X]] = 0$ implies 
$X = 0$. Realizing $\sp(\cH)$ as a closed subalgebra of 
$\gl(\cH_\C) \cong \gl(\cH \oplus \cH)$ (cf.\ Example~\ref{ex:liegrp}(h)), 
we have 
$$ \fp = \Big\{ \pmat{ 0 & a \\ a^* & 0} \: a^\top = a  \in B(\cH)\Big\} 
\quad \mbox{ and } \quad 
I = \pmat{ i & 0 \\ 0 & -i}. $$
For $X = \pmat{ 0 & a \\ a^* & 0}$ this leads to 
\begin{align*}
[X,[I,X]] 
&= \Big[\pmat{ 0 & a \\ a^* & 0}, \pmat{ 0 & 2ia \\ -2ia^* & 0}\Big]
= -4i\pmat{ aa^*  & 0 \\ 0& a^* a}. 
\end{align*}
If this operator vanishes, then $a^* a = 0$ implies that 
$a = 0$, so that $X = 0$. 

If $\fii \subeq \fp$ is a Lie algebra ideal, 
then we have for each $X \in \fii$ the relation 
$[X,[I,X]] \in \fii \cap \fu(\cH) = \{0\}$, so that 
the preceding argument shows that $X = 0$. The same argument 
applies to $\fp_2 = \fp \cap \sp_{\rm res}(\cH)$. 
\end{prf}

\begin{lem} \mlabel{lem:posproj} Let $X \in \sp(\cH)$. Then 
the projection 
$p_\fk \: \sp(\cH) \to \fu(\cH)$
onto the $\C$-linear component satisfies 
$$ p_\fk(\cO_X) \subeq \{ y \in \fu(\cH) \: i y \leq 0 \} 
= \oline{C_{\fu(\cH)}}, $$
if and only if $H_X \geq 0$. A corresponding 
statement holds for $\sp_{\rm res}(\cH)$. 
\end{lem}

\begin{prf} Since $p_{\fk}$ is the fixed point projection for the 
action of the torus $e^{\R \ad I}$, it preserves the closed convex cones 
$\pm\{ Z \in \sp(\cH) \: H_Z \geq 0\}$ (Proposition~\ref{prop:project}), 
so that the relation 
$H_X \geq 0$ implies $i p_\fk(\cO_X)  \leq 0$ 
(Theorem~\ref{thm:coneconj}(i)). 
Suppose, conversely, that $i p_\fk(\cO_X)  \leq 0$. 
If $H_X \leq 0$, then the preceding argument 
shows that $i p_\fk(\cO_X) \geq 0$, which leads to 
$p_\fk(\cO_X) = \{0\}$, so that the closed span of $\cO_X$ is 
an ideal of $\g$ contained in $\fp = \ker p_\fk$. 
In view of Lemma~\ref{lem:ideal}, this leads to $X = 0$. 

We may therefore assume that  $H_X$ is indefinite. 
Hence there exists a finite dimensional 
complex subspace $\cH_1 \subeq \cH$ on which $H_X$ is indefinite. 
Then we have a Hilbert space direct sum 
$\cH = \cH_1 \oplus \cH_1^\bot$ and accordingly we write 
operators on $\cH$ as $(2 \times 2)$-block matrices: 
$$ X = \pmat{ X_{11} & X_{12}\\ X_{21} & X_{22}}. $$
For $v \in \cH_1$ we then have
$\Im \la Xv,v \ra = \Im \la X_{11}v,v \ra,$
so that $X_{11} \in \sp(\cH_1)$ and $H_{X_{11}}$ is indefinite. 

Identifying $\Sp(\cH_1)$ in the natural way with a subgroup 
of $\Sp(\cH)$, Lemma~\ref{lem:convall} implies that 
$\conv(\Ad(\Sp(\cH_1)X))$ contains an element 
$Y$ with $Y_{11} = - i \1$. 
From 
$$ p_\fk(Y) 
= \pmat{ p_\fk(Y_{11}) & p_\fk(Y_{12})\\ 
p_\fk(Y_{21}) & p_\fk(Y_{22})} 
= \pmat{ Y_{11} & p_\fk(Y_{12})\\ 
p_\fk(Y_{21}) & p_\fk(Y_{22})} $$ 
we now derive that $i p_\fk(Y) \not\leq 0$, contradicting 
our assumption on $X$. This completes the proof of the 
first assertion. 

Since $\Sp(\cH_1) \subeq \Sp_{\rm res}(\cH)$, the preceding 
argument also implies the assertion for the restricted 
Lie algebra $\sp_{\rm res}(\cH)$. 
\end{prf}

Now we are ready to show that the metaplectic representation 
is semibounded. We write 
$$q \: \hat\sp_{\rm res}(\cH) \to \sp_{\rm res}(\cH)$$ 
for the quotient map and $\z$ for its kernel. 
Since the central extension is trivial over the subalgebra 
$\fu(\cH)$ of $\sp_{\rm res}(\cH)$, we have 
$$ \hat\fu(\cH) := q^{-1}(\fu(\cH)) \cong \z \oplus \fu(\cH), $$
where the $\fu(\cH)$-complement is uniquely determined by the 
property that it is the commutator algebra. 
Here we use that the Lie algebra $\fu(\cH)$ is perfect 
(cf.~\cite[Lemma~I.3]{Ne02a}). 
Accordingly, we may identify $\fu(\cH)$ in a natural 
way with a subalgebra of $\hat\sp_{\rm res}(\cH)$.

\begin{rem} \mlabel{rem:9.13b} The momentum set $I_{\pi_s}$ 
is completely 
determined by the restriction of its support function 
$s_{\pi_s}$ to $W_{\pi_s} = B(I_{\pi_s})^0$ (Remark~\ref{rem:1.1}(b)) 
which is invariant under the adjoint action. 
We shall see below that $q(W_{\pi_s}) \subeq W_{\sp_{\rm res}(\cH)}$, 
so that Theorem~\ref{thm:coneconj}(ii) leads to 
$$ W_{\pi_s} \subeq \Ad(\hat\Sp_{\rm res}(\cH))(\fz \times C_{\fu(\cH)}), $$
where $\fz \cong \R$ denotes the center. This 
in turn entails that $I_{\pi_s}$ is already determined 
by the restriction of $s_{\pi_s}$ to $\z \times C_{\fu(\cH)}$, 
which we determine below. 

This restriction is the support 
function for the momentum set of the restriction of 
${\pi_s}$ to the subgroup $\hat\U(\cH) \cong \T \times \U(\cH)$ 
(Proposition~\ref{prop:repcone}(iii)), which is a 
direct sum of the representations on the subspaces $S^n(\cH)$, on 
which we have 
$$ {\pi_s}(z,g)(v_1 \vee \cdots \vee v_n) 
= z gv_1 \vee \cdots \vee g v_n \quad \mbox{ for } \quad 
v_1,\ldots, v_n \in \cH. $$
For the representation of $\U(\cH)$ on $S^n(\cH)$, the momentum 
set $I_{S^n(\cH)}^{\U(\cH)}$ is simply given by 
$$ I_{S^n(\cH)}^{\U(\cH)} = n I_{\cH}^{\U(\cH)} $$ 
(Example~\ref{ex:7.3}(b)), so that 
$$ I_{\pi_s}^{\hat U(\cH)} 
= \oline\conv\big(\bigcup_{n \in \N_0} \{1\} \times n I_\cH^{\U(\cH)}\big), $$
which is a convex cone with vertex $(1,0)$, and for 
$x \in\fu(\cH)$, we have 
\begin{equation}
  \label{eq:spi}
s_{\pi_s}(t,x) = 
\begin{cases}
-t & \text{ for } ix \leq 0 \\ 
\infty & \text{ else}. \\ 
\end{cases}
\end{equation}

The representation of $\hat\Sp_{\rm res}(\cH)$ on 
$S(\cH)$ decomposes into two irreducible pieces 
$$ S^{\rm even}(\cH) = \hat\oplus_{n \in \N_0} S^{2n}(\cH)
\quad \mbox{ and } \quad 
S^{\rm odd}(\cH) = \hat\oplus_{n \in \N_0} S^{2n+1}(\cH) $$ 
(Remark~\ref{rem:9.4}).  
Writing $\pi_\pm$ for the corresponding representations, we obtain 
$$ I_{\pi_+}^{\hat U(\cH)} = I_{{\pi_s}}^{\hat U(\cH)} 
= \oline\conv\Big(\bigcup_{n \in \N_0} \{1\} \times 2n I_\cH^{\U(\cH)}\Big), $$
and
$$ I_{\pi_-}^{\hat U(\cH)} 
= \oline\conv\Big(\bigcup_{n \in \N_0} \{1\} \times (2n+1)I_\cH^{\U(\cH)}
\Big). $$ 
Now it is easy to derive that 
$W_{\pi_\pm} = W_{\pi_s}$, $s_{\pi_+} = s_{\pi_s}$, and that 
$$ s_{\pi_-}(t,x) = -t + \sup(\Spec(ix)) \quad 
\mbox{ for } \quad  ix \leq 0. $$
\end{rem}

\begin{thm} \mlabel{thm:metaplec-sembo} 
The metaplectic representation 
$({\pi_s}, S(\cH))$ of 
$\hat\Sp_{\rm res}(\cH))$ is semibounded with 
$$ W_{\pi_s} = W_{\hat\sp_{\rm res}(\cH)} := q^{-1}(W_{\sp_{\rm res}(\cH)}). $$
\end{thm}

\begin{prf} From Example~\ref{ex:4.9} we know that 
$$ s_{\pi_s}(x) \leq 0 \quad \mbox{ for } \quad x \in C_{\fu(\cH)},
\quad \mbox{ so that } \quad 
s_{\pi_s}(t,x) = -t + s_{\pi_s}(x) < \infty $$
for $(t,x) \in \fz \times C_{\fu(\cH)}$. 
From the invariance of $s_{\pi_s}$ under the adjoint action 
and Theorem~\ref{thm:coneconj}(ii) it now follows that 
\begin{equation}
  \label{eq:contain}
\Ad(\hat\Sp_{\rm res}(\cH))(\fz \times C_{\fu(\cH)})
= q^{-1}(\Ad(\Sp_{\rm res}(\cH))C_{\fu(\cH)})
= q^{-1}(W_{\sp_{\rm res}(\cH)}) \subeq B(I_{\pi_s}) 
\end{equation}
(cf.\ Theorem~\ref{thm:coneconj}(ii)). 
In particular,  $B(I_{\pi_s})$ has interior points, 
so that ${\pi_s}$ is semibounded with $W_{\hat\sp_{\rm res}(\cH)} \subeq W_{\pi_s}$.

To prove equality, we note that \eqref{eq:spi} in 
Remark~\ref{rem:9.13b} implies that 
$$ W_{\pi_s} \cap \fu(\cH) \subeq 
\{ x \in \fu(\cH) \: s_{\pi_s}(x) < \infty\} 
= \{ x \in \fu(\cH) \: ix \leq 0\} = \oline{C_{\fu(\cH)}}.  $$ 
As $W_{\pi_s}$ is invariant under the projection 
$p_\fk \: \hat\sp_{\rm res}(\cH) \to \hat\fu(\cH)$ 
onto the fixed point space of the compact group 
$e^{\R \ad I}$ (Proposition~\ref{prop:project}), 
we obtain for each $x \in W_{\pi_s}$ the relation 
$i p_\fk(\cO_{q(x)})\leq 0,$ 
so that Lemma~\ref{lem:posproj} leads to 
$q(x) \in \oline{W_{\sp_{\rm res}(\cH)}}$. As $W_{\pi_s}$ is open and 
$W_{\sp_{\rm res}(\cH)}$ coincides with the interior of its closure 
(Lemma~\ref{lem:bou}), it follows that 
$$ q(W_{\pi_s}) \subeq W_{\sp_{\rm res}(\cH)}. $$ 
Combining this with \eqref{eq:contain}, 
we obtain $q(W_{\pi_s}) = W_{\sp_{\rm res}(\cH)}.$ 

As the center $\fz\cong \R$ of $\hat\sp_{\rm res}(\cH)$ 
acts by multiples of the identity on $S(\cH)$, 
we have $\z \subeq H(W_{\pi_s})$, which finally leads to 
$W_{\pi_s} = q^{-1}(q(W_{\pi_s})) =  W_{\hat\sp_{\rm res}(\cH)}.$ 
\end{prf}

\subsection{The momentum set of the metaplectic representation} 
\mlabel{sec:9.4}

Now that we have determined the cone $W_{\pi_s}$ for the metaplectic 
representation, we now apply the tools from Section~\ref{sec:5} 
to determine the momentum set for the representation of the 
central extension $\hat\HSp_{\rm res}(\cH)$ acting on 
$S(\cH)$. In particular, we show that the momentum set is the 
closed convex hull of a single coadjoint orbit. 

\begin{defn} We define 
$\hat\HSp_{\rm res}(\cH)$ 
\glossary{name={$\hat\HSp_{\rm res}(\cH)$},description={central extension of 
$\cH \rtimes \Sp_{\rm res}(\cH)$ containing the metaplectic group}}
as the quotient of the semidirect product 
$$ \Heis(\cH) \rtimes_\alpha \hat\Sp_{\rm res}(\cH), \quad 
\alpha(g,u)((t,v)) := (t,gv),  $$
by the central subgroup 
$$ S := \{ (t,(\1,z\1)) \in \R \times Z \: e^{it}z = 1\} \cong \R $$
which acts trivially on $S(\cH)$. 
This means that $\hat\HSp_{\rm res}(\cH)$ is a central extension 
of $\cH \rtimes \Sp_{\rm res}(\cH)$ by $\T$. The corresponding 
Lie algebra cocycle is given by 
$$((v,x), (v',x')) \mapsto \Im\la v,v'\ra + \eta(x,x') 
= \Im \la v,v'\ra + \shalf \Im \tr(x_2 x'_2) $$ 
(cf.\ Theorem~\ref{thm:metasmooth}). 
\end{defn} 

\begin{rem}\mlabel{rem:9.4.1}
The representations $W$ of $\Heis(\cH)$ and $\pi_s$ of 
$\hat\Sp_{\rm res}(\cH)$ combine to a representation, also denoted 
$\pi_s$ of 
$\hat\HSp_{\rm res}(\cH)$ on $S(\cH)$. Since $(W,S(\cH))$ is irreducible,
the extension to $\hat\HSp_{\rm res}(\cH)$ is also irreducible. 
To see that it is smooth, it suffices to show that the function 
$$ \Heis(\cH) \rtimes \hat\Sp_{\rm res}(\cH)\to \C, \quad 
(h,g) \mapsto 
\la W(h) \pi_s(g)\Omega, \Omega \ra 
= \la \pi_s(g)\Omega, W(h)^{-1}\Omega \ra  $$ 
is smooth (Theorem~\ref{thm:unitautsmooth}). 
This follows from the smoothness 
of the vector $\Omega$ for $\Heis(\cH)$ and $\hat\Sp_{\rm res}(\cH)$ 
(cf.\ \eqref{eq:matcoeff} and Theorem~\ref{thm:metasmooth}). 
\end{rem}

\begin{prop} The metaplectic representation $(\pi_s, S(\cH))$ 
of $\hat\HSp_{\rm res}(\cH)$ is semibounded with 
$$ W_{\pi_s} = q^{-1}(\cH \times W_{\sp_{\rm res}(\cH)}) 
= \Ad(\hat\HSp_{\rm res}(\cH))(\z \times C_{\fu(\cH)}), $$
where $q \: \hat\hsp_{\rm res}(\cH) \to \cH \rtimes \sp_{\rm res}(\cH)$ 
is the quotient map whose kernel is the center~$\fz \cong \R$.  
\end{prop}

\begin{prf} Let 
$$p \: \hsp_{\rm res}(\cH) = \heis(\cH) \rtimes \sp_{\rm res}(\cH) 
\to \cH \rtimes \sp_{\rm res}(\cH), \quad 
(t,v,x) \mapsto (v,x) $$
denote the quotient map, so that Theorem~\ref{thm:coneconj}(iv) implies that 
$$ W_{\hsp_{\rm res}(\cH)} = p^{-1}(W_{\sp_{\rm res}(\cH)})
= \Ad(\HSp_{\rm res}(\cH)) W_{\sp_{\rm res}(\cH)}. $$
Modulo the center, this relation leads to 
\begin{equation}
  \label{eq:conj-hsp}
\cH \rtimes W_{\sp_{\rm res}(\cH)} 
= \Ad(\cH \rtimes \Sp_{\rm res}(\cH))W_{\sp_{\rm res}(\cH)} 
= \Ad(\cH \rtimes \Sp_{\rm res}(\cH))C_{\fu(\cH)}. 
\end{equation}

From Theorem~\ref{thm:metaplec-sembo} we know that 
$s_{\pi_s}(x)$ is finite if $x \in \fz \times C_{\fu(\cH)}$, 
so that \eqref{eq:conj-hsp} implies that 
it is also finite if $q(x) \in \cH \rtimes W_{\sp_{\rm res}(\cH)}$, i.e., 
$$ W_{\pi_s} \supeq q^{-1}(\cH \rtimes W_{\sp_{\rm res}(\cH)}). $$
This proves already that $\pi_s$ is semibounded 
because $W_{\pi_s}$ has interior points (Proposition~\ref{prop:seqcrit}). 
We further conclude that $\heis(\cH) \subeq H(W_{\pi_s})$, so that 
we obtain with Theorem~\ref{thm:metaplec-sembo} that 
\begin{align*}
W_{\pi_s} 
&= \heis(\cH) + (W_{\pi_s} \cap \hat\sp_{\rm res}(\cH)) 
\subeq \heis(\cH) + q^{-1}(W_{\sp_{\rm res}(\cH)}) \\
&= q^{-1}(\cH \times W_{\sp_{\rm res}(\cH)}). 
\end{align*}
This proves the desired equality. 

Finally, we note that \eqref{eq:conj-hsp} implies that every element 
$x \in W_{\pi_s}$ is conjugate to an element $y \in q^{-1}(C_{\fu(\cH)})$, 
which means that $y \in \fz \times C_{\fu(\cH)}$. 
\end{prf}

\begin{thm} \mlabel{thm:9.14} 
The momentum set $I_{\pi_s}$ of the metaplectic representation 
of $\hat\HSp_{\rm res}(\cH)$ is the closed convex hull of the coadjoint 
orbit of $\lambda := \Phi_{\pi_s}([\Omega])$, and this linear functional 
is given by 
$$\lambda(t,x) = t\quad \mbox{ on } \quad 
\hat\hsp_{\rm res}(\cH) = \R \oplus (\cH \rtimes \sp_{\rm res}(\cH)).$$ 
\end{thm}

\begin{prf} Since the metaplectic representation of 
$G := \hat{\HSp}_{\rm res}(\cH)$ is irreducible, 
we want to apply Theorem~\ref{thm:5.7}. 
For the subgroup $K := \hat\U(\cH)\subeq G$ with Lie algebra 
$\fk = \z \oplus \fu(\cH)$, we have 
$$ \g_\C \cong \fp_+ \oplus \fk_\C \oplus \fp_-, 
\quad \mbox{ where } \quad 
\fp_\pm = \ker(\ad I \mp i\1)\oplus \ker(\ad I \mp  2i\1),  $$
$I \in \fu(\cH)$ is the multiplication with $i$ on $\cH$, 
$$ \heis(\cH)_\C = \fz_\C \oplus \cH_+ \oplus \cH_-, \quad 
\cH_\pm := \ker(\ad I \mp i\1)\subeq \cH_\C. $$

Since the Lie algebra $\fk$ of $K$ is complemented by the closed 
subspace $\cH \oplus \fp_2$ in $\g$, the coset space 
$G/K$ carries the structure of a Banach homogeneous space. 
Moreover, the closed $K$-invariant subalgebra 
$\fk_\C \oplus \fp_+$ determines on $G/K$ the structure of a 
complex manifold for which the tangent space in the base point 
can be identified with $\g_\C/(\fk_\C + \fp_+) \cong \fp_-$ 
(\cite[Thm.~6.1]{Bel06}). 

Let $\Omega \in S^0(\cH) \subeq S(\cH)$ be the vacuum vector. 
Then $\Omega$ is a smooth vector (Remark~\ref{rem:9.4.1}), 
which is an eigenvector for $K$. On $S(\cH)$, the operator 
$-i\dd\pi(I)$ is diagonal, and $S^n(\cH)$ is the eigenspace corresponding 
to the eigenvalue $n$. As $-i\ad I$ acts on $\fp_-$ with the 
eigenvalues $-1$ and $-2$, it follows that 
$\dd\pi_s(\fp_-)\Omega = \{0\},$
so that Theorem~\ref{thm:5.7} provides a holomorphic equivariant map 
$$ \eta \: G/K\to \bP(S(\cH)'), \quad 
g \mapsto \la \cdot, \pi_s(g)^{-1} \Omega \ra $$ 
and a realization of $(\pi_s, S(\cH))$ by holomorphic sections of a 
complex line bundle over $G/K$. Since the representation 
$(\pi_s, S(\cH))$ of $G$ is irreducible, $\Omega$ is a cyclic vector. 

It remains to show that for each $x \in W_\pi$ the flow generated by 
$-x$ on $G/K$ extends holomorphically to $\C_+$. Since every such 
element $x$ is conjugate to an element of $\fz \times C_{\fu(\cH)}$, we may 
w.l.o.g.\ assume that $x \in C_{\fu(\cH)}$. 

To get more information on $G/K$, 
we note that the choice of the complex 
structure implies that $G/K$ contains 
$\Heis(\cH)/(K \cap \Heis(\cH)) \cong \cH_-\cong \cH'$ 
as a complex submanifold. 
The translation action of $\Heis(\cH)$ by 
$(t,v).z := v + z$ on $\cH$ factors through an action of the additive group 
$\cH$ and extends naturally to a holomorphic action of the 
complexified group $\cH_\C \cong \cH_+ \times \cH_-$. 
As the $G$-action on $G/K$ induces a transitive 
action on the set of $\Heis(\cH)$-orbits on $G/K$, 
the action of $\Heis(\cH)$ extends to a holomorphic action of 
$\Heis(\cH)_\C$ on $G/K$. Therefore the action of the one-parameter group
$t \mapsto \exp(-tx)$ extends to a holomorphic action of $\C_+$ 
on $G/K$ if and only if the same holds for the action 
on the quotient space $\Sp_{\rm res}(\cH)/\U(\cH)$. 
As we have seen in Example~\ref{ex:liegrp}(h), 
this space can be identified with the symmetric Hilbert domain  
$$ \cD_s = \{ z \in B_2(\cH) \: z^\top = z, zz^* \< 1  \},  $$
but here it is endowed with the opposite complex structure. 

The group $\U(\cH)$ acts on $\cD_s$ by 
$u.z = uzu^\top$, and for $x \in C_{\fu(\cH)}$ 
we obtain for our choice of complex structure the relation 
$e^{-ix}.z  = e^{ix} ve^{ix^\top}.$
Now $\|e^{ix}\| = \|e^{ix^\top}\|< 1$ 
implies that the action of the one-parameter group 
$t \mapsto \exp(-tx)$ extends to $\C_+$, and this completes the 
proof. 

Finally, we derive from  Theorem~\ref{thm:5.7}(b) that 
$I_{\pi_s} = \oline{\conv}(\cO_\lambda)$ 
holds for 
$$ \lambda(x) = - i \la \dd\pi_s(x)\Omega, \Omega \ra. $$
On $\fk = \hat\fu(\cH)$ we have $\lambda(t,x) = t$, and $\lambda$ vanishes 
on $\fp_\C$ because $\fp_-$ annihilates $\Omega$, and 
$\la \dd\pi_s(\fp_+)\Omega , \Omega \ra 
= \la \Omega , \dd\pi_s(\fp_-)\Omega \ra = 0$. 
This completes the proof. 
\end{prf}

With similar and even easier arguments as in the proof of the 
preceding theorem, we also obtain: 
\begin{thm} \mlabel{thm:9.15} 
The momentum set $I_{\pi_s^+}$ of the even metaplectic representation 
$(\pi_s^+, S^{\rm even}(\cH))$ of $\hat\Sp_{\rm res}(\cH)$ 
is the closed convex hull of the coadjoint 
orbit of $\Phi_{\pi_s^+}([\Omega])$. 
\end{thm}

\section{The spin representation} \mlabel{sec:10} 

In this section we take a closer look 
at the spin representation 
$(\pi_a, \Lambda(\cH))$ of the 
central extension $\hat\rO_{\rm res}(\cH)$ of $\rO_{\rm res}(\cH)$ 
on the fermionic Fock space $\Lambda(\cH)$. 
This representation arises from self-intertwining 
operators of the Fock representation of the $C^*$-algebra 
$\Car(\cH)$. Here we also show that $\pi_a$ is semibounded and 
determine the corresponding cone $W_{\pi_a}$. 
For the irreducible representation of the 
identity component on the even part 
$\Lambda^{\rm even}(\cH)$, we show that the momentum set 
is the weak-$*$-closed convex hull of a single coadjoint orbit. 
 
\subsection{Semiboundedness of the  spin representation} \mlabel{sec:10.1} 

Let $\cH$ be a complex Hilbert space and 
write $\{a,b\} := ab + ba$ for the anticommutator of two elements 
of an associative algebra. The
{\it CAR-algebra} of $\cH$ is a $C^*$-algebra
$\Car(\cH)$, together with a continuous antilinear map
$a \: \cH \to \Car(\cH)$ satisfying the
{\it canonical anticommutation relations}
\begin{equation}
  \label{eq:car}
\{a(f), a(g)^*\}  = \la g,f \ra \1
\quad \hbox{ and } \quad  \{a(f),a(g)\} = 0
\quad \mbox{ for } \quad f, g \in \cH
\end{equation}
and which is generated by the image of $a$. This determines
$\Car(\cH)$ up to natural isomorphism (\cite[Thm.~5.2.8]{BR97}). 
We also write $a^*(f) := a(f)^*$, which defines a complex 
linear map $a^* \: \cH \to \Car(\cH)$. 

The orthogonal group $\OO(\cH)$ of the underlying real 
Hilbert space 
$$(\cH_\R,\beta), \quad \beta(v,w) = \Re \la v,w\ra,$$ 
acts by automorphisms on this $C^*$-algebra as follows. Writing 
a real linear isometry as 
$g = g_1 + g_2$, where $g_1$ is linear and $g_2$ is antilinear, the 
relations $gg^* = g^*g = \1$ turn into 
$$ g_1 g_1^* + g_2 g_2^* = \1 = g_1^*g_1 + g_2^* g_2 
\quad \mbox{ and } \quad 
g_1 g_2^* + g_2 g_1^* = 0 = g_2^* g_1 + g_1^*g_2. $$
These relations imply that 
$$ a_g \: \cH \to \CAR(\cH), \quad 
f \mapsto a(g_1 f) + a^*(g_2 f) $$ 
satisfies the same anticommutation relations, so that the 
universal property of $\CAR(\cH)$ implies the existence of a unique 
automorphism $\alpha_g$ with 
$$ \alpha_g(a(f)) = a_g(f) \quad \mbox{ for } \quad f \in \cH.$$
These automorphisms of the $\CAR(\cH)$ are called 
{\it Bogoliubov automorphisms}. They define an action of 
$\rO(\cH)$ on $\Car(\cH)$. In particular, the unitary group 
$\U(\cH) \subeq \OO(\cH)$ acts on $\CAR(\cH)$ by 
$\alpha_g(a(f)) = a(gf)$ for $f \in \cH$. 

\begin{rem} \mlabel{rem:10.1} 
The $C^*$-algebra $\Car(\cH)$ has a natural irreducible 
representation 
$(\pi_0, \Lambda(\cH))$ 
on the antisymmetric Fock space $\Lambda(\cH)$ (\cite[Prop.~5.2.2(3)]{BR97}). 
The image $a_0(f) := \pi_0(a(f))$ 
acts by $a_0(f)\Omega=0$ and 
$$ a_0(f)(f_1 \wedge \cdots\wedge f_n) 
= \sum_{j = 1}^n (-1)^{j-1} \la f_j, f \ra 
f_1 \wedge \cdots \wedge f_{j-1} \wedge f_{j+1} \wedge \cdots \wedge f_n.$$
Accordingly, we have 
$$ a_0^*(f)\Omega = f \quad \mbox{ and } \quad 
a_0^*(f)(f_1 \wedge \cdots \wedge f_n) 
= f \wedge f_1 \wedge \cdots \wedge f_n. $$

Let $Iv = iv$ denote the complex structure on $\cH$. The action of the 
restricted orthogonal group 
$$\OO_{\rm res}(\cH) 
\:= \{ g \in \OO(\cH) \: \|[g,I]\|_2 < \infty \} $$
on $\Car(\cH)$ preserves the equivalence class of the representation 
$\pi_0$, 
so that there is a projective unitary representation 
$\pi_a \: \rO_{\rm res}(\cH) \to \U(\Lambda(\cH))$, called the 
{\it spin representation}, satisfying 
\begin{equation}
  \label{eq:spin-int}
\pi_a(g) \pi_0(A)\pi_a(g)^* = \pi_0(\alpha_g A) \quad \mbox{ for }
 \quad g \in \OO_{\rm res}(\cH), A \in \Car(\cH) 
\end{equation}
(cf.\ \cite[Thm.~3, p.~35]{Ot95}, \cite{ShSt65}). 
\end{rem}

In analogy with Theorem~\ref{thm:metasmooth} we here obtain: 

\glossary{name={$(\pi_a, \Lambda(\cH))$},description={spin representation of metagonal group}}
\glossary{name={$\hat\rO_{\rm res}(\cH)$},description={metagonal group}}

\begin{thm}\mlabel{thm:spinsmooth} 
$\hat\rO_{\rm res}(\cH)$ is a Lie group and 
the spin representation is smooth. 
A Lie algebra cocycle $\eta$ defining $\hat\fo_{\rm res}(\cH)$ 
as an extension of $\fo_{\rm res}(\cH)$ by $\R$ is given by 
$$ \eta(x,y) = -\frac{1}{2i} \tr([x_2, y_2]). $$
\end{thm}

\begin{prf} With completely analogous arguments as in the proof of 
Theorem~\ref{thm:metasmooth} we find 
with Appendix~\ref{app:c} that any 
$g \in \rO_{\rm res}(\cH)$ for which $g_1$ is invertible 
has a unique lift $\pi_a(g) \in \U(\Lambda(\cH))$ with 
$$ \pi_a(g) \Omega= c(g) e^{-\hat T(g)} \quad \mbox{ for } \quad 
T(g) := g_2 g_1^{-1} \in \Aherm(\cH)_a,  
c(g) = \|e^{-\hat T(g)}\|^{-1} $$
(cf.\ Lemma~\ref{lem:hs2}). 

This implies that $\hat\rO_{\rm res}(\cH)$ is a Lie group, 
its representation $\pi_a$ on $\Lambda(\cH)$ is smooth, and 
that $\Omega$ is a smooth vector. 
Since $\rO_{\rm res}(\cH)$ acts smoothly on $\Car(\cH)$, 
the space $\Lambda(\cH)^\infty$ of smooth vectors is invariant under 
$\pi_0(\Car(\cH))$, so that the irreducibility of the 
representation of $\Car(\cH)$ on $\Lambda(\cH)$ 
implies the smoothness of $\pi_a$. 

For the Lie algebra cocycle defining $\hat\fo_{\rm res}(\cH)$, 
we find as in the proof of Theorem~\ref{thm:metasmooth} 
with Lemma~\ref{lem:hs2}(iii) 
\begin{align*}
 \eta(x,y) 
= 2 \Im \la \hat x_2, \hat y_2 \ra
= -\Im \tr(x_2y_2) = -\frac{1}{2i} \tr([x_2,y_2]). 
\end{align*}
\end{prf}

\begin{rem} \mlabel{rem:sigma} As in Subsection~\ref{sec:9.2}, 
we obtain a Banach--Lie group 
\begin{align*}
\rO_{1,2}(\cH) 
&:= \Big\{ g = g_1  + g_2 \in \rO(\cH)\: 
\|g_2 \|_2 < \infty, \|\1 - g_1\|_1 < \infty\} 
\end{align*}
with 
$$ \rO_{\rm res}(\cH) \cong (\rO_{1,2}(\cH) \rtimes \U(\cH))/N, 
\quad N \cong \U(\cH) \cap \rO_{1,2}(\cH) = \U_1(\cH) $$
(cf.\ \cite[Def.~IV.7]{Ne02a}). From \cite[Rem.~IV.14]{Ne02a} 
we recall that $\rO_{\rm res}(\cH)$ has two connected components 
and that its fundamental group is trivial. 
On the other hand 
$$ \pi_0(\rO_{1,2}(\cH)) \cong \Z/2\Z \cong \pi_1(\rO_{1,2}(\cH)) $$
(\cite[Prop.~12.4.2]{PS86}, \cite[Prop.~III.15]{Ne02a}), 
so that there exists a simply connected $2$-fold covering group 
$q \: \Spin_{1,2}(\cH)\to \rO_{1,2}(\cH)_0$ 
(cf.\ \cite{dH72} for a discussion of the smaller group 
$\Spin_1(\cH)$ covering $\SO_1(\cH) = \rO_1(\cH)_0$). 
As the 
inclusion $\U_{1}(\cH) \to \rO_{1,2}(\cH)$ induces a surjective homomorphism 
$$ \pi_1(\U_1(\cH)) \cong \Z \to \pi_1(\rO_{1,2}(\cH)) \cong \Z/2, $$ 
$\hat\U_1(\cH) := q^{-1}(\U_1(\cH))$ is the unique 
$2$-fold connected covering of 
$\U_1(\cH)$ from Subsection~\ref{sec:9.2}. 

We further have an embedding 
\begin{equation}
  \label{eq:o12emb}
\sigma \: \fo_{1,2}(\cH) \to \hat\fo_{\rm res}(\cH), \quad 
\sigma(x) := \Big(-\frac{1}{2i} \tr(x_1), x\Big)
\end{equation}
of Banach--Lie algebras. 
On the subgroup $\hat\U_1(\cH) \subeq \Spin_{1,2}(\cH)$, 
$\sigma$ integrates to a group homomorphism 
$\sigma_G(g) := (g, \sqrt{\det}(g)^{-1}\pi_s(g))$  
(cf.\ Subsection~\ref{sec:9.2})   
and $\sigma$ integrates to a morphism of Banach--Lie groups 
$\sigma_G \: \Spin_{1,2}(\cH) \to \hat\rO_{\rm res}(\cH).$
Combining this map with the canonical inclusion 
$\U(\cH) \into \hat\rO_{\rm res}(\cH)$, the equivariance of 
$\sigma_G$ under conjugation with unitary operators implies the 
existence of a homomorphism 
$$ \mu \: \Spin_{1,2}(\cH) \rtimes \U(\cH) \to \hat\rO_{\rm res}(\cH), 
\quad (g,u) \mapsto \sigma_G(g)u.$$ 
\end{rem} 

The following proposition is proved as Proposition~\ref{prop:olinemu}, 
using the representation $\pi_a$ instead of $\pi_s$. 
\begin{prop} The homomorphism $\mu$ factors through an 
isomorphism 
$$ \oline\mu \: (\Spin_{1,2}(\cH) \rtimes \U(\cH))/\ker \mu 
\to \hat\rO_{\rm res}(\cH)_0$$ 
of connected Banach--Lie groups with 
$\ker \mu \cong \SU_1(\cH).$ 
\end{prop}

\begin{rem} \mlabel{rem:10.6} 
The embedding $\fo_{1,2}(\cH) \into \hat\fo_{\rm res}(\cH)$ 
restricts in particular to an embedding $\fo_1(\cH) \into 
\hat\fo_{\rm res}(\cH)$. One can show that the operators 
$\dd\pi_a(x)$ are bounded for $x \in \fo_1(\cH)$ 
(cf.\ \cite{AW64}, \cite[Sect.~V]{Ne98}), so that we obtain 
a morphism of Banach--Lie algebras 
$$ \fo_1(\cH) \to \fu(\Lambda(\cH)), \quad 
x \mapsto  \dd\pi_a(x) - \frac{1}{2} \tr(x_1)\1. $$
\end{rem}

\begin{defn} To determine the open cone $W_{\pi_a}$ 
for the spin representation, we have to take a closer 
look at natural cones in $\fo_{\rm res}(\cH)$. 
Let $\fu_\infty(\cH)$ denote the ideal of compact skew-hermitian 
operators in $\fu(\cH)$ and recall the {\it Calkin algebra} 
$\Cal(\cH) := B(\cH)/B_\infty(\cH)$, where 
$B_\infty(\cH)$ denotes the ideal of compact operators 
in the $C^*$-algebra $B(\cH)$. Then the surjection 
$B(\cH) \to \Cal(B(\cH))$ induces an isomorphism 
$$ \fu(\cH)/\fu_\infty(\cH) \cong \fu(\Cal(\cH)). $$
Next we note that 
$\{ x \in \fo_{\rm res}(\cH) \: x^* - x \in \fu_\infty(\cH)\}$ 
is a closed ideal of $\fo_{\rm res}(\cH)$ 
which defines a quotient morphism 
$q \: \fo_{\rm res}(\cH) \to \fu(\Cal(\cH)).$
We thus obtain an open invariant cone in $\fo_{\rm res}(\cH)$ 
by 
$$ W_{\fo_{\rm res}(\cH)} := q^{-1}(C_{\fu(\Cal(\cH))}) 
= \{ x \in \fo_{\rm res}(\cH) \: -i q(x) \< 0 \} $$
(cf. Example~\ref{ex:liealgcon}). We write 
$W_{\hat\fo_{\rm res}(\cH)}$ for its inverse image in the 
central extension $\hat\fo_{\rm res}(\cH)$. 
\end{defn}

\begin{thm} \mlabel{thm:spin-sembo} 
The spin representation $(\pi_a, \Lambda(\cH))$ of 
the connected Lie group $\hat \rO_{\rm res}(\cH)$ is semibounded 
with $W_{\pi_a} = W_{\hat\fo_{\rm res}(\cH)}$. 
\end{thm} 

\begin{prf} First we note that under 
$\hat\U(\cH)\cong \T \times \U(\cH)$ the spin representation decomposes 
into the irreducible subspaces $\Lambda^n(\cH)$ with 
$$ {\pi_a}(z,g)(v_1 \wedge \cdots \wedge v_n) 
= z gv_1 \wedge \cdots \wedge g v_n \quad \mbox{ for } \quad 
v_1,\ldots, v_n \in \cH. $$
If $i x \leq 0$, then the corresponding operator 
on $\Lambda^n(\cH) \subeq \cH^{\otimes n}$ is also 
$\leq 0$, and this implies that 
$$ s_{\pi_a}(t,x) = -t + s_{\pi_a}(x) \leq - t < \infty 
\quad \mbox{ for } \quad 
x \in C_{\fu(\cH)}. $$
Writing $\z \cong \R$ for the center of $\hat\fo_{\rm res}(\cH)$, 
we thus find that 
$\fz \times C_{\fu(\cH)} \subeq B(I_{\pi_a}).$

Next we consider the decomposition of $\g = \hat\fo_{\rm res}(\cH) 
= \fk \oplus \fp_2$, where $\fk = \fz \times \fu(\cH)$ and 
$\fp_2 = \{ x \in \fo_{\rm res}(\cH) \: Ix = - xI\}.$
For the map 
$$ F \: \fk \times \fp_2 \to \g, \quad (x,y) \mapsto 
e^{\ad y} x \quad \mbox{ we then have }\quad 
\dd F(x,0)(a,b) = a + [b,x], $$
which is invertible if and only if $\ad x \: \fp_2 \to \fp_2$ is 
invertible. This is in particular the case for 
$x = I$, the complex structure of $\cH$. Therefore the image of 
\break $(\fz \times C_{\fu(\cH)}) \times \fp_2$ has interior points, 
and since it is contained in the invariant subset 
$B(I_{\pi_a})$, the cone 
$W_{\pi_a} = B(I_{\pi_a})^0$ is non-empty. 
With Proposition~\ref{prop:seqcrit} we now see that the 
representation $(\pi_a, \Lambda(\cH))$ is semibounded. 

To determine the cone $W_{\pi_a}$, we recall from Remark~\ref{rem:10.6} 
that for each $x \in \fo_1(\cH)$, the corresponding 
operator $\dd\pi_a(x)$ on $\Lambda(\cH)$ is bounded. This means that 
$\fo_1(\cH) \subeq B(I_{\pi_a}) \subeq \oline{W_{\pi_a}}$, and 
since $H(W_{\pi_a})$ is a closed ideal (Lemma~\ref{lem:limcone}(v)), 
it follows that 
$$ \fn := \fz \oplus \fu_\infty(\cH) \oplus \fp_2 \subeq H(W_{\pi_a}). $$
In view of  $\hat\fo_{\rm res}(\cH) = \fu(\cH) + \fn$, this proves that 
$$W_{\hat\fo_{\rm res}(\cH)} = C_{\fu(\cH)} + \fn \subeq W_{\pi_a} 
= \fn + (W_{\pi_a} \cap \fu(\cH)). $$
It therefore remains to show that 
$$ W_{\pi_a} \cap \fu(\cH) \subeq C_{\fu(\cH)} + \fu_\infty(\cH). $$

To verify this assertion, let $A = A^*$ be a hermitian 
operator on $\cH$ and $\dd\pi_a(A)$ denote the corresponding operator 
on $\Lambda^n(\cH)$. Suppose that there exists an $\eps > 0$ 
such that the range of the spectral projection 
$P([\eps,\infty[)$ is infinite dimensional. Then there exists 
for each $n$ an orthonormal subset $v_1, \ldots, v_n$ in this space. 
We then have 
\begin{align*}
&\ \ \ \ \la \dd\pi_a(A)(v_1 \wedge \cdots \wedge v_n), 
v_1 \wedge \cdots \wedge v_n \ra \\
&= \sum_{j = 1}^n \la v_1 \wedge \cdots \wedge Av_j \wedge \cdots \wedge v_n, 
v_1 \wedge \cdots \wedge v_n \ra 
= \sum_{j = 1}^n \la Av_j, v_j \ra \geq n \eps, 
\end{align*}
and this implies that $\dd\pi_a(A)$ is not bounded from above. 
If, conversely, for every $\eps > 0$ 
the spectral projection $P([\eps,\infty[)$ has finite dimensional 
range, then $A = A_- + A_+$, where 
$A_- \leq 0$ and $A_+$ is compact, i.e., 
$iA \in \oline{C_{\fu(\cH)}} + \fu_\infty(\cH)$.  
This implies that the open cone $W_{\pi_a} \cap \fu(\cH)$ 
is contained in $C_{\fu(\cH)} + \fu_\infty(\cH)$.  
\end{prf}

\begin{rem} \mlabel{rem:10.8} 
(a) With a similar argument as in Remark~\ref{rem:9.4} ones 
argues that under the action of the identity component 
$\hat\rO_{\rm res}(\cH)_0$ on the space $\Lambda(\cH)$ 
decomposes into two irreducible subrepresentations 
$\Lambda^{\rm even}(\cH)$ and $\Lambda^{\rm odd}(\cH)$. 
However, the action of the full group $\hat\rO_{\rm res}(\cH)$ 
is irreducible because the elements $g \in \hat\rO_{\rm res}(\cH)$ 
not contained in the identity component exchange 
the two subspaces 
$\Lambda^{\rm even}(\cH)$ and  
$\Lambda^{\rm odd}(\cH)$ (\cite[p.~239]{PS86}). 

(b) There is also an analog of the Banach--Lie algebra 
$\hsp_{\rm res}(\cH)$ acting irreducibly on $\Lambda(\cH)$ 
in the fermionic case. Here it is an infinite dimensional 
analog of the odd orthogonal Lie algebra $\so_{2n+1}(\R)$. 
To see this Lie algebra, we define linear and antilinear 
rank-one operators on $\cH$ by 
$$ P_{v,w}(x) := \la x,w \ra v \quad \mbox{ and } \quad 
\oline P_{v,w}(x) := \la w,x\ra v $$ 
and observe that $P_{v,w}^* = P_{w,v}$ and 
$\oline P_{v,w}^* = \oline P_{w,v}$. 
Therefore $Q_{v,w} := P_{v,w} - P_{w,v} \in \fu_1(\cH)$, 
and a direct calculation 
yields 
$$ \dd\pi_a(Q_{v,w}) = a_0^*(v)a_0(w) - a_0^*(w)a_0(v). $$ 
As $\tr(Q_{v,w}) = 2i \Im \la v,w\ra$, we obtain with the map 
$\sigma$ from \eqref{eq:o12emb} (Remark~\ref{rem:sigma}) 
$$ \dd\pi_a(\sigma(Q_{v,w})) = \dd\pi_a(-\Im \la v,w \ra, Q_{v,w}). $$

For the antilinear operators $\oline Q_{v,w} := 
\oline P_{v,w} - \oline P_{w,v}$ we have 
$\oline Q_{v,w} \in \fo_1(\cH) \subeq \fo_{1,2}(\cH)$ and 
$$ \dd\pi_a(\oline Q_{v,w}) 
= a_0^*(v)a_0^*(w) - a_0(w) a_0(v) = a_0^*(v)a_0^*(w) + a_0(v) a_0(w). $$
Next we observe that the operators $\rho(v) := \frac{1}{\sqrt 2}
(a_0(v) -a_0^*(v))$ satisfy 
\begin{align*}
&\ \ \ \ [\rho(v), \rho(w)]
= \shalf [a_0(v) - a_0^*(v), a_0(w) - a_0^*(w)] \\
&= (a_0(v)a_0(w) + a_0^*(v)a_0^*(w)) 
+ (a_0^*(w) a_0(v) - a_0^*(v)a_0(w)) + i\Im \la v, w \ra \\ 
&= \dd\pi_a(\Im \la v, w \ra, \oline Q_{v,w}-Q_{v,w})= 
 \dd\pi_a(\sigma(\oline Q_{v,w}- Q_{v,w})).
\end{align*}
This calculation implies that we obtain on the direct sum 
$\cH \oplus \fo_{1,2}(\cH)$ a Banach--Lie algebra structure 
with the bracket 
$$ [(v,X), (v',X')] 
:= (Xv' - X'v, [X,X'] + \oline Q_{v,w}-Q_{v,w}), $$ 
and 
$$ (v,X) \mapsto \rho(v) + \dd\pi_a(\sigma(X)) $$
defines a representation by unbounded operators on $\Lambda(\cH)$, 
where the subalgebra $\cH \oplus \fo_1(\cH)$ is represented by bounded 
operators. A closer inspection shows that 
$\cH \oplus \fo_1(\cH) \cong \fo_1(\cH \oplus \R)$ is an 
infinite dimensional version of $\fo_{2n+1}(\R)$ 
(cf.\ \cite[p.~215]{Ne98}). 
\end{rem}

\begin{defn} Let $\cI_\beta \subeq \GL(\cH_\R)$ be the set of 
orthogonal real-linear complex structures on the real 
Hilbert space $\cH_\R$. This set parameterizes the 
complex Hilbert space structures on $\cH_\R$ compatible 
with the given real Hilbert space structure. 
\end{defn}

\begin{lem}\mlabel{lem:o.1} The following assertions hold: 
  \begin{description}
  \item[\rm(i)] $\cI_\beta = \rO(\cH) \cap \fo(\cH) 
= \{ g \in \GL(\cH_\R) \: g^\top = g^{-1} = -g \}$ is a submanifold 
of $\rO(\cH)$. 
  \item[\rm(ii)] The conjugation action of $\rO(\cH)$ on 
$\cI_\beta$ leads to a diffeomorphism $\cI_\beta \cong \rO(\cH)/\U(\cH)$. 
\item[\rm(iii)] $\cI_\beta^{\rm res} := \Ad(\rO_{\rm res}(\cH))I 
=  \{ J \in \cI_\beta \: \|I - J \|_2 < \infty\} 
\cong \rO_{\rm res}(\cH)/\U(\cH)$.  
  \end{description}
\end{lem}

\begin{prf} (i) It only remains to show that $\cI_\beta$ is a 
submanifold of $\rO(\cH)$. For $J \in \cI_\beta$ we parameterize 
a neighborhood of $J$ by the map 
$\fo(\cH) \to \rO(\cH),x \mapsto J e^x$
which is a diffeomorphism on some open $0$-neighborhood $U 
\subeq \fo(\cH)$, which we may assume to be invariant under 
$\Ad(J)$. That $J e^x$ is a complex structure is equivalent to 
$$ -Je^x = (Je^x)^{-1} = e^{-x} J^{-1} = - J e^{-J^{-1}xJ}, $$
which is equivalent to 
$Jx = - xJ$. We conclude that $\cI_\beta$ is a submanifold 
of $\rO(\cH)$ whose tangent space can be identified 
with the set of $J$-antilinear elements in $\fo(\cH)$. 

(ii) Let $I \in \cI_\beta$ be the canonical complex structure 
given by $Iv = iv$ and $J \in \cI_\beta$. 
Then $\cH = (\cH_\R, I)$ and $(\cH_\R, J)$ are two complex Hilbert 
spaces whose underlying real Hilbert spaces are isomorphic. 
This implies that they have the same complex Hilbert dimension, 
i.e., there exists a unitary isomorphism 
$g \: (\cH_\R, I) \to (\cH_\R, J)$, i.e., $g \in \rO(\cH)$ with 
$gI = Jg$. Therefore $\rO(\cH)$ acts transitively on 
$\cI_\beta$ and the stabilizer of $I$ is the subgroup 
$\U(\cH)$ (cf.\ \cite[Lemma 1]{BMV68}). 
Since its Lie algebra $\fu(\cH)$ is complemented by 
the closed subspace 
$$\fp := \{ x \in \fo(\cH) \: Ix = - xI\},$$ 
it follows that $\rO(\cH)/\U(\cH)$ is a Banach homogeneous space 
diffeomorphic to~$\cI_\beta$.

(iii) For $g = g_1 + g_2 \in \rO(\cH)$ the operator 
$[I,g] = [I, g_2] = 2 I g_2$ is Hilbert--Schmidt if and only if 
$g_2$ is, i.e., $g \in \rO_{\rm res}(\cH)$. 
This in turn is equivalent to 
$gIg^{-1} - I = [g,I]g^{-1}$ being Hilbert--Schmidt, so that 
$$ \rO_{\rm res}(\cH) = \{ g \in \rO(\cH) \: 
\|gIg^{-1} - I \|_2 < \infty\} $$
and therefore (ii) implies (iii). 
\end{prf}

\begin{rem} A priori, the $C^*$-algebra 
$\Car(\cH)$ depends on the complex structure on $\cH$, but 
it can also be expressed as the $C^*$-algebra generated by 
the hermitian elements $b(f) := \frac{1}{\sqrt 2}(a(f) + a^*(f))$ satisfying 
$$ \{b(f), b(g)\} 
= \shalf \{a(f), a^*(g)\} + \shalf \{a^*(f), a(g)\} 
= \shalf(\la g, f \ra + \la f, g\ra) \1 =  \beta(f,g) \1. $$
In this sense it is the $C^*$-envelope of the Clifford algebra 
of $(\cH_\R,\beta)$. 
From this point of view it is even more transparent why the group 
$\rO(\cH)$ acts by automorphisms. 

Now we can think of the Fock representation $\pi_0$ as 
depending on the complex structure $I$ on $\cH$, and any other 
representation $\pi_0 \circ \alpha_g$ is the Fock representation 
corresponding to the complex structure $gIg^{-1} \in \cI_\beta$. 
We thus obtain a map of $\cI_\beta$ into the set $\Ext(S(\Car(\cH)))$ 
of pure states of $\Car(\cH)$, mapping 
$J := gIg^{-1}$ to $\omega_J(A) := \la \pi_0(\alpha_g^{-1}(A))\Omega, 
\Omega\ra$. These states are called the {\it Fock states} of 
$\Car(\cH)$, and in \cite[Thm.~3]{BMV68} they are essentially 
characterized as the pure quasi-free states (cf.\ Subsection~\ref{sec:10.2} below). 
\end{rem}

\begin{thm} \mlabel{thm:spinmom} For the irreducible subrepresentation 
$(\pi_a^+, \Lambda^{\rm even}(\cH))$ of $\hat\rO_{\rm res}(\cH)$ 
the momentum set is the closed convex hull of the coadjoint orbit of 
$\lambda := \Phi_{\pi_a}([\Omega])$. This functional is given 
by 
$$\lambda(t,x) = t\quad \mbox{ on } \quad 
\hat\fo_{\rm res}(\cH) = \R \oplus \fo_{\rm res}(\cH). $$
\end{thm}

\begin{prf} We want to apply Theorem~\ref{thm:5.7}. 
For the subgroup $K := \hat\U(\cH)\subeq G := \hat\rO_{\rm res}(\cH)_0$ 
with Lie algebra 
$\fk = \z \oplus \fu(\cH)$, we have 
$$ \g_\C \cong \fp_+ \oplus \fk_\C \oplus \fp_-, 
\quad \mbox{ where } \quad 
\fp_\pm = \ker(\ad I \mp  2i\1)  $$
and $I \in \fu(\cH)$ is the multiplication with $i$ on $\cH$. 

Since the Lie algebra $\fk$ of $K$ is complemented by the closed 
subspace $\fp_2$ of antilinear elements in $\g$, the coset space 
$G/K$ carries the structure of a Banach homogeneous space 
(cf.\ Lemma~\ref{lem:o.1}). 
Moreover, the closed $K$-invariant subalgebra 
$\fk_\C \oplus \fp_+$ determines on $G/K$ the structure of a 
complex manifold for which the tangent space in the base point 
can be identified with $\g_\C/(\fk_\C + \fp_+) \cong \fp_-$ 
(\cite[Thm.~6.1]{Bel06}). 

In view of \eqref{eq:ocomplex}, we obtain with 
\cite[Prop.~V.8, Rem.~V.10(c)]{Ne02a} that the 
group $\rO_{\rm res}(\cH)$ acts transitively on the homogeneous 
space $\rO_{\rm res}(\cH_\C,\beta_\C)/P$, 
where $P$ is the stabilizer of the subspace $\cH \oplus \{0\} 
\subeq \cH_\C$ (cf.\ Example~\ref{ex:liegrp}(d),(i)). 
From $P \cap \rO_{\rm res}(\cH) = \U(\cH)$ it now follows that 
the action of the Banach--Lie group $G$ on $G/K$ actually extends to a 
holomorphic action of a complex group. 

Let $\Omega \in \Lambda^0(\cH)$ be the vacuum vector. 
Then $\Omega$ is a smooth vector 
which is an eigenvector for $K$. On $\Lambda(\cH)$, the operator 
$-i\dd\pi(I)$ is diagonal, and $\Lambda^n(\cH)$ is the eigenspace 
corresponding to the eigenvalue $n$. As $-i\ad I$ acts on $\fp_-$ with the 
eigenvalue $-2$, it follows that 
$\dd\pi_a(\fp_-)\Omega = \{0\},$
so that Theorem~\ref{thm:5.7} provides a holomorphic equivariant map 
$$ \eta \: G/K\to \bP(\Lambda(\cH)'), \quad 
g \mapsto \la \cdot, \pi_a(g)^{-1} \Omega \ra $$ 
and a realization of $(\pi_a, \Lambda^{\rm even}(\cH))$ 
by holomorphic sections of a complex line bundle over $G/K$. 
Since the representation 
$(\pi_a, \Lambda^{\rm even}(\cH))$ of $G$ 
is irreducible, $\Omega$ is a cyclic vector. 
As we have argued above, the action of every one-parameter subgroup of 
$G$ extends to a holomorphic action of $\C$ on $G/K$, 
so that Theorem~\ref{thm:5.7}(c) implies that 
$I_\pi$ is the closed convex hull of the coadjoint orbit 
of $\lambda$. That $\lambda$ has the desired form follows 
as in the proof of Theorem~\ref{thm:9.14}. 
\end{prf}

\subsection{Quasi-free representations} \mlabel{sec:10.2}

Let $P = P^* = P^2$ be an 
orthogonal projection on $\cH$, 
$\cH_- := \im P$, $\cH_+ := \ker P$, and 
$\Gamma \: \cH \to \cH$ be an isometric antilinear involution  
commuting with $P$. Then 
$\tau_P := (\1 - P) + \Gamma P \in \OO(\cH)$. 
We write $a_0 \: \cH \to B(\Lambda(\cH))$ for the map corresponding to 
the Fock representation of $\Car(\cH)$ 
(Remark~\ref{rem:10.1}). 
Twisting with the Bogoliubov automorphism defined by $\tau_P$, we 
obtain an irreducible 
representation $(\pi_P, \Lambda(\cH))$ of $\Car(\cH)$ by 
$$ a_P(f) := a_0((\1-P)f) + a_0^*(\Gamma Pf) 
= \pi_0(a_{\tau_P}(f))  \quad \mbox{ for } 
\quad f \in \cH, $$  
i.e., $\pi_P = \pi_0 \circ \alpha_{\tau_P}$. 
These representations are called {\it quasi-free}. 
For $P =0$ we recover the Fock representation defined by~$a_0$.  
Two quasi-free representations $a_P$ and $a_Q$ are equivalent if and only if 
$\|P - Q\|_2 < \infty$ (\cite{PoSt70}). 

\begin{rem} The physical interpretation of $P$ is that its 
range $\cH_-$ consists 
of the negative energy states and its kernel $\cH_+$ consists of the 
positive energy states. 
For $f \in \cH_+$ we have 
$a_P(f)\Omega = 0$ and likewise 
$a_P^*(f)\Omega = 0$ for $f \in \cH_-$. 
This is interpreted in such a way that the 
creation operator $a_P^*(f)$ cannot create an additional negative 
state from $\Omega$ because all negative states are already 
filled (Pauli's Principle). 
Likewise, the annihilation operator $a_P(f)$, corresponding 
to the positive energy vector $f$, cannot extract any positive 
energy state from $\Omega$. 
\end{rem}

The restricted unitary group $\U_{\rm res}(\cH,P)$ 
(cf.\ Example~\ref{ex:liegrp}(c)) is a subgroup of $\rO(\cH)$, 
and for all projections $Q = gPg^{-1}$, $g \in \U_{\rm res}(\cH,P)$ 
we have \break $\|P - Q\|_2 < \infty$, so that the equivalence of 
$a_P$ and $a_Q$ leads to a projective unitary representation 
of $\U_{\rm res}(\cH)$ on $\Lambda(\cH)$ determined by 
$$ \pi_a^P(g) a_P(f) \pi_a^P(g)^* = a_P(gf), \quad f \in \cH, 
g \in \U_{\rm res}(\cH). $$
Let $\hat\U_{\rm res}(\cH)$ denote the corresponding central 
extension and write $\pi_a^P$ for its unitary representation 
on $\Lambda(\cH)$. To see that we thus obtain a semibounded 
representation of a Lie group, we first note that 
we have an embedding 
$$\iota \: \U(\cH) \to \rO(\cH), \quad 
g \mapsto \tau_P g \tau_P $$ 
for which $\U_{\rm res}(\cH,P)$ is precisely the inverse image 
of $\rO_{\rm res}(\cH)$ and that for $g \in \U_{\rm res}(\cH,P)$ we have 
\begin{align*}
\pi_a(\tau_P g \tau_P)a_P(f)\pi_a(\tau_P g \tau_P)^*
&= \pi_a(\tau_P g \tau_P)\pi_0(\alpha_{\tau_P} a(f))\pi_a(\tau_P g \tau_P)^*\\
&= \pi_0(\alpha_{\tau_P g} a(f)) = \pi_P(a_g(f)) = a_P(gf), 
\end{align*}
so that the projective representation of $\U_{\rm res}(\cH,P)$ on 
$\Lambda(\cH)$ coincide with $\pi_a \circ \iota$ 
(cf.\ \cite[p.~53]{Ot95}). It follows in particular that 
$\hat\U_{\rm res}(\cH,P) \cong \iota^*\hat\rO_{\rm res}(\cH)$ 
is a Lie group, $\pi_a^P$ is a smooth representation, 
and since $\im(\L(\iota))$ intersects the cone $C_{\fu(\cH)} 
\subeq W_{\pi_a}$, the representation $\pi_a^P$ is semibounded 
(cf.\ Proposition~\ref{prop:repcone}(iv)). 
More precisely, each element $x = (x_+, x_-) \in 
\fu(\cH_+) \oplus \fu(\cH_-)$ with 
$i x_+ \< 0$ and $ix_- \> 0$ lies in $W_{\pi_a^P}$. 
This holds in particular for $\tau_p I \tau_p = I(\1 - 2P)$. 

In the physics literature, the 
corresponding selfadjoint operator 
$$Q := -i \dd\pi_a^P(I) = \dd\pi_a(\1 - 2P)$$ is called the {\it 
charge operator} and $\1 - 2P$ the {\it one-particle charge operator}. 
Its $1$-eigenspace is $\cH_+$ and its $-1$-eigenspace is $\cH_-$. 

\begin{rem} The situation described above can be viewed as a 
``second quantization'' procedure that can be used to 
turn a self-adjoint operator $A$ on the single particle space 
$\cH$ into a non-negative operator $\hat A$ on the many particle 
space. In fact, let $P := P(]-\infty,0[)$ denote the 
spectral projection of $A$ corresponding to the open negative axis.
Then $iA$ generates a strongly continuous one-parameter group 
$\gamma_A(t) := e^{itA}$ 
of $\U(\cH_+) \times \U(\cH_-) \subeq \U_{\rm res}(\cH,P)$, 
and $\pi_a^P(\gamma_A(t)) = e^{it \hat A}$ is a one-parameter group 
of $\U(\Lambda(\cH))$ whose infinitesimal generator 
$\hat A$ has non-negative spectrum. 
\end{rem}

For more details on the complex manifolds 
$\rO_{\rm res}(\cH)/\U(\cH)$ (the isotropic restricted Gra\3mannian) 
and $\U_{\rm res}(\cH,P)/(\U(P(\cH)) \times \U(\cH_+))$ 
(the restricted or Sato--Segal--Wilson Gra\3mannian) we refer to 
\cite{PS86} and \cite{SW98}, where one also finds a detailed 
discussion of the corresponding complex line bundles, 
the Pfaffian line bundle ${\rm Pf}$ over $\rO_{\rm res}(\cH)/\U(\cH)$
whose dual permits the even spin representation as a space of 
holomorphic sections, and the determinant line bundle ${\rm Det}$ 
over the restricted Gra\3mannian whose dual provides a 
realization of the representations of $\hat\U_{\rm res}(\cH,P)$ 
mentioned above. Physical aspects of highest weight representations 
of $\U_{\rm res}(\cH)$ are discussed in \cite{CL02}. 

\section{Perspectives}  \mlabel{sec:11}

{\bf Classification problems:} 
In the preceding sections we mainly discussed three 
prototypical classes of semibounded representations: 
highest weight representations of the Virasoro group, 
the metaplectic representation on the bosonic Fock space 
and the spin representation on the fermionic Fock space. 

These representations are of fundamental importance in 
mathematical physics and homomorphisms from a Lie group 
$G$ to $\Sp_{\rm res}(\cH)$,  $\rO_{\rm res}(\cH)$ or 
$\U_{\rm res}(\cH,P)$ 
can be used to obtain semibounded representations of a central 
extension by pulling back the representations discussed above. 
This construction has been a major source of (projective) representations 
for loop groups $G = C^\infty(\bS^1,K)$ and $\Diff(\bS^1)$ 
(cf.\ \cite{SeG81}, \cite{PS86}, \cite{Ca83}, \cite{CL02}).

What still remains to be developed is a better global perspective 
on semibounded representations, including classification 
results on irreducible ones and the existence of direct integral 
decompositions. 

\begin{prob} Classify all irreducible semibounded unitary representations 
of the groups $\hat\Sp_{\rm res}(\cH)$, $\hat\rO_{\rm res}(\cH)$ 
and $\hat\U_{\rm res}(\cH,P)$, where $\cH$ is a complex Hilbert space and 
$P$ an orthogonal projection on $\cH$. 

These classification problems are special cases of the more general 
problem of the classification of the projective semibounded 
unitary representations of the automorphism groups 
of hermitian Hilbert symmetric spaces. For the three groups above, 
the corresponding spaces of 
$\cI_\omega^{\rm res}$, $\cI_\beta^{\rm res}$, resp., 
the restricted Gra\3mannian $\Gr_{\rm res}(\cH,P) := 
\{ gPg^{-1} \: g \in \U_{\rm res}(\cH,P)\}$. A crucial 
difference between these spaces is that the first one 
is equivalent to a symmetric Hilbert domain 
(of negative curvature) and the latter two are positively 
curved spaces. This  difference is also reflected in the difference 
between the cones $W_{\pi_s}$ and $W_{\pi_s}$ 
in $\hat\sp_{\rm res}(\cH)$ and $\hat\fo_{\rm res}(\cH)$  
(Theorems~\ref{thm:metaplec-sembo} and \ref{thm:spin-sembo}). 

An even larger class of groups arises as automorphism groups 
of Hilbert flag manifolds such as the orbit of a finite 
flag in $\cH$ under the group $\U_2(\cH)$. For a systematic 
discussion of these manifolds and the topology of the 
corresponding real and complex groups we refer to 
\cite[Sect.~V]{Ne02a} and for corresponding representations 
to \cite{Ne04}. 
\end{prob}

For loop groups of the form $C^\infty(\bS^1,K)$, 
$K$ a compact connected Lie group, the irreducible 
projective positive energy representations can be identified 
as highest weight representations (\cite[Thm.~11.2.3]{PS86}, 
\cite[Cor.~VII.2]{Ne01b}, and in particular 
\cite[Prop.~3.1]{SeG81} for $K = \T$). 
With the convexity theorems in \cite{AP83} and \cite{KP84} it should 
be possible to show that these representations 
are semibounded and one can hope for an analog of 
Theorem~\ref{thm:8.22} asserting that every irreducible 
semibounded representation of the corresponding double extension 
either is a highest weight representation or its dual. 
In \cite{Ne09b} it is shown on the algebraic level that there 
also exist many interesting unitary representations 
of infinite rank generalizations of twisted and untwisted 
loop algebras, resp., their double extensions 
(so-called locally affine Lie algebras). On the group level 
they correspond to (double extensions of) groups 
of the form $C^\infty(\bS^1, K)$, where 
$K$ is a Hilbert--Lie group, such as $\U_2(\cH)$. 
It seems quite likely that all these representations 
are semibounded with momentum sets generated by a single 
coadjoint orbit. 

{\bf New sources of semibounded representations:} 
Although many interesting classes of semibounded representations 
are known, a systematic understanding of the geometric sources 
of these representations is still lacking. 
As we have seen above, on the Lie algebra level the existence 
of open invariant cones in a trivial central extension  
is necessary (Remark~\ref{rem:4.7}). However, as the example 
$\cV(\bS^1)$ shows, it is not sufficient 
(Theorem~\ref{thm:sembotriv}).

Clearly, the circle $\bS^1$ plays a special role in many constructions, 
as the rich theory for central extensions of loop groups and the 
Virasoro group shows. Beyond $\bS^1$, it seems that 
Lie algebras of conformal vector fields 
(as generalizations of $\cV(\bS^1)$) (cf.\ 
Example~\ref{ex:conform}, \cite{MdR06} and \cite{Se76}) and Lie algebras of 
sections of vector bundles over Lorentzian manifolds  
(or more general ``causal'' spaces, \cite{HO96}) 
are natural candidates to be investigated 
with respect to the existence of semibounded representations. 
The latter class of Lie algebras is a natural 
generalization of loop algebras. 
Here an interesting point is that, although the conformal 
groups $\Conf(\bS^1) = \Diff(\bS^1)$ of the circle and the conformal 
group $\Conf(\bS^{1,1})$ of the Lorentzian torus 
are infinite dimensional , 
for Lorentzian manifolds of dimension $\geq 3$ the conformal groups are finite 
dimensional (cf.\ \cite{Sch97}). In particular, it is contained in the list 
of \cite{MdR06}. This leads to the well-studied class of hermitian 
Lie groups (see in particular \cite{Se76}, \cite{SJOSV78}, \cite{HO96}). 


{\bf Coadjoint orbits:} There is a symplectic version of semibounded 
representations, namely Hamiltonian actions 
$\sigma \: G \times M \to M$ with a momentum map 
$\Phi \: M \to \g'$ for which the image of the momentum map 
is semi-equicontinuous. If $G$ is finite dimensional 
and $\Phi(M)$ is closed, then Proposition~\ref{prop:locomp} 
implies that there exists an $x \in \g$ for which the 
Hamiltonian function $H_x(m) := \Phi(m)(x)$ is proper. 
In particular, $\Phi$ is a proper map. In this sense 
the semi-equicontinuity of the image of $\Phi$ is 
a weakening of properness, which is a useful property 
as far as convexity properties 
are concerned (\cite{HNP94}). 

Even though we do not know in general to which extent coadjoint 
orbits of infinite dimensional Lie groups are manifolds, the orbit 
$\cO_\lambda$ can always be viewed as the range of a momentum map 
$\Phi \: G \to \g', g \mapsto \Ad^*(g)\lambda$ corresponding to the 
left action of $G$ on itself preserving the left invariant 
closed $2$-form $\Omega$ with $\Omega_1(x,y) = \lambda([x,y])$. 
This action is semibounded if and only if $\cO_\lambda$ is 
semi-equicontinuous. 

Since the momentum sets of semibounded representations always consists  
of semi-equicontinuous orbits, the identification of the set 
$\g'_{\rm seq}$ of semi-equicontinuous coadjoint orbits of a given 
Lie algebra is already a solid first step towards the understanding 
of corresponding semibounded representations. 
As we have seen in many situations above, a useful tool to 
study convexity properties of coadjoint orbits 
are projection maps 
$p_\ft \: \g \to \ft$, where $\ft \subeq \g$ is a 
``compactly embedded'' subalgebra, i.e., 
the action of $e^{\ad \ft}$  on $\g$ factor through 
the action of a compact abelian group. 
As we have seen in Section~\ref{sec:9}, sometimes one 
does not want to project to abelian subalgebras and 
one has to study projections $p_\fk \: \g \to \fk$, 
where $\fk$ is a subalgebra for which $e^{\ad \fk}$ 
leaves a norm on $\g$ invariant 

This is of particular interest to understand open invariant 
cones $W \subeq \g$ because they often have the form 
$$ W = \Ad(G) (W \cap \fk) \quad \mbox{ with  } \quad 
p_\fk(W) = W \cap \fk $$ 
(cf.\ Theorem~\ref{thm:coneconj}). 
In this situation on needs convexity theorems of the type 
$$ p_\fk(\cO_x) \subeq \oline{\conv}(\Ad(N_G(\fk))x) + C, $$
where $C$ is a certain invariant convex cone in $\fk$ and
$N_G(\fk) \subeq G$ is the normalizer of $\fk$ in $G$ 
(cf.\ \cite{Ne00} for finite dimensional Lie algebras). 
If $\fk = \ft$ is abelian, then $N_G(\ft)$ is an analog 
of the Weyl group. 

For infinite dimensional Lie algebras, not many convexity theorems 
are known. Of relevance for semibounded representations 
is the particular case of affine Kac--Moody Lie algebras 
(\cite{AP83}, \cite{KP84}), and for Lie algebras of bounded 
operators on Hilbert spaces 
the infinite dimensional version of Kostant's Theorem by 
A.~Neumann (\cite{Neu02}) is crucial. What is still lacking 
is a uniform framework for results of this type.

\appendix
\section{Smooth vectors for representations} \mlabel{app:a}

Let $G$ be a Lie group with Lie algebra $\g$ and exponential 
function $\exp_G \: \g \to G$. 
Further, let $V$ be a locally convex space and 
$\pi \: G \to \GL(V)$ be a homomorphism 
defining a continuous action of  $G$ on $V$. We write  
$\pi^v(g) := \pi(g)v$ for the orbit maps and 
$$V^\infty := \{ v \in V \: \pi^v \in C^\infty(G,V) \} $$
for the space of smooth vectors. 
In this appendix we collect some results of 
\cite{Ne10} that are used in the present paper.
Let 
$$\dd \pi \: \g \to \End(V^\infty), \quad 
\dd\pi(x)v := \derat0 \pi(\exp tx)v $$
denote the derived  action of $\g$ on $V^\infty$. 
That this is indeed a representation of $\g$ follows by observing that 
the map $V^\infty \to C^\infty(G,V), v \mapsto \pi^v$ 
intertwines the action of $G$ with the right translation action on 
$C^\infty(G,V)$, and this implies that the derived action 
corresponds to the action of $\g$ on $C^\infty(G,V)$ by left invariant
 vector fields (cf.\ \cite[Rem.~IV.2]{Ne01b} for details). 

\begin{defn} Let $G$ be a Banach--Lie group and 
write $\cP(V)$ for the set of continuous seminorms on $V$. 
For each $p \in \cP(V)$ and $n \in \N_0$ we define a seminorm 
$p_n$ on $V^\infty$ by 
$$ p_n(v) := \sup \{ p(\dd\pi(x_1) \cdots \dd\pi(x_n)v) \: 
\|x_i\|\leq  1\} $$
and endow $V^\infty$ with the locally convex topology defined by these 
seminorms. 
\end{defn}

\begin{thm} \mlabel{thm:4.1.3} 
If $(\pi, V)$ is a representation of the Banach--Lie group 
$G$ on the locally convex space $V$ defining a continuous action 
of $G$ on $V$, then the 
action $\sigma(g,v) := \pi(g)v$ of $G$ on $V^\infty$ is smooth. 
If $V$ is a Banach space, then $V^\infty$ is complete, 
i.e., a Fr\'echet space. 
\end{thm}

\begin{thm} \mlabel{thm:unitautsmooth} 
If $(\pi, \cH)$ is a unitary representation of a Lie group 
$G$, then $v \in \cH$ is a smooth vector if and only if the 
corresponding matrix coefficient $\pi^{v,v}(g) := \la \pi(g)v,v\ra$ 
is smooth on a $\1$-neighborhood in $G$. If, in addition, 
$v$ is cyclic, then the representation is smooth. 
\end{thm}

In the following we write $\oline u$ for the image of a unitary 
operator $u \in \U(\cH)$ in the projective unitary group 
$\PU(\cH) :=\U(\cH)/\T \1$. 

\begin{thm} \mlabel{thm:unirep-cenext} Let $G$ be a connected Lie group, 
$\cH$ a complex Hilbert space and  
$\pi \: G \to \U(\cH)$ be a map with $\pi(\1) = \1$ 
for which the corresponding map 
$\oline\pi \: G \to \PU(\cH)$ is a group homomorphism. 

If there exists a $v \in \cH$ for which 
the function $(g,h) \mapsto \la \pi(g)\pi(h)v, v\ra$ is 
smooth on a neighborhood of $(\1,\1)$ in $G \times G$, 
then the central 
extension 
$$\hat G := \oline\pi^*\U(\cH) = \{ (g,u) \in G \times \U(\cH)  
\: \oline\pi(g) = \oline u \} $$
of $G$ by $\T$ is a Lie group and 
$v$ is a smooth vector for the 
representation $(\hat\pi, \cH)$ of $\hat G$ by 
$\hat\pi(g,u) := u$. 
\end{thm}

\begin{prf} To exhibit $\hat G$ as a Lie group, we have to show that 
there exists a section $\sigma \: G \to \hat G$ for which the 
corresponding $2$-cocycle 
$$ f_\sigma(g_1, g_2) 
= \sigma(g_1) \sigma(g_2) \sigma(g_1 g_2)^{-1} $$ 
is smooth in a neighborhood of 
$(\1,\1)$ (\cite[Prop.~4.2]{Ne02b}). 
Here we use that $G$ is connected. 
A particular section $\sigma \: G \to \hat G$ is given by 
$\sigma(g) = (g,\pi(g))$. 

Let $U \subeq G$ be an open $\1$-neighborhood such that 
$$ \la \pi(g)v,\pi(h)^{-1}v \ra \not=0 \quad \mbox{ for } \quad g,h \in U. $$
Its existence follows from our continuity assumption. 
If $U'\subeq U$ is an open $\1$-neighborhood with $U'U' \subeq U$, we then 
have for $g_1, g_2 \in U'$ 
$$ f_\sigma(g_1, g_2) \la \pi(g_1 g_2)v, v\ra 
= \la \pi(g_1) \pi(g_2) v, v\ra, $$
which leads to 
$$ f_\sigma(g_1, g_2) = \frac{\la \pi(g_1)\pi(g_2) v, v\ra}
{\la \pi(g_1 g_2)v, v\ra}. $$
Therefore $f_\sigma$ is smooth in a neighborhood of $(\1,\1)$.
This shows that $\hat G$ is a Lie group and the multiplication map 
$G \times \T \to \hat G, (g,t) \mapsto (g, t\pi(g))$
is smooth in a neighborhood of $\1$ (\cite[Prop.~4.2]{Ne02b}). 
The representation 
$\hat\pi$ now satisfies 
$$ \la \hat\pi(g,t\pi(g))v, v  \ra   
=  \la t\pi(g)v, v  \ra = t \la \pi(g)v,v\ra, $$
which is smooth in a neighborhood of $\1$. Now 
Theorem~\ref{thm:unitautsmooth} implies that $v$ is a smooth 
vector for the representation $\hat\pi$.
\end{prf}

\begin{rem} The assumption that $G$ is connected in the preceding theorem 
can be removed if $G$ is a Banach--Lie group. In this case we assume, 
in addition, that $\oline\pi$ is continuous in the sense that 
all functions $g \mapsto |\la \pi(g)v,w \ra|$ on $G$ are continuous 
(cf.\ \cite[p.~175]{Mag92}). 
Then $\hat G$ is the pullback of a central extension 
of topological groups $\U(\cH) \to \PU(\cH)$, hence in particular 
a topological group with respect to the topology inherited from 
$G \times \U(\cH)$. Under the assumptions of Theorem~\ref{thm:unirep-cenext}, 
the central extension $\hat G_0$ of the identity component 
$G_0$ of $G$ carries a natural Lie group structure compatible 
with the given topology on $\hat G$. For each $g \in \hat G$, the conjugation 
map $c_g$ induces a continuous automorphism of the Banach--Lie group 
$\hat G_0$, and since continuous homomorphisms of Banach--Lie groups 
are automatically smooth (\cite[Thm.~IV.1.18]{Ne06}), 
$c_g$ also defines a Lie automorphism of 
$\hat G_0$. This implies that $\hat G$ carries a unique Lie group 
structure which coincides on the open subgroup $\hat G_0$ with the given 
one (see also \cite[Rem.~4.3]{Ne02b}). 
\end{rem}

\begin{rem} \mlabel{rem:a.11} 
In the situation of Theorem~\ref{thm:unirep-cenext}, 
the Lie algebra cocycle defining the central extension 
$\hat\g = \L(\hat G)$ of $\g$ by $\R$ in the sense that 
$$ \hat\g = \R \oplus_{\eta} \g, \quad 
[(t,x), (t',x')] = (\eta(x,x'), [x,x']) $$
can be calculated as follows. 
\glossary{name={$\R \oplus_{\eta} \g$},description={central extension defined by 
cocycle $\eta$}}

If $\dd\pi \: \g \to \End(\cH^\infty)$ is the map obtained 
from the representation $\dd\hat\pi$ via 
$\dd\pi(x) := \dd\hat\pi(0,x)$, then we have for each unit 
vector $v \in \cH^\infty$ the relation 
$$ [\dd\pi(x), \dd\pi(y)]v 
= \dd\pi([x,y])v + i \eta(x,y)v, $$
so that 
\begin{align*}
 \eta(x,y) 
&= \Im \la [\dd\pi(x), \dd\pi(y)]v, v\ra 
- \Im \la \dd\pi([x,y])v, v \ra \\
&= 2\Im \la \dd\pi(x)v, \dd\pi(y) v\ra + i \la \dd\pi([x,y])v, v \ra.
\end{align*}
\end{rem}

\section{The cone of positive definite forms on a Banach space} 
\mlabel{app:b}

Let $V$ be a Banach space and 
$\Sym^2(V,\R)$ be the Banach space of continuous symmetric bilinear 
maps $\beta \: V \times V \to \R$, endowed with the norm 
$$\|\beta\| := \sup \{ |\beta(v,v)| \: \|v\| \leq 1\}. $$
Clearly, the set $\Sym^2(V,\R)_+$ of positive semidefinite 
bilinear maps is a closed convex cone in $\Sym^2(V,\R)$. 

\begin{lem}
  \mlabel{lem:b.1} 
The cone $\Sym^2(V,\R)_+$ has interior points if and only if 
$V$ is  topologically isomorphic to a Hilbert space. If this is 
the case, then its interior consists of all those positive definite 
forms $\beta$ for which the norm $\|v\|_\beta := \sqrt{\beta(v,v)}$ 
is equivalent to the norm on $V$. 
\end{lem}

\begin{prf} Suppose first that $V$ carries a Hilbert space 
structure $\beta$. Replacing the original norm by an equivalent 
Hilbert norm, we may assume that $V$ is a real Hilbert space. 
Then $\Sym^2(V,\R)$ can be identified with the space 
$\Sym(V)$ of symmetric operators on $V$ by assigning to 
$A \in \Sym(V)$ the form $\beta_A(v,v) = \la Av,w\ra$, satisfying 
$$ |\beta_A(v,v)| = |(Av,v)| \leq \|A\|\|v\|^2, 
\quad \mbox{ so that } \quad \|\beta_A\| \leq \|A\|.$$ 
The polarization identity 
$$ (Av,w) = \beta_A(v,w) = \frac{1}{4} (\beta_A(v+w,v+w) - 
\beta_A(v-w,v-w)) $$ 
further implies that $\|A\| \leq 2 \beta_A$, so that the 
Banach spaces $\Sym(V)$ and \break 
$\Sym^2(V,\R)$ are topologically isomorphic. 
The identity $\id_V = \1$ is an interior point in the cone 
of positive operators, $\Sym^2(V,\R)_+$ has interior points. 

Suppose, conversely, that $\beta \in \Sym^2(V,\R)_+$ is an 
interior point. Then \break $\beta - \Sym^2(V,\R)_+$ is a $0$-neighborhood 
in $\Sym^2(V,\R)$, which implies the existence of some 
$c > 0$ such that $\|\gamma\| \leq c$ implies 
$\gamma(v,v) \leq \beta(v,v)$ for all $v \in V$. 
Fixing $v_0 \in V$, we pick $\alpha\in V'$ with 
$\|\alpha\| = 1$ and $\alpha(v_0) = \|v_0\|$. 
Then $\gamma(v,w) := c\alpha(v)\alpha(w)$ is a symmetric bilinear form 
with $\|\gamma\| = c\|\alpha\|^2 = c$. We therefore obtain 
$$ c \|v_0\|^2 = \gamma(v_0, v_0) \leq \beta(v_0, v_0) 
\leq \|\beta\| \|v_0\|^2, $$ 
showing that $\beta$ is positive definite and that the norm 
$\sqrt{\beta(v,v)}$ is equivalent to the norm on $V$. 
\end{prf}

\begin{rem} \mlabel{rem:b.2}
Let $(V,\omega)$ be a weakly symplectic Banach space, i.e., 
$\omega$ is a non-degenerate skew-symmetric bilinear form. 
Then 
$$ \sp(V,\omega) := \{ X \in B(V) \: (\forall v,w\in V)\, 
\omega(Xv,w) + \omega(v,Xw) = 0\} $$
is the corresponding symplectic Lie algebra. In particular, 
$X \in B(V)$ belongs to $\sp(V,\omega)$ if and only if 
the bilinear form $\omega(Xv,w)$ is symmetric. The corresponding 
quadratic function $H_X(v) := \shalf\omega(Xv,v)$ is called the 
{\it Hamiltonian function defined by $X$}. 

If there exists an $X \in \sp(V,\omega)$ for which 
$(v,w) := \omega(Xv,w)$ defines a Hilbert space structure 
on $V$, then each continuous linear 
functional $\alpha \in V'$ is of the form $i_v\omega$ for 
some $v \in V$. This means that 
$(V,\omega)$ is {\it strongly symplectic}, 
i.e., the map 
$\Phi_\omega \: V \to V', v \mapsto i_v\omega$ is surjective, 
hence a topological isomorphism by the Open Mapping Theorem. 
We further see that $X$ is injective. To see that it is also surjective, 
let $u \in V$ and represent the continuous linear functional 
$i_u\omega$ as $(w,\cdot)$ for some $w \in V$. Then 
$\omega(u,v) = (w,v) = \omega(Xw,v)$ for all $v \in V$, 
and thus $Xw = u$. Hence $X$ is a topological isomorphism, 
satisfying $\omega(u,v) = (X^{-1}u,v)$ for $v,w \in V$. This further implies 
that $X$ is skew-symmetric.  
\end{rem} 

\begin{prop} \mlabel{prop:b.3} 
Let $(V,\omega)$ be a strongly symplectic Banach space. 
Then the convex cone $\{ X \in \sp(V,\omega) \: H_X \geq 0\}$ has interior 
points if and only if $V$ is topologically isomorphic to a 
Hilbert space. If this is the case, 
then there exists a complex structure $I \in \sp(V,\omega)$ for 
which $H_I$ defines a compatible Hilbert space structure on $V$, i.e., 
$$ \omega(v,w) = \Im \la v, w\ra $$
for the underlying complex Hilbert space. 
\end{prop} 

\begin{prf} (cf.\ \cite[Thm.~3.1.19]{AM78}) 
If $\beta \: V \times V \to \R$ is a symmetric bilinear form, 
then $\Phi_\beta \: V \to V', v \mapsto i_v\beta$ also is a continuous 
linear map, and 
$B := \Phi_\omega^{-1} \circ \Phi_\beta \: V \to V$ is continuous 
linear with 
$$ \omega(Bv,v) = \Phi_\omega(Bv)(v) 
= \Phi_\beta(v)(v) = \beta(v,v), $$
so that every symmetric bilinear form can be represented by an 
element of $\sp(V,\omega)$, and we obtain a topological isomorphism 
$\sp(V,\omega) \cong \Sym^2(V,\R)$. Combining this with 
Lemma~\ref{lem:b.1} proves the first assertion. 

We now assume that $V$ is a real Hilbert space with the scalar product 
$(v,w) = \omega(Xv,w)$, where $X$ is as above. Then $A := X^{-1}$ is an 
invertible skew-symmetric operator with $\omega(v,w) = (Av,w)$. 
Then the complex linear extension $A_\C$ to the complex Hilbert 
space $V_\C$ yields a self-adjoint operator $iA_\C$. 
The complex conjugation $\sigma$ of $V_\C$ with respect to $V$ 
now satisfies $\sigma \circ iA_\C \circ \sigma = - i A_\C$, so 
that $\Spec(iA_\C)$ is a compact symmetric subset of $\R$, not 
containing $0$. We therefore have an orthogonal decomposition 
$V_\C = V_+ \oplus V_-$ into the positive and negative spectral 
subspaces of $iA_\C$. Since $\sigma(V_\pm) = V_\mp$, we obtain 
an isomorphism $V \cong V_+$, and hence a complex structure 
$I$ on $V$, corresponding to multiplication by $i$ on $V_+$. 
This means that $v_+ \in V_+$ with $v = v_+ + \sigma(v_+)$ 
satisfies 
\begin{align*}
\omega(Iv,v) 
&= (AIv,v) 
= (AI(v_+ + \sigma(v_+)), v_+ + \sigma(v_+)) \\
&= (A_\C(i v_+ - i \sigma(v_+)), v_+ + \sigma(v_+)) 
= i (A_\C(v_+ - \sigma(v_+)), v_+ + \sigma(v_+)) \\
&= i \big((A_\C v_+, v_+) - (v_+, A_\C v_+)\big) 
= 2 (iA_\C v_+, v_+), 
\end{align*}
so we obtain a complex structure $I$ on $V$ 
for which the Hamiltonian $H_I$ defines a Hilbert structure on 
$V$. Since $I$ is skew-symmetric, $(V,I)$ inherits the structure 
of a complex Hilbert space with respect to the scalar product 
$$ \la v, w \ra 
:= \omega(Iv,w) + i \omega(Iv,Iw) = \omega(Iv,w) + i \omega(v,w). $$
Therefore our assumption leads to the representation 
of $\omega$ as $\omega(v,w) = \Im \la v, w \ra$ for a complex 
Hilbert space structure on $V$. 
\end{prf}

\section{Involutive Lie algebras with root decomposition} 
\mlabel{app:e}

\begin{defn} \label{def:basic} (a) We call an abelian subalgebra 
$\ft$ of the real Lie algebra $\g$ a {\it compactly embedded 
Cartan subalgebra} if $\ft$ is maximal abelian and $\ad \ft$ 
is simultaneously diagonalizable on the complexification $\g_\C$ 
with purely imaginary eigenvalues. 
Then we have a {\it root decomposition} 
$$ \g_\C = \ft_\C + \sum_{\alpha \in\Delta} \g_\C^\alpha, $$
where $\g_\C^\alpha = \{ x \in \g_\C \: (\forall h \in \ft_\C) [h,x]= \alpha(h)
x\}$ and 
$$\Delta :=  \{ \alpha \in \ft_\C^* \setminus \{0\} \: \g_\C^\alpha
\not= \{0\}\}$$ is the corresponding {\it root system}. 

If $\sigma \: \g_\C \to \g_\C$ denotes the complex conjugation with 
respect to $\g$, we write $x^* := -\sigma(x)$ for $x \in \g_\C$, so that 
$\g = \{ x \in \g_\C \: x^* = -x\}$. 
We then have 
\begin{itemize}
\item[\rm(I1)] $\alpha(x) \in \R$ for $x \in i \ft$. 
\item[\rm(I2)] $\sigma(\g_\C^\alpha) = \g_\C^{-\alpha}$ for $\alpha \in \Delta.$ 
\end{itemize}
\end{defn}

\begin{lem}
  \mlabel{lem:e.1} For $0 \not= x_{\alpha} \in \g_\C^{\alpha}$ the subalgebra 
$\g_\C(x_\alpha) := \Spann\{x_\alpha, x_{\alpha}^*,
[x_\alpha, x_{\alpha}^*]\}$
is $\sigma$-invariant and of one of the following types: 
\begin{description}
\item[\rm(A)] The abelian type: $[x_\alpha, x_{\alpha}^*] = 0$, i.e., 
$\g_\C(x_\alpha)$ is two dimensional abelian. 
\item[\rm(N)] The nilpotent type: $[x_\alpha, x_{\alpha}^*] \not= 0$ 
and $\alpha([x_\alpha, x_{\alpha}^*]) = 0$, i.e., 
$\g_\C(x_\alpha)$ is a three dimensional Heisenberg algebra. 
\item[\rm(S)] The simple type: $\alpha([x_\alpha, x_{\alpha}^*]) \not= 0$,
i.e., $\g_\C(x_\alpha) \cong \fsl_2(\C)$. In this case we distinguish the
two cases: 
\item[\rm(CS)] $\alpha([x_\alpha, x_{\alpha}^*]) > 0$, i.e., 
$\g_\C(x_\alpha) \cap \g \cong \su_2(\C)$, and 
\item[\rm(NS)] $\alpha([x_\alpha, x_{\alpha}^*]) < 0$, i.e., 
$\g_\C(x_\alpha) \cap \g \cong \su_{1,1}(\C) \cong \fsl_2(\R)$. 
\end{description}
\end{lem} 

\begin{prf}
First we note that, in view of $x_\alpha^* \in \g_\C^{-\alpha}$,  
\cite[Lemma~I.2]{Ne98} applies, and we see that $\g_\C(x_\alpha)$ is of one of
the three types (A), (N) and (S). We note that 
$\alpha([x_\alpha, x_\alpha^*]) \in \R$ 
because of (I2) and $[x_\alpha, x_\alpha^*] \in i\ft$. 
Now it is easy to check that 
$\g_\C(x_\alpha) \cap \g$ is of type (CS), resp., (NS), according to the sign
of this number. 
\end{prf} 

The following proposition provides useful information for the 
analysis of invariant cones and orbit projections. 
Here we write $p_\ft \: \g \to \ft$ for the 
projection along $[\ft,\g] = \sum_{\alpha} (\g_\C^\alpha + \g_\C^{-\alpha}) \cap \g$, 
and $p_{\ft^*} \: \g^* \to \ft^*$ for the restriction map.

\begin{prop}
  \mlabel{prop:orb-pro} 
Let $x \in \ft$, $x_\alpha \in \g_\C^\alpha$ and 
$\lambda \in \ft^*$. Then the following assertions hold: 
\begin{description}
\item[\rm(i)] 
$p_\ft(e^{\R\ad (x_\alpha - x_\alpha^*)}x) 
= x  + 
\begin{cases}
\R^+ \alpha(x) [x_\alpha^*, x_\alpha]  & \text{ for }\  
\alpha([x_\alpha, x_\alpha^*]) \leq 0 \\ 
[0,2] { \frac{\alpha(x)}{\alpha([x_\alpha, x_\alpha^*])}} 
[x_\alpha^*, x_\alpha] & \text{ for }\ \alpha([x_\alpha, x_\alpha^*]) > 0. 
\end{cases}$
\item[\rm(ii)] 
$p_{\ft^*}(e^{\R\ad^* (x_\alpha - x_\alpha^*)}\lambda) 
= \lambda  + 
\begin{cases}
\R^+ \lambda([x_\alpha^*, x_\alpha]) \alpha   & \text{ for }\  
\alpha([x_\alpha, x_\alpha^*]) \leq 0 \\ 
[0,2] \frac{\lambda([x_\alpha^*, x_\alpha])}{\alpha([x_\alpha, x_\alpha^*])} 
\alpha & \text{ for }\ \alpha([x_\alpha, x_\alpha^*]) > 0. 
\end{cases}$
\end{description}
\end{prop}

\begin{prf} (i) is an immediate consequence of \cite[Lemma~VII.2.9]{Ne00}, 
and (ii) follows from (i) and the relation 
$p_{\ft^*}(\lambda)(e^{\ad y}x)
= \lambda(p_\ft(e^{\ad y}x)).$
\end{prf}

\section{Some facts on Fock spaces} \mlabel{app:c}

Let $\cH$ be a complex Hilbert space. We endow the $n$-fold tensor 
product with the canonical Hilbert structure defined by 
$$ \la v_1 \otimes \cdots \otimes v_n, w_1 \otimes \cdots \otimes w_n \ra 
:= \prod_{j = 1}^n \la v_j, w_j \ra, $$
and form the Hilbert space direct sum 
$T(\cH) := \hat\bigoplus_{n \in \N_0} \cH^{\otimes n}$. 
In $\cH^{\otimes n}$ we write $S^n(\cH)$ for the closed subspace 
generated by the symmetric tensors and 
$\Lambda^n(\cH)$ for the closed subspace generated by the 
alternating tensors. We thus obtain subspaces 
$$ S(\cH) := \hat\bigoplus_{n \in \N_0} S^n(\cH) 
\quad \mbox{ and } \quad 
 \Lambda(\cH) := \hat\bigoplus_{n \in \N_0} \Lambda^n(\cH) $$
of $T(\cH)$ and write $P_s$, resp., $P_a$ for the corresponding 
orthogonal projections: 
$$ P_s(f_1 \otimes \cdots \otimes f_n) 
= \frac{1}{n!} \sum_{\sigma \in S_n} f_{\sigma(1)} \otimes \cdots \otimes 
f_{\sigma(n)} $$  
and 
$$ P_a(f_1 \otimes \cdots \otimes f_n) 
= \frac{1}{n!} \sum_{\sigma \in S_n} \sgn(\sigma) 
f_{\sigma(1)} \otimes \cdots \otimes f_{\sigma(n)}. $$  

The dense subspace $S(\cH)_0 := \sum_{n \geq 0} S^n(\cH)$ of $S(\cH)$ 
carries an associative 
algebra structure, given by 
$$ f_1 \vee \cdots \vee f_n := \sqrt{n!}P_s(f_1 \otimes \cdots \otimes f_n) $$
and likewise $\Lambda(\cH)_0 := \sum_{n \geq 0} \Lambda^n(\cH)$ 
inherits an algebra structure defined by 
$$ f_1 \wedge  \cdots \wedge f_n := \sqrt{n!}
P_a(f_1 \otimes \cdots \otimes f_n). $$
A unit vector $\Omega$ in the one-dimensional space 
$S^0(\cH) = \Lambda^0(\cH)$ is called a {\it vacuum vector}. 

\begin{lem}
  \mlabel{lem:a.c.1} We have
\begin{equation}
  \label{eq:symprod}
 T \vee S 
= \sqrt{{n+m \choose n}} P_s(T \otimes S) \quad  
\mbox{ for } \quad T \in S^n(\cH), S \in S^m(\cH) 
\end{equation}
and 
\begin{equation}
  \label{eq:skewprod}
 T \wedge S 
= \sqrt{{n+m \choose n}} P_a(T \otimes S) \quad  
\mbox{ for } \quad T \in \Lambda^n(\cH), S \in \Lambda^m(\cH). 
\end{equation}
\end{lem} 

\begin{prf} For the symmetric case we first note that 
$f^n = \sqrt{n!} f^{\otimes n}$, so that 
$$ f^n \vee g^m 
= \sqrt{(n+m)!} P_s(f^{\otimes n} \otimes g^{\otimes m}) 
= \sqrt{{n+m \choose n}} P_s(f^n \otimes g^m). $$
Since the elements $f^n$ generate $S^n(\cH)$ topologically, 
\eqref{eq:symprod} follows. 

For the alternating case we obtain for 
$T = f_1 \wedge \cdots \wedge f_n$ and 
$S = g_1 \wedge \cdots \wedge g_m$ the relation 
\begin{align*}
&\ \ \ \ \frac{1}{\sqrt{(n+m)!}}T \wedge S 
= P_a(f_1 \otimes \cdots \otimes g_m)\\
&= P_a(
P_a(f_1 \otimes \cdots \otimes f_n) \otimes 
P_a(g_1 \otimes \cdots \otimes g_m)) 
= \frac{1}{\sqrt{n!m!}} P_a(T \otimes S). 
\end{align*}
\end{prf}

\begin{rem} \mlabel{rem:fock-alg} 
(a) For the norms of the product of $T \in S^n(\cH)$ and $S \in S^m(\cH)$, 
we obtain with \eqref{eq:symprod} 
$$ \frac{1}{\sqrt{(n+m)!}} \|T \vee S\| 
= \frac{1}{\sqrt{n!}\sqrt{m!}} \|P_s(T \otimes S)\|
\leq \frac{1}{\sqrt{n!}\sqrt{m!}}  \|T \otimes S\| 
=  \frac{\|T\|}{\sqrt{n!}} \frac{\|S\|}{\sqrt{m!}}.$$
From that we derive for $T \in S^k(\cH)$ and $n \in \N$ the relation
$$ \frac{1}{\sqrt{(kn)!}} \|T^n\| 
\leq \Big(\frac{1}{\sqrt{k!}}\|T\|\Big)^n, $$
and for $T \in S^2(\cH)$ we find in particular 
\begin{equation}
  \label{eq:expesti}
\|e^T\|^2 
= \sum_{n = 0}^\infty \frac{1}{(n!)^2} \|T^n\|^2 
\leq \sum_{n = 0}^\infty \frac{(2n)!}{(n!)^2 2^n} \|T\|^{2n}, 
\end{equation}
which converges for $\|T\|^2 < \frac{1}{2}$. 

(b) For $T \in \Lambda^n(\cH)$ and $S \in \Lambda^m(\cH)$ 
we likewise obtain with \eqref{eq:skewprod} 
$$ \frac{1}{\sqrt{(n+m)!}} \|T \wedge S\| 
\leq  \frac{\|T\|}{\sqrt{n!}} \frac{\|S\|}{\sqrt{m!}},$$ 
which leads for $T \in \Lambda^2(\cH)$ to 
\begin{equation}
  \label{eq:expestialt}
\|e^T\|^2 
= \sum_{n = 0}^\infty \frac{1}{(n!)^2} \|T^n\|^2 
\leq \sum_{n = 0}^\infty \frac{(2n)!}{(n!)^2 2^n} \|T\|^{2n}. 
\end{equation}
\end{rem}

\begin{lem} \mlabel{lem:hs2} 
Let $A$ be an antilinear Hilbert--Schmidt operator on $\cH$ 
and define $A^*$ by $\la A^*v,w \ra = \la Aw,v\ra$. 
If $A^* = A$, then there exists a unique element 
$\hat A \in S^2(\cH)$ with
$$ \la \hat A, f_1 \vee f_2 \ra = \la Af_1, f_2 \ra f
\quad \mbox{ for } \quad 
f_1, f_2 \in \cH,  $$ 
and if $A^* = - A$, there exists a unique element $\hat A \in \Lambda^2(\cH)$ 
with 
$$ \la \hat A, f_1 \wedge  f_2 \ra = \la Af_1, f_2 \ra 
\quad \mbox{ for } \quad 
f_1, f_2 \in \cH. $$ 
Moreover, the following assertions hold: 
\begin{description}
\item[\rm(i)] $\|\hat A\|^2 = \shalf \|A\|_2^2$. 
\item[\rm(ii)] If $A$ and $B$ are hermitian and antilinear, then  
$\la \hat A, \hat B \ra = \shalf \tr(AB)$. 
\item[\rm(iii)] If $A$ and $B$ are skew-hermitian and antilinear, then  
$\la \hat A, \hat B \ra = -\shalf \tr(AB)$. 
\end{description}
\end{lem}

\begin{prf} (a)  First we consider the case where $A^* = A$. 
Let $(e_j)_{j \in J}$ be an orthonormal basis of $\cH$. 
Then we have in $S^2(\cH)$ the relations 
$$ \|e_j^2\|^2 = 2\| e_j \otimes e_j\|^2 = 2 
\quad \mbox{ and } \quad 
\|e_j \vee e_k\|^2 
= 2\|\shalf(e_j \otimes e_k + e_k \otimes e_j)\|^2 = 1, \quad 
j\not=k. $$
If $<$ denotes a linear order on $J$, we thus obtain the 
orthonormal basis 
$$\frac{1}{\sqrt 2} e_j^2, \quad e_j \vee e_k, \quad j < k, \quad 
\mbox{ of } \quad S^2(\cH).$$ 
This leads to 
\begin{align*}
\|\hat A\|^2 
&= \sum_j \shalf |\la \hat A, e_j^2 \ra|^2 
+ \sum_{j<k} |\la \hat A, e_j \vee e_k \ra|^2 \\
&= \sum_j \shalf |\la Ae_j, e_j \ra|^2 
+ \sum_{j<k} |\la Ae_j, e_k \ra|^2 
= \shalf \|A\|_2^2. 
\end{align*}

For $A^* = -A$ we similarly get 
\begin{align*}
\|\hat A\|^2 
&= \sum_{j<k} |\la \hat A, e_j \vee e_k \ra|^2
= \sum_{j<k} |\la Ae_j, e_k \ra|^2 = \shalf \|A\|_2^2. 
\end{align*}

(b) Let $\Herm_2(\cH)_a$ denote the complex subspace 
of hermitian antilinear 
Hilbert--Schmidt operators. 
Then the prescription $\la A, B \ra := \tr(AB) = \tr(AB^*)$ defines a 
sesquilinear form on this space with 
$$ \oline{\tr(AB)} = \tr((AB)^*) = \tr(B^*A^*) = \tr(BA), $$
so that it is hermitian. Its restriction to the diagonal 
satisfies 
$$\la A,A \ra = \tr(A^2)
= \sum_{j \in J} \la A^2e_j, e_j \ra 
= \sum_{j \in J} \la Ae_j, Ae_j \ra = \|A\|_2^2. $$
In view of (a), polarization implies that 
$\tr(AB) = \la A, B \ra = 2 \la \hat A, \hat B\ra$ 
for $A,  B \in \Herm_2(\cH)_a$. 

(c) For the space $\Aherm_2(\cH)_a$ of skew-hermitian antilinear 
Hilbert--Schmidt operators the same argument works with 
$\la A, B \ra := -\tr(AB) = \tr(AB^*)$. 
\end{prf}


\begin{thebibliography}{aaaaaaaa} 

\bibitem[AM78]{AM78} Abraham, R., and J. E. Marsden, ``Foundations of 
Mechanics,'' 2nd Ed., Benjamin/Cummings Publ. Comp., 1978 

\bibitem[AW64]{AW64} Araki, H., and W. Wyss, {\it Representations 
of the canonical anticommutation relations}, Helv. Phys. Acta 
{\bf 37} (1964), 136--159 

\bibitem[AP83]{AP83} Atiyah, M. F., and A. N. Pressley, 
{\it Convexity and loop groups}, in ``Arithmetic and Geometry, Vol. II'', 
Birkh\"auser, Basel, 1983 


\bibitem[Bak07]{Bak07} Bakalov, B., N. M. Nikolov, 
K.-H. Rehren, I. Todorov, 
{\it Unitary positive energy representations 
of scalar bilocal quantum fields}, Comm. Math. Phys. 
{\bf 271:1} (2007), 223--246 

\bibitem[BMV68]{BMV68} Balslev, E., J. Manuceau, and A. Verbeure, 
{\it Representations of the anticommutation relations and Bogoliubov 
transformations}, Comm. Math. Phys. {\bf 8} (1968), 315--326 

\bibitem[Bel06]{Bel06} Belti\c{t}\u{a}, D.,  ``Smooth Homogeneous Structures 
in Operator Theory,''  Chapman and Hall, CRC Monographs and 
Surveys in Pure and Applied mathematics, 2006 

\bibitem[BN08]{BN08} Beltita, D., and K.-H. Neeb, {\it A non-smooth 
continuous unitary representation of a Banach--Lie group}, 
J. Lie Theory {\bf 18} (2008), 933-936 

\bibitem[BN10]{BN10} Beltita, D., and K.-H. Neeb, {\it 
Schur--Weyl Theory for $C^*$-algebras}, in preparation 

\bibitem[BoSi71]{BoSi71} Bochnak, J., Siciak, J., 
{\it Analytic functions in topological vector spaces}, 
Studia Math. \textbf{39}  (1971), 77--112

\bibitem[Bo96]{Bo96} Borchers, H. -J., ``Translation group and 
particle representations in quantum field theory,'' 
Lecture Notes in Physics, Springer, 1996 

\bibitem[Bou89]{Bou89} Bourbaki, N., ``General Topology. Chaps.~1-4'', 
Springer-Verlag, Berlin, 1989 

\bibitem[Bou07]{Bou07} Bourbaki, N., ``Espaces vectoriels topologiques. 
Chap.1 \`a 5'', Springer-Verlag, Berlin, 2007 

\bibitem[BR87]{BR87} Bowick, M. J., and S. G. Rajeev, 
{\it The holomorphic geometry of the closed bosonic string 
theory and $\Diff(\bS^1)/\bS^1$}, Nuclear Physics {\bf B293} (1987), 
348--384 

\bibitem[BR97]{BR97} Bratteli, O., Robinson, D.~W., 
``Operator Algebras and Quantum Statistical Mechanics 2,'' 
2nd Ed., Texts and Monographs in Physics, Springer--Verlag, New York 1997 

\bibitem[Ca83]{Ca83} Carey, A. L., {\it Infinite Dimensional Groups and 
Quantum Field Theory}, Act. Appl. Math. {\bf 1} (1983), 321--333 

\bibitem[CHa96]{CHa96} Carey, A., and K. C. Hannabus, 
{\it Infinite dimensional groups and Riemann surface field theory}, 
Comm. Math. Phys. {\bf 176} (1996), 321--351 

\bibitem[CL02]{CL02} Carey, A., and E. Langmann, 
{\it Loop Groups and Quantum Fields}, in 
``Geometric Analysis and Applications to Quantum Field Theory,'' 
Eds. P. Bouwknegt, S. Wu, Progr. in Math. {\bf 205}, Birkh\"auser 
Verlag, 2002; 45--94 

\bibitem[CR87]{CR87} Carey, A., and S. 
N. M. Ruijsenaars, {\it On fermion gauge
groups, current algebras, and Kac--Moody algebras}, Acta
Appl. Math. {\bf 10} (1987), 1--86 


\bibitem[DF77]{DF77} Delbourgo, R., and J. R. Fox, {\it 
Maximum weight vectors possess minimal uncertainty}, 
J. Phys. A: Math. Gen. {\bf 10:12} (1977), 233--235  

\bibitem[Dix64]{Dix64} Dixmier, J., ``Les $C^*$-alg\`ebres et leurs
repr\'esentations,'' Gauthier-Villars, Paris, 1964 

\bibitem[FH05]{FH05} Fewster, Chr., and S. Hollands, {\it 
Quantum energy inequalities in two-dimensional conformal field theory}, 
Rev. Math. Phys.  {\bf 17:5} (2005), 577--612 

\bibitem[GV64]{GV64} Gel'fand, I. M., and N. Ya Vilenkin, ``Generalized 
Functions. Vol. 4: Applications of Harmonic Analysis,'' 
Translated by Amiel Feinstein, Academic Press, New York, London, 1964

\bibitem[Ge07]{Ge07} Georgescu, V., {\it On the spectral analysis 
of quantum field Hamiltonians}, arXiv:math-ph/0604072v1, 28. April 2006 

\bibitem[Gl02]{Gl02} Gl\"ockner, H., 
{\it Algebras whose groups of units are Lie groups}, 
Studia Math. {\bf 153} (2002), 147--177 

\bibitem[GN09]{GN09} Gl\"ockner, H., and K.-H. Neeb, ``Infinite dimensional 
Lie groups, Vol. I, Basic Theory and Main Examples,'' book in preparation 

\bibitem[GO86]{GO86} Goddard, P., and D. Olive, 
{\it Kac--Moody and Virasoro algebras in relation 
to quantum physics}, Internat. J. Mod. Phys. {\bf A1} (1986), 303--414 

\bibitem[GW84]{GW84} Goodman, R., and N. R. Wallach, {\it Structure and unitary
cocycle representations of loop groups and the group of
diffeomorphisms of the circle}, J. reine ang. Math. {\bf 347} (1984), 69--133

\bibitem[GW85]{GW85} ---, {\it Projective 
unitary positive energy representations of
$\Diff(\bS^1)$}, J. Funct. Anal. {\bf 63} (1985), 299-312 

\bibitem[GrNe09]{GrNe09} Grundling, H., and K.-H. Neeb, {\it 
Full regularity for a 
$C^*$-algebra of the canonical commutation relations}, 
in Reviews in Math. Physics, to appear 

\bibitem[GR07]{GR07}
Guieu, L., and Roger, C., ``L'Alg\`ebre et le Groupe de Virasoro: 
aspects g\'eom\'etriques et alg\'ebriques, g\'en\'eralisations'', 
Les Publications CRM, 2007 

\bibitem[HK64]{HK64} Haag, R., and D. Kastler, {\it 
An algebraic approach to quantum field theory}, J. Math. Phys. 
{\bf 5} (1964), 848--861

\bibitem[Ha82]{Ha82} Hannabus, K. C., {\it On a property of highest 
weight vectors}, Quart. J. Math. Oxford (2) {\bf 33} (1982), 91--96 

\bibitem[dH72]{dH72} de la Harpe, P., 
{\it The Clifford Algebra and the Spinor Group of a Hilbert Space}, 
Compos. Math. {\bf 25:3} (1972), 245--261 

\bibitem[He71]{He71} Hegerfeldt, G. C., {\it G\aa{}rding domains 
and analytic vectors for quantum fields}, J. Math. Phys. {\bf 13:6} (1972), 
821--827 

\bibitem[HiHo85]{HiHo85} Hilgert, J., and K. H. Hofmann, {\it Lorentzian cones 
in real Lie algebras}, Monatshefte f\"ur Math. {\bf 100} (1985), 183--210 

\bibitem[HHL89]{HHL89} Hilgert, J., 
K.H. Hofmann, and J.D. Lawson, ``Lie Groups, Convex 
Cones, and Semigroups,'' Oxford University Press, 1989 

\bibitem[HNP94]{HNP94} Hilgert, J., K.-H. Neeb, 
and W. Plank, {\it Symplectic convexity 
theorems and coadjoint orbits}, Comp. Math.\ {\bf 94} (1994), 129--180 

\bibitem[HO96]{HO96} Hilgert, J., and G. \'Olafsson, 
``Causal Symmetric Spaces, Geometry and Harmonic Analysis,'' 
Acad. Press, 1996 

\bibitem[HoMo98]{HoMo98} Hofmann, K.\ H., and S.\ A.\ Morris, ``The Structure of
Compact Groups,'' Studies in Math., de Gruyter, Berlin, 1998. 

\bibitem[KP84]{KP84} Kac, V. G., and D. H. Peterson, {\it Unitary structure in 
representations of infinite dimensional groups and a convexity theorem}, 
Invent. Math. {\bf 76} (1984), 1--14 

\bibitem[KR87]{KR87} Kac, V. G., and A. K. Raina, ``Highest weight representations of 
infinite dimensional Lie algebras,'' Advanced Series in Math. Physics, 
World Scientific, Singapore, 1987

\bibitem[Ka83]{Ka83} Kaup, W., {\it \"Uber
die Klassifikation der symmetrischen
her\-mi\-te\-schen Mannigfaltigkeiten unendlicher Dimension I, II},
Math. Annalen {\bf 257}(1981), 463--486; {\bf 262}(1983), 57--75

\bibitem[Ka97]{Ka97} Kaup, W., {\it On real Cartan factors}, 
manuscripta math. {\bf 92} (1997), 191--222 

\bibitem[Ki87]{Ki87} Kirillov, A., {\it K\"ahler structure on the 
$K$-orbits of the group 
of diffeomorphisms of a circle}, Funct. Anal. Appl. {\bf 21:2} 
(1987), 122--125 

\bibitem[Ki98]{Ki98} ---, {\it Geometric approach to 
discrete series of unirreps for Vir}, 
J. Math. Pures Appl. {\bf 77:8}  (1998),  735--746

\bibitem[KY87]{KY87} 
Kirillov, A.~A., and D.~V.~Yuriev, {\it K\"ahler geometry of
the infinite-dimensional homogeneous space $M =
\Diff_+(\bS^1)/\Rot(\bS^1)$}, Funct.\ Anal.\ and Appl. {\bf 21:4}
(1987), 284--294 

\bibitem[KY88]{KY88} ---, {\it Representations of the Virasoro algebra 
by the orbit method}, J. Geom. Phys. {\bf 5:3} (1988), 351--363 

\bibitem[KS88]{KS88} Kostant, B., and S. Sternberg, {\it 
The Schwarzian derivative and the conformal geometry of the 
Lorentz hyperboloid}, in ``Quantum Theories and Geometry,'' 
M. Cahen and M. Flato eds., Kluwer, 1988; 113--1225 

\bibitem[Lm94]{Lm94} Langmann, E., {\it Cocycles for boson and fermion 
Bogoliubov transformations},  J. Math. Phys. {\bf 35:1} (1994), 96--112 

\bibitem[Le95]{Le95} Lempert, L., {\it The Virasoro group as a complex manifold}, 
Math. Res. Letters {\bf 2} (1995), 479--495 

\bibitem[MdR06]{MdR06} Mack, G., and M. de Riese, 
{\it Simple Space-time Symmetries: Generalized Conformal Field Theory}, 
arXiv:hep-th/0410277v2, 25. Jan.\ 2006

\bibitem[Mag92]{Mag92} Magyar, M., ``Continuous Linear Representations,'' 
North-Holland, Math.\ Studies {\bf 168}, 1992

\bibitem[MR99]{MR99} Marsden, J. E., and T. Ratiu, ``Introduction to 
Mechanics and Symmetry,'' 2n Ed., Texts in Appl. Math. {\bf 17}, 
Springer, 1999 

\bibitem[MR85] {MR85} Medina, A., and P. Revoy, {\it Alg\`ebres de Lie et produit scalaire 
invariant}, Ann. scient. \'Ec. Norm. Sup. $4^e$ s\'erie {\bf 18}(1985), 
533-561

\bibitem[Mick87]{Mick87} Mickelsson, J., {\it Kac--Moody groups, topology of the Dirac determinant bundle, and fermionization}, Commun. Math. Phys. {\bf 110} (1987), 173--183 

\bibitem[Mick89]{Mick89} ---, ``Current algebras and groups,'' Plenum Press,
New York, 1989 

\bibitem[Mil84]{Mil84} Milnor, J., {\it Remarks on 
infinite dimensional Lie groups}, in DeWitt, B., Stora, R. (eds), 
``Relativit\'{e}, groupes et topologie II (Les Houches, 1983), North Holland, Amsterdam, 1984; 1007--1057 

\bibitem[Ne94]{Ne94} Neeb, K.-H., {\it The classification 
of Lie algebras with invariant cones}, 
J.\ Lie Theory {\bf 4:2} (1994), 1--47 

\bibitem[Ne98]{Ne98} ---, {\it Holomorphic highest weight representations
of infinite dimensional complex classical groups}, 
J.\ Reine Angew. Math.\ {\bf 497} (1998), 171--222  

\bibitem[Ne00]{Ne00} ---, ``Holomorphy and Convexity in Lie Theory,'' 
Expositions in Mathematics {\bf 28}, de Gruyter Verlag, Berlin, 1999 

\bibitem[Ne01a]{Ne01a} ---, {\it Compressions of infinite dimensional bounded symmetric
domains}, Semigroup Forum {\bf 63:1} (2001), 71--105

\bibitem[Ne01b]{Ne01b} Neeb, K.-H., {\it Representations of infinite dimensional
  groups}, pp.~131--178; in ``Infinite Dimensional K\"ahler Manifolds,'' 
Eds. A. Huckleberry, T. Wurzbacher, DMV-Seminar {\bf 31}, 
Birkh\"auser Verlag, 2001 

\bibitem[Ne02a]{Ne02a} ---, {\it Classical Hilbert--Lie groups, 
their extensions and their
homotopy groups}, in ``Geometry and Analysis on Finite-
and Infinite-Dimensional Lie Groups,'' 
 Eds. A.~Strasburger et al., 
Banach Center Publications {\bf 55}, Warszawa 2002; 87--151 

\bibitem[Ne02b]{Ne02b} ---, 
{\it Central extensions of infinite-dimensional
Lie groups}, Annales de l'Inst. Fourier {\bf 52} (2002), 1365--1442 

\bibitem[Ne04]{Ne04} K.-H. Neeb, {\rm Infinite-dimensional Lie groups and their
represen\-ta\-tions}, in ``Lie Theory: Lie Algebras and Representations'', 
Progress in Math. {\bf 228}, 
Ed. J.~P.~Anker, B. \O{}rsted, Birkh\"auser Verlag, 2004; 213--328 
   
\bibitem[Ne06]{Ne06} ---, {\it Towards a Lie theory of locally convex 
groups}, Jap. J. Math. 3rd ser. {\bf 1:2} (2006), 291--468 

\bibitem[Ne08]{Ne08} ---, {\it A complex semigroup approach to group 
algebras of infinite dimensional Lie groups}, Semigroup Forum 
{\bf 77} (2008), 5--35

\bibitem[Ne09a]{Ne09a} ---, 
{\it Semibounded unitary representations of infinite dimensional Lie groups}, to appear in 
``Infinite Dimensional Harmonic Analysis IV'', 
Eds. J.\ Hilgert et al, World Scientific, 2009; 209--222

\bibitem[Ne09b]{Ne09b} ---, {\it Unitary highest weight modules of 
locally affine Lie algebras}, in ``Proceedings of 
the Workshop on Quantum Affine Algebras, Extended Affine Lie Algebras
and Applications (Banff, 2008)'', Eds. Y. Gao et al, Contemporary Math., 
Amer. Math. Soc., to appear 

\bibitem[Ne09c]{Ne09c} Neeb, K.-H., {\it Invariant convex cones 
in infinite  dimensional Lie algebras}, in preparation 

\bibitem[Ne10]{Ne10} Neeb, K.-H., {\it 
On differentiable vectors for representations of infinite dimensional 
Lie groups}, in preparation 

\bibitem[NO98]{NO98} Neeb, K.-H., B. \O{}rsted, {\it Unitary highest 
weight representations in 
Hilbert spaces of holomorphic functions on infinite 
dimensional domains}, J. Funct.\ Anal. {\bf 156} (1998), 263--300

\bibitem[Ner90]{Ner90} Neretin, Y., {\it Holomorphic extension 
of representations of the group of diffeomorphisms of the 
circle}, Math. USSR Sbornik {\bf 67:1} (1990), 75--97 



\bibitem[Neu99]{Neu99} Neumann, A., 
{\it The classification of symplectic structures of
convex type}, Geom. Dedicata {\bf 79:3} (2000), 299--320 

\bibitem[Neu02]{Neu02} Neumann, A., 
{\it An infinite dimensional version of the
Kostant Convexity Theorem}, J.~Funct.\ Anal. {\bf 189}
(2002), 80--131 

\bibitem[Ols82]{Ols82} Ol'shanski\u\i, G.I., 
{\it Invariant cones in Lie algebras,
Lie semigroups, and the holomorphic discrete series}, Funct. Anal.
Appl. {\bf 15} (1982), 275--285 

\bibitem[Ot95]{Ot95} Ottesen, J. T., 
``Infinite Dimensional Groups and Algebras in Quantum 
Physics,'' Springer-Verlag, Lecture Notes in 
Physics {\bf m 27}, 1995 

\bibitem[Ov01]{Ov01}  Ovsienko, V. Yu., {\it Coadjoint representation 
of Virasoro-type Lie algebras and differential operators on 
tensor-densities}, in ``Infinite dimensional K\"ahler manifolds'' 
(Oberwolfach, 1995),  DMV Sem. {\bf 31}, Birkhäuser, Basel, 2001; 231--255 

\bibitem[Pe89]{Ped89}
Pedersen, G. K., ``{\rm C}$^*$--Algebras and their Automorphism Groups,''
  London, Academic Press,  1989

\bibitem[PoSt70]{PoSt70} Powers, R. T., and E. St\o{}rmer, {\it Free states 
of the Canonical Anticommutation Relations}, Comm. Math. Phys. 
{\bf 16} (1970), 1--33 

\bibitem[PS86]{PS86} Pressley, A., and G. Segal, ``Loop Groups," Oxford University Press, 
Oxford, 1986

\bibitem[Re69]{Re69} Reed., M. C., {\it A G\aa{}rding domain 
for quantum fields}, Comm. Math. Phys. {\bf 14} (1969), 336--346 

\bibitem[RS75]{RS75} Reed, S., and B. Simon, 
``Methods of Modern Mathematical Physics I: Functional Analysis,'' 
Academic Press, New York, 1973

\bibitem[Ru77]{Ru77} Ruijsenaars, S. N. M., {\it 
On Bogoliubov transformations for systems of relativistic charged 
particles}, J. Math. Physics {\bf 18:3} (1977), 517--526 

\bibitem[Sa91]{Sa91} Samoilenko, Y. S., ``Spectral Theory of Families 
of Self-Adjoint Operators,'' Kluwer Acad. Publ., 1991 

\bibitem[Sch97]{Sch97} Schottenloher, M., ``A Mathematical Introduction 
to Conformal Field Theory,'' Lecture Notes in Physics {\bf m 43}, 
Springer, 1997 

\bibitem[Sw73]{Sw73} Schwartz, L., ``Th\'eorie des distributions,'' 
Hermann, Paris, 1973 

\bibitem[SeG81]{SeG81} Segal, G., {\it Unitary representations of some
infinite-dimensional groups}, Comm.\ Math.\ Phys. {\bf 80} (1981), 301--342 

\bibitem[Se57]{Se57} Segal, I.E., 
{\it The structure of a class of representations of the unitary 
group on a Hilbert space}, Proc. Amer. Math. Soc. {\bf 81} (1957), 
197--203 

\bibitem[Se58]{Se58} ---, {\it Distributions in Hilbert spaces and canonical
systems of operators}, Trans.\ Amer.\ Math.\ Soc.\ 
{\bf 88} (1958), 12--41 

\bibitem[Se59]{Se59} ---, {\it Foundations 
of the theory of dynamical systems of
infinitely many degrees of freedom. I}, Mat. Fys. Medd. Danske
Vid. Selsk. {\bf 31:12} (1959), 1--39

\bibitem[Se67]{Se67} ---, {\it Positive-energy particle models with 
mass splitting}, Proc. Nat. Acad. Sci. U.S.A. {\bf 57} (1967), 194--197 

\bibitem[Se76]{Se76} ---, ``Mathematical Cosmology and Extragalactical Astronomy,'' 
Academic Press, New York, San Francisco, London, 1976 

\bibitem[Se78]{Se78} ---, {\it The Complex-wave Representation of the Free
Boson Field}, in ``Topics in Funct. Anal.,'' Adv.\ in Math. Suppl. 
Studies {\bf 3} (1978), 321--343 

\bibitem[SJOSV78]{SJOSV78} Segal, I. E., H. P. Jakobsen, B. \O{}rsted, 
B. Speh, and M. Vergne, {\it Symmetry and causality 
properties of physical fields}, Proc. Nat. Acad. Sci. U.S.A. 
{\bf 75:4} (1978), 1609--1611 

\bibitem[Sh62]{Sh62}  Shale, D., {\it Linear symmetries of free boson fields}, 
Trans.\ Amer. Math. Soc. {\bf 103} (1962), 146--169

\bibitem[ShSt65]{ShSt65} Shale, D., and W. F. Stinespring, {\it 
Spinor representations of infinite orthogonal groups}, 
J. Math. Mech.  14  1965 315--322

\bibitem[SW98]{SW98} Spera, M., and T. Wurzbacher, {\it Determinants, 
Pfaffians and quasi-free representations of the CAR algebra}, 
Rev. Math. Phys.  {\bf 10:5}  (1998), 705--721

\bibitem[Th09]{Th09} Thill, M., {\it Representations of hermitian 
commutative $*$-algebras by unbounded operators}, arXiv:math.OA.0908.3267, 
22. Aug. 2009 

\bibitem[TL99a]{TL99a} Toledano Laredo, V., {\it Positive energy 
representations of the loop groups of non-simply connected Lie groups}, 
Comm. Math. Phys. {\bf 207:2} (1999),  307--339

\bibitem[TL99b]{TL99b} ---,  {\it Integrating unitary representations
of infinite dimensional Lie groups}, Journal of Funct.\ Anal. {\bf
161} (1999), 478--508 

\bibitem[Up85]{Up85} Upmeier, H., ``Symmetric Banach Manifolds and Jordan 
$C^*$-algebras,'' North Holland Mathematics Studies, 1985 

\bibitem[Ve77]{Ve77} Vergne, M., {\it Groupe symplectique et seconde
quantification}, C.\ R.\ Acad.\ Sc.\ Paris {\bf 285} (1977), 191--194 

\bibitem[Ve90]{Ve90} Vershik, A. M., {\it Metagonal 
groups of finite and infinite dimension}, 
in ``Representations of Lie Groups and Related Topics,'' A.~ M.~Vershik
and D.~P.~Zhelobenko eds., Gordon and Breach Science Publishers, 1990;
1--37 

\bibitem[Vin80]{Vin80} Vinberg, E.B., 
{\it Invariant cones and orderings in Lie 
groups}, Funct. Anal. Appl. {\bf 14} (1980), 1--13 

\end{thebibliography}
\end{document}